\documentclass{article}
\usepackage[a4paper,top=2cm,bottom=2cm,left=2cm,right=2cm]{geometry}

\usepackage{amsmath}
\usepackage{amssymb}
\usepackage{amsthm}
\usepackage{mathrsfs}
\usepackage[scr=rsfs,cal=boondox]{mathalfa} 
\usepackage{tikz-cd}
\usetikzlibrary{mindmap,trees}
\usepackage{stmaryrd}
\usepackage{hyperref}
\usepackage{authblk}
\usetikzlibrary{positioning}
\usepackage{braids}
\usepackage{t-angles}
\usepackage[all]{xy}
\def\pulb{\ar@{}[dr]|(0.2){\mbox{\Large{$\lrcorner$}}}}

\usepackage{comment}

\usepackage{mathtools}

\usepackage[toc,page]{appendix}

\usepackage[bbgreekl]{mathbbol}
\usepackage{bbm}

\newtheorem{proposition}{Proposition}[section]
\newtheorem{definition}[proposition]{Definition}
\newtheorem{lemma}[proposition]{Lemma}
\newtheorem{theorem}[proposition]{Theorem}
\newtheorem{remark}[proposition]{Remark}
\newtheorem{example}[proposition]{Example}
\newtheorem{corollary}[proposition]{Corollary}

\newcommand{\ten}{\otimes}
\newcommand{\id}{\mathrm{id}}
\newcommand{\op}{\mathrm{op}}
\newcommand{\cC}{\mathcal{C}}

\allowdisplaybreaks

\title{{\bf Quantization of infinitesimal braidings and pre-Cartier quasi-bialgebras}}

\author[a]{Chiara Esposito}
\author[b]{Andrea Rivezzi}
\author[c]{Jonas Schnitzer}
\author[b]{Thomas Weber}

\affil[a]{{\sl Dipartimento di Matematica, Università degli Studi di Salerno}

{\sl Via Giovanni Paolo II, 132, 84084 Fisciano, Italy}

{~}}

\affil[b]{{\sl
{Mathematical Institute of Charles University}}

{\sl Sokolovsk\'a 49/83, 186 75 Prague 8 , Czech Republic}

{~}}

\affil[c]{{\sl
{Dipartimento di Matematica ``Felice Casorati", Università degli Studi di Pavia}}
	
{\sl Via Ferrata 5, 27100 Pavia, Italy}}

\date{}

\begin{document}
	
\maketitle
	
\begin{abstract}
\noindent
In this paper we extend Cartier's deformation theorem of braided monoidal categories admitting an infinitesimal braiding to the non-symmetric case.  The algebraic counterpart of these categories is the notion of a pre-Cartier quasi-bialgebra, which extends the well-known notion of quasitriangular quasi-bialgebra given by Drinfeld. Our result implies that one can quantize the infinitesimal $\mathcal{R}$-matrix of any Cartier quasi-bialgebra.
We further discuss the emerging concepts of infinitesimal quantum Yang-Baxter equation and Cartier ring, the latter containing braid groups with additional generators that correspond to infinitesimal braidings. Explicit deformations of the representation categories of the gauge deformed quasitriangular quasi-bialgebras $E(n)$ are provided.

\end{abstract}
\tableofcontents

\section{Introduction}

This article focuses on the interplay of braided monoidal categories
\cite{JS93} and quasitriangular (quasi-)bialgebras \cite{DrQG}. We
discuss their first-order deformations, as proposed in \cite{ABSW},
and address the resulting deformation problem, alongside providing
solutions of the infinitesimal braid relations and the infinitesimal
Yang-Baxter equation. There is a well-known correspondence between
categorical data, given by braided monoidal categories
$(\mathcal{C},\otimes,I,\sigma)$, and algebraic data, given
by quasitriangular bialgebras $(H,\Delta,\varepsilon,\mathcal{R})$:
the representation category ${}_H\mathcal{M}$ of a bialgebra $H$ is
monoidal, with the tensor product representation and the trivial
representation determined by the coproduct $\Delta\colon H\to H\otimes
H$ and counit $\varepsilon\colon H\to\Bbbk$, respectively. Conversely,
if a tensor category $(\mathcal{C},\otimes,I)$ admits a
representable strong monoidal functor to vector spaces, one is able to
obtain a bialgebra $H$ such that the functor factors through the
representation category of $H$ via Tannaka--Krein reconstruction
\cite{SaaRiv}.  In these cases, the representation category of $H$ is
braided if and only if $H$ is endowed with a quasitriangular structure
$\mathcal{R}\in H\otimes H$, so that the braiding $\sigma$ is
determined by the action of $\mathcal{R}$, and it follows that the
bialgebra reconstructed from a braided monoidal category is
quasitriangular, as well (see \cite[Proposition
  9.4.2]{MajidFoundation} for the latter statement). 
As $\mathcal{R}$-matrices provide solution of the quantum Yang-Baxter equation \cite{Baxter}, quasitriangular bialgebras admit various applications, e.g. in the construction of topological invariants of $3$-manifolds \cite{ReTu} and the quantum inverse scattering method \cite{FaTa}. Explicit constructions of $\mathcal{R}$-matrices involve Drinfeld twists \cite{DrKZ}, such as the famous Moyal--Weyl twist. There is active research on the applications of Drinfeld twists and their induced $\mathcal{R}$-matrices to noncommutative geometry and deformation quantization, see e.g. \cite{PaoloWess,ChJoSt,Thomas}

\vspace{0.3cm}

\noindent
Thus, one expects
that additional categorical structures of braided monoidal categories
should correspond to additional algebraic structures on the level of
quasitriangular bialgebras.
One such instance is an \emph{infinitesimal braiding} $t$ of a braided
monoidal category $(\mathcal{C},\otimes,I,\sigma)$. The
former is a natural transformation $t\colon\otimes\Rightarrow\otimes$
of the tensor product, satisfying linearized versions of the hexagon
relations (see Definition \ref{def:preCart}). Identifying
$\mathcal{C}$ with the representation category of a quasitriangular
bialgebra $H$, as above, the infinitesimal braiding corresponds to an
\emph{infinitesimal $\mathcal{R}$-matrix} $\chi\in H\otimes H$, which
can be interpreted as a first-order deformation
\begin{equation}\label{eq:Rtilde}
	\tilde{\mathcal{R}}=\mathcal{R}(1\otimes 1+\hbar\chi+\mathcal{O}(\hbar^2))
\end{equation}
of the $\mathcal{R}$-matrix $\mathcal{R}$ with a formal parameter
$\hbar$. More precisely, $\chi$ commutes with the coproduct and
satisfies linearized versions of the hexagon relations (see Definition
\ref{def:preCartquasibialg}). The notions of infinitesimal braiding of
a braided monoidal category and infinitesimal $\mathcal{R}$-matrix of
a quasitriangular bialgebra have been introduced in \cite{ABSW} under
the name \emph{pre-Cartier} category and bialgebra, respectively. The
attribution to P. Cartier stems from the foundational work \cite{Car},
where infinitesimal braidings of symmetric categories have been
introduced. Most remarkably, it has been shown in \cite{Car} that
symmetric monoidal categories
$(\mathcal{C},\otimes,I,\sigma)$ endowed with an
infinitesimal braiding $t$ such that
\begin{equation}\label{eq:sigmat}
	\sigma_{M,N}\circ t_{M,N}=t_{N,M}\circ\sigma_{M,N}
\end{equation}
on all objects $M,N$, admit a \emph{deformation} as braided monoidal
category in the following sense: given $\Psi$, a Drinfeld associator
\cite{DR89}, there is a deformed associativity constraint on the
$\hbar$-adic completed category
$\tilde{\mathcal{C}}=(\mathcal{C},\tilde{\otimes},I)$, such
that
\begin{equation}\label{eq:sigmaexpt}
	\tilde{\sigma}_{M,N}:=\sigma_{M,N}\circ\exp\bigg(\frac{\hbar}{2}t_{M,N}\bigg)
\end{equation}
structures $\tilde{\mathcal{C}}$ as a braided monoidal category. As
pointed out for example in \cite[Chapter XX.6]{Kassel95}, this
quantization procedure recovers the $\hbar$-adic version of the
Drinfeld--Jimbo quantum groups $U_q(\mathfrak{g})$ (see
\cite{DrHopfAlg,Jimbo}) of a complex semisimple Lie algebra
$\mathfrak{g}$, where the infinitesimal braiding is determined by the
Casimir of the universal enveloping algebra $U(\mathfrak{g})$. The
same machinery has been utilized in \cite{EtiKaz,Sev} for
the universal quantization of Lie bialgebras and in \cite{PulSev} for the universal quantization of Poisson Hopf algebras. Furthermore,
infinitesimally braided categories and their deformation recently
attracted attention in the contexts of homology \cite{FGS}, algebraic
topology \cite{KKMP} and higher category theory \cite{CirMart,KeLaSc,Kemp}. Infinitesimal braidings for symmetric categories further appear in the study of curved Lie bialgebras \cite{HeckVen}.

\noindent
In the course of \cite{ABSW}, new examples of infinitesimally braided
monoidal categories and pre-Cartier bialgebras appeared, which do not
fit the framework of \cite{Car} since they are either not based on
symmetric categories, or not complying with \eqref{eq:sigmat}, or
neither of both. For example, the representation category of
Sweedler's $4$-dimensional Hopf algebra is symmetric, but with
infinitesimal braidings not satisfying \eqref{eq:sigmat}. More general
examples in all dimensions, satisfying some or none of the above
assumptions, appeared in \cite{BRS24}.

\vspace{0.3cm}
\noindent
The natural question arises if a pre-Cartier category admits a
deformation whose braiding is given by \eqref{eq:sigmaexpt}, i.e. if a
generalization of Cartier's theorem exists. Or, in the sense of the
previously mentioned duality with bialgebras, if a pre-Cartier
bialgebra $(H,\Delta,\varepsilon,\mathcal{R},\chi)$ admits a
deformation $\tilde{H}=H[[\hbar]]$ with deformed coproduct $\tilde{\Delta}$
such that \eqref{eq:Rtilde} is a quasitriangular structure on
$\tilde{H}$. In this paper we give a positive answer to all the above
questions in case of braided monoidal categories with infinitesimal
braidings satisfying \eqref{eq:sigmat} or such that the commutator
$[t\otimes\mathrm{id},\mathrm{id}\otimes t]$ vanishes. This addresses
Question 2.10 raised in \cite{ABSW}. In particular, it turns out that
the associativity constraint of the category has to be deformed in
terms of the infinitesimal braiding via a chosen Drinfeld
associator. On the level of algebra this means that the deformation
constitutes a quasitriangular \emph{quasi}-bialgebra, i.e. there is a
re-associator, measuring the failure of coassociativity.  This tells us
that the deformation problem is best understood in the context of
quasi-bialgebras and thus we formulate the notion of pre-Cartier
structure for quasi-bialgebras. On the categorical level this means
that one obtains non-trivial associativity constraints after
strictification and we formulate the notion of pre-Cartier category in
this setup accordingly. Apart from the aforementioned deformation
results, and in fact as a preparation for them, we prove that every
pre-Cartier category admits two solutions of the infinitesimal braid
relations (collapsing into one solution in case of a trivial
re-associator) and the emerging infinitesimal version of the quantum
Yang-Baxter equation is discussed. We further introduce Cartier rings, which are represented by
pre-Cartier quasi-bialgebras. They have two sets of generators,
corresponding to the actions of the $\mathcal{R}$-matrix and
infinitesimal $\mathcal{R}$-matrix, respectively.

\vspace{0.3cm}

\noindent
\emph{Application and further research.} In Section \ref{sec:3.5} we briefly introduce Cartier rings and their representations. In the same way a module over a quasitriangular bialgebra gives rise to representations of the braid groups, a module over a pre-Cartier quasi-bialgebra induces representations of the Cartier rings. We believe that this interaction can be as fruitful as it is for the braid groups and quasitriangular bialgebras, and we plan to investigate this in the future.

\noindent
A very obvious extension of our work is given by general deformations of
pre-Cartier categories, i.e. without the assumptions \eqref{eq:sigmat}
or $[t\otimes\mathrm{id},\mathrm{id}\otimes t]=0$. In doing so, we
believe that the exponential ansatz \eqref{eq:sigmaexpt} is too
restrictive. Nevertheless, we are planning to consider this general
situation in future work.

\noindent
Another possible application we want to mention is the dequantization
procedure from \cite{EtiKazII}, since there the authors consider the
category of Yetter-Drinfeld modules over a quantized enveloping algebra, which is a
truly non-symmetric braided monoidal category. The dequantization procedure
is now based on finding a curve of compatible Drinfeld associators and
braidings connecting the given braiding to a symmetric one using the Grothendieck-Teichmüller 
group. The formal derivative of this curve of braidings should give rise to an infinitesimal braiding
in a non-symmetric monoidal category. The obvious question one can ask is now, what the curve of the 
monoidal structure has to do with the application of the Drinfeld associator to this infinitesimal
braiding. We plan to explore this relation in light of our main theorem in this context 
in a future work.  

\vspace{0.3cm}

\noindent
\emph{Structure of the paper.} We first recall the concept of braided
monoidal category and its corresponding algebraic counterpart given by
quasitriangular quasi-bialgebras in Sections \ref{sec:2.1} and
\ref{sec:2.2}, respectively. We further recall the definition and
basic properties of Drinfeld associators, together with the related
concept of infinitesimal braid relation \cite{Kohno} in Section
\ref{sec:2.3}. Adding infinitesimal braidings and infinitesimal
$\mathcal{R}$-matrices leads to the notions of pre-Cartier category
(Section \ref{sec:pC}) and pre-Cartier quasi-bialgebra (Section
\ref{sec:pCbialg}), which generalizes the account \cite{ABSW}. Here,
we focus on examples and fundamental constructions, such as twist
deformations. In Section~\ref{sec:3.3} we continue to prove that every
pre-Cartier quasi-bialgebra admits \emph{two} solutions of the
infinitesimal braid relations. This phenomenon is purely due to the
non-trivial re-associator, since in the coassociative setup both
solutions coincide. An emerging infinitesimal version of the quantum
Yang-Baxter equation is discussed in Section \ref{sec:3.4}, while the
notion of Cartier ring and its corresponding representations are
studied in Section \ref{sec:3.5}. In Section~\ref{sec:4.1} we prove
that the categorical counterparts of the previous solutions of the
infinitesimal braid relations exist for arbitrary pre-Cartier
categories. We emphasize that this result functions as the main lemma
to obtain the subsequent deformation results. Namely, we prove in
Section \ref{sec:CartierCase} that a braided monoidal category with
infinitesimal braiding satisfying \eqref{eq:sigmat} admits a
deformation in the previously discussed sense. As shown in
Section~\ref{sec:preCartierCase}, condition \eqref{eq:sigmat} can be
disregarded if one requires the vanishing of a certain commutator of
$t$ instead. Note that this results in an undeformed associativity
constraint, while the braiding is deformed non-trivially via
\eqref{eq:sigmaexpt}. Examples of deformation of the representation
categories of the $E(n)$ algebras and their quasi-bialgebra versions
obtained from a gauge twisting are discussed in detail. Finally, in
Section \ref{sec:4.4}, the algebraic counterparts of the deformation
result in Sections \ref{sec:CartierCase} and \ref{sec:preCartierCase} are given.

\subsubsection*{Conventions and notations}

We fix a field $\Bbbk$ and we assume the characteristic of $\Bbbk$ to be zero, if not specified otherwise. All algebras are understood as $\Bbbk$-algebras
and all maps are assumed to be $\Bbbk$-linear. For the cardinality of sets we write $\#$. We denote by ${}_H\mathcal{M}$ the category of left modules of an algebra $H$. The braiding of the category of vector spaces (i.e. the usual tensor flip) is denoted by $\tau$. Throughout this paper $\hbar$ is a formal parameter. Moreover, we use Sweedler's short notation as well as Einstein sum convention over repeated indices.

\section*{Acknowledgments}
We are grateful to Matteo Misurati and Fabio Renda for clarifying some aspects of the quasi-bialgebra theory.
C.E. thanks her Aikido master for his teachings, which help her rediscover the calm and focus essential to research, and her dojo for the friendship and inspiring atmosphere that have supported and enriched this work. Additionally, we are grateful to the referee for valuable suggestions and comments.
C.E. has been supported from Italian Ministerial grant PRIN 2022
“Cluster algebras and Poisson Lie groups'', n. 20223FEA2E - CUP
D53C24003330006. 
C.E. and
A.R. were supported by the National Group for Algebraic and Geometric
Structures, and their Applications (GNSAGA -- INdAM).  A.R.  is
supported by GA\v{C}R/NCN grant Quantum Geometric Representation
Theory and Noncommutative Fibrations 24-11728K. T.W. is supported by
the GAČR PIF 24-11324I.  This publication is based upon work from
\emph{COST Action CaLISTA CA21109} supported by COST (European
Cooperation in Science and Technology), www.cost.eu.
%

\section{Preliminaries}

Here we recall some basic notions regarding braided monoidal categories, needed throughout
this paper, as well as the algebraic counterpart consisting of quasitriangular quasi-bialgebras.
%

\subsection{Braided monoidal categories}\label{sec:2.1}

Mainly to fix notation we recall the following categorical notions.
\begin{definition}[\cite{MacLan63}\cite{JS93}]
  A \textbf{braided monoidal category} is a datum $ (\mathcal{C},
  \otimes, I, a, \ell, r,\sigma)$, where $\mathcal{C}$ is a category,
  $\otimes : \mathcal{C} \times \mathcal{C} \to \mathcal{C}$ is a
  functor (called the tensor product), $I$ is an object in
  $\mathcal{C}$ (called the tensor unit), and
  \begin{equation*}
    \begin{split}
      a&: \ten \circ (\ten \times \id) \Rightarrow \ten \circ (\id \times \ten) \\
      \ell&: \ten \circ (I \times \id) \Rightarrow \id \\
      r &: \ten \circ (\id \times I) \Rightarrow \id \\
      \sigma &: \ten \Rightarrow \ten^\op
    \end{split}
  \end{equation*}
  are natural isomorphisms --called respectively the associativity
  constraint, the left unit constraint, the right unit constraint, and
  the braiding-- making commutative the following diagrams for all
  objects $X,Y,Z,W$ in $\mathcal{C}$:
  \begin{enumerate}
  \item[i.)] (pentagon axiom)
    \begin{equation}
      \label{eq:pentagon-axiom}
      \begin{tikzcd}
        \big(X \ten( Y \ten Z)\big) \ten W \arrow[dd, "{a_{X,Y \ten Z,W}}"']     &  & \big((X\ten Y) \ten Z\big) \ten W \arrow[d, "{a_{X \ten Y,Z,W}}"] \arrow[ll, "{a_{X,Y,Z} \ten \id_W}"'] \\
        &  & (X \ten Y) \ten (Z \ten W) \arrow[d, "{a_{X,Y,Z \ten W}}"]                                      \\
        X \ten \big((Y \ten Z) \ten W\big) \arrow[rr, "{\id_X \ten a_{Y,Z,W}}"'] &  & X \ten \big(Y \ten (Z \ten W)\big)                                                                     
      \end{tikzcd}
    \end{equation} 
  \item[ii.)] (triangle axiom)
    \begin{equation}
      \label{eq:triangle-axiom}
      \begin{tikzcd}
        (X \ten I) \ten Y \arrow[rd, "r_X \ten \id_Y"'] \arrow[rr, "{a_{X,I,Y}}"] &          & X \ten(I \ten Y) \arrow[ld, "\id_X \ten \ell_Y"] \\
        & X \ten Y &                                                 
      \end{tikzcd}
    \end{equation} 
  \item[iii.)] (hexagon axioms)
    \begin{equation}
      \label{eq:hexagon-axiom-one}
      \begin{tikzcd}
        & X \ten (Y \ten Z) \arrow[r, "\sigma_{X,Y \ten Z}"]  & (Y \ten Z) \ten X \arrow[rd, "a_{Y,Z,X}"] &  \\
        (X \ten Y) \ten Z \arrow[ru,"a_{X,Y,Z}" ] \arrow[rd, " \sigma_{X,Y} \ten \id_Z"']&  &  & Y \ten(Z \ten X) \\
        & (Y \ten X) \ten Z \arrow[r, "a_{Y,X,Z}"'] & Y \ten (X \ten Z) \arrow[ru, "\id_Y \ten \sigma_{X,Z}"'] & 
      \end{tikzcd}
    \end{equation}
    \begin{equation}
      \label{eq:hexagon-axiom-two}
      \begin{tikzcd}
        & (X \ten Y) \ten Z \arrow[r, "\sigma_{X \ten Y, Z}"]  & Z \ten (X  \ten Y) \arrow[rd, "a^{-1}_{Z,X,Y}"] &  \\
        X \ten (Y \ten Z) \arrow[ru,"a^{-1}_{X,Y,Z}" ] \arrow[rd, "\id_X \ten \sigma_{Y,Z}"']&  &  & (Z \ten X) \ten Y \\
        & X \ten (Z \ten Y) \arrow[r, "a^{-1}_{X,Z,Y}"'] & (X \ten Z)\ten Y \arrow[ru, "\sigma_{X,Z} \ten \id_Y"'] & 
      \end{tikzcd}
    \end{equation}
  \end{enumerate}
  If moreover $\sigma^{-1}_{X,Y} =\sigma_{Y,X} $ for all objects $X,Y$
  of $\mathcal{C}$, we say that $\mathcal{C}$ is a \textbf{symmetric braided monoidal category}.
\end{definition}
\noindent
We now introduce braided monoidal functors, which is the appropriate
notion of functor between braided monoidal categories.
\begin{definition}
  \label{definition-monoidal-functor}
  Let $(\mathcal{C}, \ten, I, a,\ell,r,\sigma)$ and $(\mathcal{C}', \ten', I', a',\ell',r',\sigma')$ be two braided monoidal categories.
  \begin{enumerate}
  \item[i.)] A \textbf{braided monoidal functor}\footnote{In the literature there is also an alternative convention, where braided monoidal functors are referred to as braided lax monoidal functors, while the term braided monoidal functor is sometimes used to indicate strong monoidal functors.} from $\cC$ to $\cC'$ is
    a triple $(F,F^2,F^0)$, where $F:\cC \rightarrow \cC'$ is a functor,
    $F^0: I' \to F(I)$ is a morphism, and $F^2$ is a natural
    transformation $F^2 : \ten \circ (F \times F) \Rightarrow F \circ
    \ten$ which is compatible with the associativity constraints, unit
    constraints, and the braidings, i.e. the following diagrams commute
    for all objects $X,Y,Z$ in $\cC$:
    \begin{equation}
      \label{eq:monoidal-functor-one}
      \begin{tikzcd}
        \big(F(X) \ten' F(Y)\big) \ten' F(Z) \arrow[rr, "{{a'}_{F(X),F(Y),F(Z)}}"] \arrow[d, "{F^2(X,Y) \ten' \id_{F(Z)}}"'] &  & F(X) \ten' \big(F(Y) \ten' F(Z)\big)  \arrow[d, "{\id_{F(X)} \ten'F^2(Y,Z) }"] \\
        F(X  \ten Y) \ten' F(Z)  \arrow[d, "{F^2(X \ten Y,Z)}"']                                                     &  & F(X) \ten' F(Y \ten Z) \arrow[d, "{F^2(X, Y \ten Z)}"]                 \\
        F\big((X \ten Y) \ten Z\big) \arrow[rr, "{F(a_{X,Y,Z})}"]                                                                  &  & F\big(X \ten (Y \ten Z)\big)                                                        
      \end{tikzcd}
    \end{equation}
    \begin{equation}
      \label{eq:monoidal-functor-two}
      \begin{tikzcd}
        I' \ten' F(X)  \arrow[r, "{\ell'}_{F(X)}"] \arrow[d, "F^0 \ten' \id_{F(X)}"'] & F(X)                                   &  & F(X) \ten' I'  \arrow[r, "r'_{F(X)}"] \arrow[d, "\id_{F(X)} \ten' F^0 "'] & F(X)                               \\
        F(I) \ten' F(X)  \arrow[r, "{F^2(I,X)}"]                                      & F(I \ten X)  \arrow[u, "F({\ell}_X)"'] &  & F(X) \ten' F(I)  \arrow[r, "{F^2(X,I)}"]                                 & F(X \ten I)  \arrow[u, "F(r_X) "']
      \end{tikzcd}
    \end{equation}
    \begin{equation}
      \label{eq:braided-monoidal-functor}
      \begin{tikzcd}
        F(X) \ten' F(Y) \arrow[r, "{F^2(X,Y)}"] \arrow[d, "{\sigma'_{F(X),F(Y)}}"'] & F(X \ten Y) \arrow[d, "{F(\sigma_{X,Y})}"] \\
        F(Y) \ten' F(X) \arrow[r, "{F^2(Y,X)}"]                                     & F(Y \ten X)                               
      \end{tikzcd}
    \end{equation}
    If both $F^2$ and $F^0$ are isomorphisms, we say that $(F,F^2,F^0)$ is \textbf{strong monoidal}.
  \item[ii.)]  A \textbf{natural braided monoidal transformation} $\eta
    : (F,F^2,F^0) \rightarrow (G,G^2,G^0)$ between braided monoidal
    functors from $\cC$ to $\cC'$ is a natural transformation $\eta : F
    \Rightarrow G$ such that the diagrams
    \begin{equation}
      \begin{tikzcd}
        F(X) \ten' F(Y)  \arrow[d, "\eta(X) \ten' \eta(Y) "'] \arrow[r, "{F^2(X,Y)}"] & F(X \ten Y) \arrow[d, "\eta(X \ten Y)"] &  &                             & I' \arrow[ld, "F^0"'] \arrow[rd, "G^0"] &      \\
        G(X) \ten' G(Y) \arrow[r, "{G^2(X,Y)}"]                                       & G(X \ten Y)                             &  & F(I) \arrow[rr, "\eta(I)"'] &                                                         & G(I)
      \end{tikzcd}
    \end{equation}
    commute for all objects $X,Y$ in $\cC$.
  \item[iii.)]  A \textbf{natural braided monoidal isomorphism} is a
    natural braided monoidal transformation that is also a natural
    isomorphism.
  \item[iv.)]  A \textbf{braided monoidal equivalence} is a braided
    monoidal functor $F: \cC \rightarrow \cC'$ such that there exist a
    braided monoidal functor $F': \cC' \rightarrow \cC$ and two natural
    braided monoidal isomorphisms $\eta: \id_{\cC'} \Rightarrow FF'$ and
    $\varepsilon: F'F \Rightarrow \id_{\cC}$. If there exists a braided
    monoidal equivalence between $\cC$ and $\cC'$ we say that $\cC,\cC'$
    are equivalent, and we write $\cC \cong \cC'$.
  \end{enumerate}
\end{definition}
\noindent
Next, we shall need the notion of strict monoidal category and the
MacLane coherence theorem.
\begin{definition}
  A braided monoidal category is said to be \textbf{strict} if the
  natural isomorphisms $a, \ell,r$ are all identity morphisms in
  $\mathcal{C}$. In particular, in a strict category one has $(X \ten Y)
  \ten Z = X \ten (Y \ten Z)$ and $I \ten X = X = X \ten I$ for all
  objects $X,Y,Z$ in $\cC$.
\end{definition}
\begin{remark}
  Note that in a strict braided monoidal category axioms
  \eqref{eq:pentagon-axiom}-\eqref{eq:triangle-axiom} hold
  automatically, while axioms
  \eqref{eq:hexagon-axiom-one}-\eqref{eq:hexagon-axiom-two} read
  \begin{align}
    \sigma_{X,Y \ten Z} &= (\id_Y \ten \sigma_{X,Z}) \circ (\sigma_{X,Y} \ten \id_Z) \label{eq:hexagon-strict-one} \\
    \sigma_{X \ten Y, Z} &= (\sigma_{X,Z} \ten \id_Y) \circ (\id_X \ten \sigma_{Y,Z})\label{eq:hexagon-strict-two}
  \end{align}
\end{remark}
\noindent
The following is one of the most important results of the theory of
monoidal categories, and it is known as the Mac Lane coherence
theorem, see \cite[\S 2.8]{EGNO} for more details:
\begin{theorem}
  Let $\cC$ be a braided monoidal category. Then there exists a strict
  braided monoidal category $\cC^{\mathrm{str}}$ such that $\cC \cong
  \cC^{\mathrm{str}}$.
\end{theorem}
\noindent
Mac Lane's coherence theorem allows to \emph{forget} about bracketings
when doing computations in a braided monoidal category. However, it is
crucial to stress out that the categories $\cC$ and
$\cC^{\mathrm{str}}$ have very different objects and tensor
products. However, P. Schauenburg \cite{Sch01} showed that in some
cases --such as the modules over a quasitriangular quasi-bialgebra--
one can build an equivalent strict category mantaining the same class
of objects. In the following, for the convenience of the reader, we
shall sometimes strictify braided monoidal categories.

\subsection{Quasitriangular quasi-bialgebras}\label{sec:2.2}

We further recall the notion of quasitriangular quasi-bialgebras
\cite{DrQH} and their deformations via gauge transformations,
following \cite[Chapter XV]{Kassel95}.

\noindent
Fix a field $\Bbbk$. All algebras are understood as $\Bbbk$-algebras and all maps are assumed to be $\Bbbk$-linear.

\begin{definition}
	\label{def:quasi-bialg}
	Let $H$ be an associative unital algebra and consider two algebra morphisms $\Delta\colon H\to H\otimes H$ and $\varepsilon\colon H\to\Bbbk$.
	\begin{enumerate}
		\item[i.)] We call $H$ a \textbf{quasi-bialgebra} if there exist invertible elements $\Phi\in H^{\otimes 3}$ and $\ell,r\in H$ such that
		\begin{align}
			(\mathrm{id}\otimes\Delta)\circ\Delta(\cdot)
			&=\Phi(\Delta\otimes\mathrm{id})\circ\Delta(\cdot)\Phi^{-1},\\
			(\varepsilon\otimes\mathrm{id})\circ\Delta(\cdot)
			&=\ell^{-1}(\cdot)\ell,\\
			(\mathrm{id}\otimes\varepsilon)\circ\Delta(\cdot)
			&=r^{-1}(\cdot)r
		\end{align}
		and
		\begin{align}
			(\mathrm{id}\otimes\mathrm{id}\otimes\Delta)(\Phi)(\Delta\otimes\mathrm{id}\otimes\mathrm{id})(\Phi)
			&=\Phi_{234}(\mathrm{id}\otimes\Delta\otimes\mathrm{id})(\Phi)\Phi_{123},\\
			(\mathrm{id}\otimes\varepsilon\otimes\mathrm{id})(\Phi)
			&=r\otimes\ell^{-1}.
		\end{align}
		The element $\Phi$ is referred to as the \textbf{re-associator} (see e.g. \cite[p. 109]{Bul} for the terminology). A morphism of quasi-bialgebras is a morphism of algebras $\alpha: H \to H'$ such that 
\[ (\alpha \ten \alpha) \circ \Delta = \Delta' \circ \alpha,\quad (\alpha \ten \alpha \ten \alpha)(\Phi) = \Phi', \quad \alpha(\ell) = \ell', \quad \alpha(r) = r'.\]		
		
		\item[ii.)] We call a quasi-bialgebra $(H,\Delta,\varepsilon,\Phi,\ell,r)$ \textbf{quasitriangular} if there is an invertible element $\mathcal{R}\in H\otimes H$, the \textbf{universal $\mathcal{R}$-matrix}, such that\footnote{Axioms \eqref{eq:qtqb2}, \eqref{eq:qtqb3} are slightly different in \cite{Kassel95}. We rather use the ones from \cite{Bul} since they guarantee the categorical correspondence.}
		\begin{align}
			\Delta^\mathrm{op}(\cdot)
			&=\mathcal{R}\Delta(\cdot)\mathcal{R}^{-1}, \label{eq:qtqb1}\\
			(\mathrm{id}\otimes\Delta)(\mathcal{R})
			&=\Phi_{312}^{-1}\mathcal{R}_{13}\Phi_{213}\mathcal{R}_{12}\Phi^{-1},\label{eq:qtqb2}\\
			(\Delta\otimes\mathrm{id})(\mathcal{R})
			&=\Phi_{231}\mathcal{R}_{13}\Phi_{132}^{-1}\mathcal{R}_{23}\Phi \label{eq:qtqb3}
		\end{align}
		hold, where we employed the well-known leg notation in $H^{\otimes 3}$, e.g. $\mathcal{R}_{12}=\mathcal{R}\otimes 1$, etc. In case $\mathcal{R}^{-1}=\mathcal{R}^\mathrm{op}$, where $\mathcal{R}^\mathrm{op}:=\tau(\mathcal{R})$, we call $(H,\mathcal{R})$ \textbf{triangular}. A \textbf{morphism of quasitriangular quasi-bialgebras} is a morphism of the underlying quasi-bialgebras $\alpha : H \to H'$ such that $(\alpha \ten \alpha)(\mathcal{R}) = \mathcal{R}'$.
	\end{enumerate}
\end{definition}

\noindent
For the trivial re-associator $\Phi=1\otimes 1\otimes 1$, Definition
\ref{def:quasi-bialg} recovers the notion of quasitriangular
bialgebra, which are prominent objects in Hopf algebra theory and thus
provide numerous examples. The generalization to quasitriangular
quasi-bialgebras turns out to be crucial for deformations via gauge
transformations. Namely, even starting with a quasitriangular
bialgebra one cannot expect to obtain a bialgebra as a deformation,
but rather a quasi-bialgebra. We sometimes omit $\ell$ and $r$ for brevity.
\begin{example} 
  We now collect some examples of quasitriangular quasi-bialgebras:
  \begin{enumerate}
  \item[i.)] Let $\mathfrak{g}$ be a complex Lie algebra together with
    an invariant and symmetric element $\Omega \in \mathfrak{g} \ten
    \mathfrak{g}$.  Then for any Drinfeld associator $\Psi$ (see the next
    section) there is a (topological) quasitriangular quasi-bialgebra
    structure on $\mathsf{U}(\mathfrak{g})[[\hbar]]$ with $\Phi =
    \Psi(\hbar \Omega_{12},\hbar  \Omega_{23})$ and $\mathcal{R} =
    e^{\frac{\hbar}{2} \Omega}$. This, in particular, applies to any
    complex semisimple Lie algebra, where the invariant and symmetric tensor is
    the universal Casimir element, see \cite[VI.22]{Hum} for more
    details.  For example, in the case of $\mathfrak{sl}_2$ we have
    $\Omega = e \ten f + f \ten e + \frac{1}{2} h \ten h$.
  \item[ii.)] Following \cite[Example 10.7]{Bul}, let $\Bbbk$ be a field of
    characteristic different from 2 containing an element $i$ such
    that $i^4 = 1$.  Denote by $H(2)$ the group algebra of the cyclic
    group generated by one element $g$, subject to the relation $g^2=1$. Then, setting $p_- =
    \frac{1}{2} (1-g) $, we have that $H(2)$ is a quasi-bialgebra with
    $\Phi = 1 \ten 1 \ten 1 - 2p_- \ten p_- \ten p_-$.  Moreover,
    $H(2)$ has exactly two quasitriangular structures given by $
    \mathcal{R}_\pm = 1 \ten 1- (1 \pm i) p_- \ten p_- $. Note also
    that $(\mathcal{R}_\pm)^{-1} = 1 \ten 1- (1 \pm i^3) p_- \ten
    p_-$.
  \item[iii.)] Another significant class of examples is discussed in
    \cite{Majid98}, where a quantum double $D(H)$ associated to any
    quasi-Hopf algebra $H$ has been constructed. The quantum double
    $D(H)$ turns out to be a quasitriangular quasi-Hopf algebra.
  \end{enumerate}
\end{example}
\noindent
A tool to construct (quasitriangular) quasi-bialgebras is given by
gauge transformations, in the terminology of \cite[Chapter XV.4]{Kassel95}. They are generalizations of Drinfeld twists
\cite{DrKZ}, where one drops the $2$-cocycle and normalization conditions.

\begin{definition}[{\cite[Section 1.2]{Doikou22}}]
	A \textbf{gauge transformation} of a quasi-bialgebra
        $(H,\Delta,\varepsilon,\Phi)$ is an invertible element
        $\mathcal{F}\in H\otimes H$.
\end{definition}
\noindent
Given a gauge transformation $\mathcal{F}\in H\otimes H$ of a
quasi-bialgebra $(H,\Delta,\varepsilon,\Phi)$, we define a
linear map $\Delta_\mathcal{F}\colon H\to H\otimes H$ by
\begin{align}
	\Delta_\mathcal{F}(\cdot):=\mathcal{F}\Delta(\cdot)\mathcal{F}^{-1}.
\end{align}
\begin{proposition}[{\cite[Proposition 1.7]{Doikou22}}]\label{propGauge}
Let $\mathcal{F}\in H\otimes H$ be a gauge transformation of a quasi-bialgebra $(H,\Delta,\varepsilon,\Phi,\ell,r)$. Then
$$
H_\mathcal{F}:=(H,\Delta_\mathcal{F},\varepsilon,\Phi_\mathcal{F},\ell_\mathcal{F},r_\mathcal{F})
$$
is a quasi-bialgebra, with
\begin{equation}
	\Phi_\mathcal{F}:=\mathcal{F}_{23}(\mathrm{id}\otimes\Delta)(\mathcal{F})\Phi(\Delta\otimes\mathrm{id})(\mathcal{F}^{-1})\mathcal{F}_{12}^{-1},\qquad\ell_\mathcal{F}:=\ell v^{-1},\qquad r_\mathcal{F}:=rw^{-1},
\end{equation}
where $v:=(\varepsilon\otimes\mathrm{id})(\mathcal{F})\in H$ and $w:=(\mathrm{id}\otimes\varepsilon)(\mathcal{F})\in H$.

\noindent
If $(H,\mathcal{R})$ is (quasi)triangular in addition, then $(H_\mathcal{F},\mathcal{R}_\mathcal{F})$ is (quasi)triangular with
\begin{equation}
	\mathcal{R}_\mathcal{F}:=\mathcal{F}^\mathrm{op}\mathcal{R}\mathcal{F}^{-1}.
\end{equation}
\end{proposition}
\noindent
In particular, the previous proposition enables us to construct
examples of quasitriangular quasi-bialgebras from quasitriangular
bialgebras. We recall the approach via invertible group-like elements,
taken from \cite[Example 1.11]{Doikou22}, in the following.
\begin{corollary}\label{cor:twist}
  Let $(H,\Delta,\varepsilon,\mathcal{R})$ be a quasitriangular bialgebra
  and $g\in H$ be an invertible group-like element. Then
  $\mathcal{F}:=1\otimes g\in H\otimes H$ is a gauge transformation
  and
  $(H,\Delta_\mathcal{F},\varepsilon,\Phi,\mathcal{R}_\mathcal{F})$ is
  a quasitriangular quasi-bialgebra with re-associator
  \begin{equation}
    \Phi:=1\otimes 1\otimes g^{-1}
  \end{equation}
  and $\ell=g$, $r=1$.
\end{corollary}
\noindent
We discuss a particular class of examples constructed via the previous
corollary. This is based on the family of $2^{n+1}$-dimensional
pointed Hopf algebras referred to as $E(n)$, introduced in
\cite{CaenDasc}. Note that $E(1)$ is the well-known 4-dimensional Hopf algebra of Sweedler.
\begin{example}\label{ex:E(n)}
  Let 
  $n\geq 1$ be an integer. We define $E(n)$ as the algebra generated by
  $g,x_1,\ldots,x_n$ modulo the relations
  \begin{equation}\label{Eq:EnRel}
    g^2=1,\qquad
    x_ig=-gx_i,\qquad
    x_jx_i=-x_ix_j
  \end{equation}
  for all $1\leq i,j\leq n$. In particular $x_i^2=0$, implying that
  $E(n)$ is a vector space of dimension $2^{n+1}$, with basis
  $\{x_1^{k_1}\ldots x_n^{k_n},gx_1^{\ell_1}\ldots
  x_n^{\ell_n}\}_{k_1,\ldots,k_n,\ell_1,\ldots,\ell_n\in\{0,1\}}$. The
  Hopf algebra structure of $E(n)$ is determined by
  \begin{equation*}
    \begin{split}
      \Delta(g)&=g\otimes g,\\
      \varepsilon(g)&=1,\\
      S(g)&=g,
    \end{split}\qquad\qquad
    \begin{split}
      \Delta(x_i)&=x_i\otimes 1+g\otimes x_i,\\
      \varepsilon(x_i)&=0,\\
      S(x_i)&=-gx_i
    \end{split}
  \end{equation*}
  for all $1\leq i\leq n$. 
  For every $n\times n$-matrix $(a_{ij})\in M_n(\Bbbk)$ with values in
  $\Bbbk$ there is a quasitriangular structure
  \begin{equation}
    \mathcal{R}_{(a_{ij})}:=\frac{1}{2}\bigg(1\otimes 1+1\otimes g+g\otimes 1-g\otimes g\bigg)\exp\bigg(\sum_{i,j=1}^na_{ij}gx_i\otimes x_j\bigg)
  \end{equation}
  on $E(n)$ (note that the above exponential converges because the exponent is nilpotent) and one shows that for every quasitriangular structure $\mathcal{R}$ of $E(n)$ there is a matrix $(a_{ij})\in M_n(\Bbbk)$ such that $\mathcal{R}=\mathcal{R}_{(a_{ij})}$. Note that 
  \begin{align*}
    \bigg(\sum_{i,j=1}^na_{ij}gx_i\otimes x_j\bigg)^k
    &=
    \sum_{i_1,\dots,i_k,j_1,\dots,j_k=1}^n  a_{i_1j_1}\cdots a_{i_kj_k} gx_{i_1}\cdots gx_{i_k}\otimes x_{j_1}\cdots x_{j_k}\\&
    =
    \sum_{i_1,\dots,i_k,j_1,\dots,j_k=1}^n (-1)^{\tfrac{k(k-1)}{2}} a_{i_1j_1}\cdots a_{i_kj_k} g^kx_{i_1}\cdots x_{i_k}\otimes x_{j_1}\cdots x_{j_k}\\&
    =\sum_{P,F\subseteq \{1,\dots, n\}} \sum_{\sigma,\tau\in S_k} (-1)^{\tfrac{k(k-1)}{2}}
    a_{p_{\sigma(1)}f_{\tau(1)}}\cdots a_{p_{\sigma(k)}f_{\tau(k)}} g^k x_{p_{\sigma(1)}}\cdots x_{p_{\sigma(k)}}\otimes 
    x_{f_{\tau(1)}}\cdots x_{f_{\tau(k)}} \\&
    =\sum_{P,F\subseteq \{1,\dots, n\}} \sum_{\sigma,\tau\in S_k} (-1)^{\tfrac{k(k-1)}{2}}
    \mathrm{sign}(\sigma)\mathrm{sign}(\tau)
    a_{p_{\sigma(1)}f_{\tau(1)}}\cdots a_{p_{\sigma(k)}f_{\tau(k)}} g^k x_{p_{1}}\cdots x_{p_{k}}\otimes 
    x_{f_1}\cdots x_{f_{k}}\\&
    =k!\sum_{P,F\subseteq \{1,\dots, n\}}  (-1)^{\tfrac{k(k-1)}{2}}
    \det((a_{p_i,f_j})) g^k x_{p_{1}}\cdots x_{p_{k}}\otimes 
    x_{f_1}\cdots x_{f_{k}},
  \end{align*}
  where we only used the relations \eqref{Eq:EnRel} and
  $P=\{p_1,\dots, p_k\}$ as well as $F=\{f_1,\dots,f_k\}$, and we
  denoted by $S_k$ the symmetric group on $k$ elements. The right hand
  side of the above recovers the expression of the
  $\mathcal{R}$-matrix appearing in \cite{BeaDasGru}.
  
\noindent
Applying Corollary \ref{cor:twist} to $(E(n),\mathcal{R})$ for the
  invertible group-like element $g\in E(n)$ gives quasitriangular
  quasi-bialgebras
  $(E(n),\Delta_\mathcal{F},\varepsilon,\Phi,\mathcal{R}_\mathcal{F})$,
  where $\mathcal{F}=1\otimes g$ and $\Phi=1\otimes 1\otimes g$. More
  explicitly, the twisted comultiplication of $E(n)$ reads
  \begin{equation}\label{twist:delta}
    \Delta_\mathcal{F}(g)=\Delta(g),\qquad\qquad
    \Delta_\mathcal{F}(x_i)=x_i\otimes 1-g\otimes x_i
  \end{equation}
  on generators, while the twisted universal $\mathcal{R}$-matrix reads
  \begin{equation}\label{twist:R}
    \begin{split}
      (\mathcal{R}_{(a_{ij})})_\mathcal{F}
      &=\mathcal{F}_{21}\mathcal{R}_{(a_{ij})}\mathcal{F}^{-1}
      =(g\otimes 1)\frac{1}{2}\bigg(1\otimes 1+1\otimes g+g\otimes 1-g\otimes g\bigg)\exp\bigg(\sum_{i,j=1}^na_{ij}gx_i\otimes x_j\bigg)(1\otimes g)\\
      &=-\frac{1}{2}\bigg(1\otimes 1-1\otimes g-g\otimes 1-g\otimes g\bigg)\exp\bigg(-\sum_{i,j=1}^na_{ij}gx_i\otimes x_j\bigg)
    \end{split}
  \end{equation}
  for all $n\times n$-matrices $(a_{ij})$.
\end{example}
\noindent
It is well-known that the representation categories of quasitriangular
quasi-bialgebras are braided monoidal. In fact, they constitute our
main class of examples. We thus connect quasitriangular
quasi-bialgebras to the notions introduced in the previous section.
\begin{proposition}\label{prop:quasi-bialgCAT}
  Let $H$ be an associative unital algebra. The following statements hold true.
  \begin{enumerate}
  \item[i.)] The representation category ${}_H\mathcal{M}$ has a monoidal structure $(\otimes,\Bbbk,a,\ell,r)$ and is endowed with a strong monoidal functor to the category of $\Bbbk$-vector spaces, if and only if $(H,\Delta,\varepsilon,\Phi,\ell,r)$ is a quasi-bialgebra, where the coproduct and counit correspond to the tensor product and trivial representation, respectively. 
In this case, the associativity and unit constraints
    read
    \begin{equation}
      \begin{split}
	a_{X,Y,Z}((x\otimes y)\otimes z)
	&=\Phi\cdot(x\otimes (y\otimes z)),\\
	\ell_X(1\otimes x)
	&=\ell\cdot x,\\
	r_X(x\otimes 1)
	&=r\cdot x
      \end{split}
    \end{equation}
    on left $H$-modules $X,Y,Z$. 
  \item[ii.)] The above monoidal category is braided if and only
    if the quasi-bialgebra $(H,\Delta,\varepsilon,\Phi,\ell,r)$ is
    quasitriangular. In this case, the braiding is determined by
    \begin{equation}
      \sigma_{X,Y}(x\otimes y)
      :=\mathcal{R}_i\cdot y\otimes\mathcal{R}^i\cdot x 
    \end{equation}
    on left $H$-modules $X,Y$, where
    $\mathcal{R}=\mathcal{R}^i\otimes\mathcal{R}_i\in H\otimes H$ is
    the corresponding universal $\mathcal{R}$-matrix of $H$. The
    category is symmetric if and only if $\mathcal{R}$ is triangular.

  \item[iii.)] In the setup of $i.)$, if $\mathcal{F}\in H\otimes H$ is a gauge
    transformation, then the monoidal categories $({}_H\mathcal{M},\otimes,\Bbbk,a,\ell,r)$ and
    $({}_{H_\mathcal{F}}\mathcal{M},\otimes_\mathcal{F},\Bbbk,a_\mathcal{F},\ell,r)$
    are equivalent via the monoidal functor
    $(\mathrm{id},F^2,\mathrm{id})$, where
    \begin{equation}
      F^2_{X,Y}(x\otimes y):=\mathcal{F}^{-1}\cdot(x\otimes y)
    \end{equation}
    for all left $H$-modules $X,Y$. Above, $\otimes_\mathcal{F}$
    denotes the $\Bbbk$-vector spaces tensor product of left
    $H_\mathcal{F}$-modules, structured as a left
    $H_\mathcal{F}$-module via $\Delta_\mathcal{F}$, and
    $a_\mathcal{F}$ is the associativity constraint corresponding to
    $\Phi_\mathcal{F}$. If $H$ is quasitriangular, then the monoidal
    equivalence is braided monoidal.
  \end{enumerate}
\end{proposition}
\noindent
For proofs of these facts we refer the reader to \cite[Chapter XV]{Kassel95}.

\begin{remark}\label{rem:topfree}
It is straightforward to generalize the entirety of this section to (quasitriangular) topological bialgebras. For instance, in Definition \ref{def:quasi-bialg} one simply considers objects and morphisms in the symmetric monoidal category of topologically-free modules instead of the category of vector spaces, see \cite[Chapter XVI]{Kassel95}. This generalization will be needed in Section \ref{sec:4}. Nevertheless, we chose to recall the theory in the vector space setup for clarity of presentation.
\end{remark}

\subsection{Infinitesimal braid relations and Drinfeld associators}\label{sec:2.3}

For the convenience of the reader we briefly recall the basics of
Drinfeld associators \cite{DR89,DrQH}. Let $n>2$ be a natural number.

\begin{definition}
  Let $H$ be an algebra. A collection of elements $\{A^{ij}\}_{1 \leq
    i \neq j \leq n}$ in $H$ satisfies the \textbf{infinitesimal braid
    relations} if
  \begin{align}
  	A^{ij}-A^{ji} &=0 \qquad \text{for }\#\{i,j\}=2 \label{eq:infbraid0}\\
    [A^{ij}, A^{ik} + A^{jk}] &= 0 \qquad \text{for } \# \{i,j,k \}=3
    \label{eq:infbraid1}\\   
    [A^{ij}, A^{kh}] &=0 \qquad \text{for } \# \{i,j,k,h \}=4. 
    \label{eq:infbraid2} 
  \end{align}
\end{definition}
\noindent
In concrete computations we usually specify the elements $A^{ij}$ with $i<j$ since, because of condition \eqref{eq:infbraid0}, the full set is determined by these elements.
The infinitesimal braid relations are the generating relations of the well-known Drinfeld--Kohno algebras \cite{Kohno}.
\begin{definition}
	\label{definition-Drinfeld-associator}
	Let $\lambda \in \Bbbk$ be an invertible element. A
        \textbf{Drinfeld associator} is a formal power
        series $\Psi(A,B) \in \Bbbk \langle \langle A,B \rangle
        \rangle$ in two noncommuting variables $A,B$ such that:
	\begin{itemize}
		\item[i.)] For any algebra $H$ together with elements
                  $\{A^{ij}\}_{1 \leq
                  		i \neq j \leq 4}$ satisfying the
                  infinitesimal braid relations, the \textbf{pentagon
                    equation}
		\begin{equation}
			\label{eq:pentagon-equation}
			\Psi(\hbar A^{12}, \hbar A^{23} + \hbar A^{24}) \Psi (\hbar A^{13} + \hbar A^{23}, \hbar A^{34}) = \Psi (\hbar A^{23}, \hbar A^{34}) \Psi (\hbar A^{12} + \hbar A^{13}, \hbar A^{24} + \hbar A^{34}) \Psi (\hbar A^{12}, \hbar A^{23})
		\end{equation}
		holds in $H[[\hbar]]^{\tilde{\ten} 4}$. 
		
		\item[ii.)] For any algebra $H$ together with elements $A,B,C$ and $\Lambda \coloneqq A+B+C$ satisfying $[\Lambda , A] = [\Lambda , B] = [\Lambda, C] = 0$,
		the \textbf{hexagon equation}
		\begin{equation}
			\label{eq:hexagon-equation}
			e^{\lambda \hbar \pi i \Lambda} = e^{\lambda \pi i \hbar A}  \Psi(\hbar C,\hbar A)e^{\lambda \hbar \pi i C} \Psi(\hbar B,\hbar C) e^{\lambda \hbar \pi i B} \Psi(\hbar A,\hbar B)
		\end{equation}
		holds in $H[[\hbar]]^{\tilde{\ten} 3}$.
		\item[iii.)] For any $A,B$ one has 
		\begin{equation}
			\label{eq:inverse-associator}
			\Psi(A,B)^{-1} = \Psi(B,A).
		\end{equation}
		\item[iv.)] $\Psi$ is a group-like element with
                  respect to the comultiplication of $ \Bbbk \langle
                  \langle A,B \rangle \rangle$ generated by $\Delta(A)
                  = 1 \otimes A + A \otimes 1$ and $\Delta(B)= 1
                  \otimes B + B \otimes 1$.
	\end{itemize}
\end{definition}
\begin{remark}
\label{remark-associator-A-B-commute}
	If $[A,B]=0$ then $\Psi(A,B)=1$.
\end{remark}
\noindent
We denote Drinfeld associators by $\Psi$, even though the symbol $\Phi$ is more common in the literature. This is in order to distinguish Drinfeld associators from re-associators.
The most prominent example of a Drinfeld associator is the
Knizhnik-Zamolodchikov one, which exists in the complex framework, see
also \cite[Chapter XIX]{Kassel95}, and the recent article \cite{BRW}
for more details. Namely, Drinfeld \cite{DR89,DrQH} showed, by
computing the monodromy of the Knizhnik-Zamolodchikov connection, that
there exists a complex Drinfeld associator. 
Moreover, it
has been proven in \cite{DR89} that there exists a Drinfeld associator
with coefficients in the field of rational numbers $\mathbb{Q}$. This
result implies the existence of a Drinfeld associator in any field of
characteristic zero.

\section{Infinitesimal braidings and pre-Cartier quasi-bialgebras}

In this chapter we add yet another ingredient to braided monoidal
categories: an infinitesimal braiding. This will give us a regime
suitable for quantization. After discussing the categorical setup in
Section \ref{sec:pC} we focus on the algebraic realm in Section
\ref{sec:pCbialg}, via the aforementioned duality.

\subsection{Pre-Cartier categories}\label{sec:pC}
We now introduce non-strict analogues of pre-Cartier braided monoidal
categories, generalizing \cite[Section 1.2]{ABSW}.
Recall that a category is pre-additive if it is enriched over the category of abelian groups.

\begin{definition}\label{def:preCart}
	A \textbf{pre-Cartier category} $(\mathcal{C}, \otimes, I, a,
        \ell,r, \sigma, t)$ is a pre-additive braided monoidal
        category \\ $(\mathcal{C}, \otimes, I, a, \ell,r, \sigma)$,
        together with a natural transformation $t: \otimes \Rightarrow
        \otimes$, called the \textbf{infinitesimal braiding}, such that the identities
		\begin{align}
			t_{X, Y\otimes Z} &= a_{X,Y,Z} \circ (t_{X,Y} \otimes \mathrm{id}_Z) \circ a_{X,Y,Z}^{-1} \nonumber \\
			& \ + a_{X,Y,Z} \circ (\sigma^{-1}_{X,Y} \otimes \mathrm{id}_Z) \circ a^{-1}_{Y,X,Z} \circ (\mathrm{id}_Y \otimes t_{X,Z}) \circ a_{Y,X,Z} \circ (\sigma_{X,Y} \otimes \mathrm{id}_Z) \circ a^{-1}_{X,Y,Z} \label{eq:pre-cartier-one} \\
			t_{X \otimes Y, Z} &= a^{-1}_{X,Y,Z} \circ (\mathrm{id}_X \otimes t_{Y,Z}) \circ a_{X,Y,Z} \nonumber \\
			& \ +a_{X,Y,Z}^{-1} \circ (\mathrm{id}_X \otimes \sigma^{-1}_{Y,Z}) \circ a_{X,Z,Y} \circ (t_{X,Z} \otimes \mathrm{id}_Y) \circ a^{-1}_{X,Z,Y} \circ (\mathrm{id}_X \otimes \sigma_{Y,Z}) \circ a_{X,Y,Z} \label{eq:pre-cartier-two}
		\end{align}
	hold true for all objects $X,Y,Z$ of $\mathcal{C}$. If, moreover, the identity
	\begin{equation}
		\label{eq:cartier-category}
		\sigma_{X,Y} \circ t_{X,Y} = t_{Y,X} \circ \sigma_{X,Y}
	\end{equation}
	holds true for all objects $X,Y$ of $\mathcal{C}$, we say that
        $(\mathcal{C}, \otimes, I, a, \ell,r, \sigma, t)$ is a
        \textbf{Cartier category}. In case $\sigma_{X,Y} =
        \sigma_{Y,X}^{-1}$ we refer to the above as
        \textbf{(pre-)Cartier symmetric categories}.
\end{definition}
\noindent
Cartier categories are named after P. Cartier, who
introduced them in the context of symmetric braidings \cite{Car}. These are usually referred to as infinitesimally braided monoidal categories and the corresponding natural transformation is also called infinitesimal braiding. We, instead, call such categories Cartier symmetric and add the
prefix "pre-" to indicate that condition
\eqref{eq:cartier-category} is not assumed in general.
\begin{example} 
  \label{examples-precartier-categories}
  The following are examples of pre-Cartier categories:
  \begin{enumerate}
  \item[i.)] Every pre-additive braided monoidal category is Cartier with respect to the trivial infinitesimal braiding, defined on objects $X$, $Y$ by $t_{X,Y}=0$. Moreover, $\mathbb{Z}$-fold multiples (or $\Bbbk$-multiples if the category is $\Bbbk$-linear) of infinitesimal braidings are infinitesimal braidings with respect to the same underlying braided monoidal category.
    
  \item[ii.)] Given a braided monoidal category $(\mathcal{C}, \otimes, I, a, \ell,r, \sigma)$, it is well-known that $(\mathcal{C}, \otimes, I, a, \ell,r, \sigma^{-1})$ is braided monoidal. If the former is pre-Cartier with respect to $t$, so is the latter, with respect to the infinitesimal braiding $t':=\sigma\circ t\circ\sigma^{-1}$. Note that this construction is involutive, meaning that performing it a second time results in the initial pre-Cartier category. Further note that in the Cartier case $t$ and $t'$ coincide.
    
  \item[iii.)] Let $(\mathfrak{b}, [\cdot,\cdot],\delta)$ be a Lie
    bialgebra. Recall that a $\mathfrak{b}$-dimodule (or
    Drinfeld-Yetter module over $\mathfrak{b}$) is a triple $(V, \pi,
    \pi^*)$, where $V$ is a vector space and $\pi: \mathfrak{b} \ten V
    \to V$, $\pi^*:V \to \mathfrak{b} \ten V$ are linear maps
    satisfying
    \begin{align*}
      \pi \circ ([\cdot, \cdot] \ten \id_V) &= \pi \circ (\id_\mathfrak{b} \ten \pi) - \pi \circ (\id_\mathfrak{b} \ten \pi) \circ (\tau_{\mathfrak{b},\mathfrak{b}} \ten \id_V) \\
      (\delta \ten \id_V) \circ \pi^* &= (\tau_{\mathfrak{b},\mathfrak{b}} \ten \id_V) \circ (\id_\mathfrak{b} \ten \pi^*) \circ \pi^* -  (\id_\mathfrak{b} \ten \pi^*) \circ \pi^* \\
      \pi^* \circ \pi &= (\id_\mathfrak{b} \ten \pi)\circ (\tau_{\mathfrak{b},\mathfrak{b}} \ten \id_V) \circ (\id_\mathfrak{b} \ten \pi^*) + ([\cdot, \cdot] \ten \id_V)\circ (\id_\mathfrak{b} \ten \pi^*) - (\id_\mathfrak{b} \ten \pi) \circ (\delta \ten \id_V),
    \end{align*}
    where $\tau$ denotes the tensor flip. One can show that category
    of $\mathfrak{b}$-dimodules is a symmetric Cartier category, where
    the constraints $a, \ell,r,\sigma$ are the standard ones (i.e. the
    same of $\mathrm{Vec}_\Bbbk$), and the infinitesimal braiding is
    \[ t_{V,W} =  (\id_V \ten \pi_W)\circ (\tau_{\mathfrak{b},V} \ten \id_W)\circ (\pi_V^* \ten \id_W) + (\pi_V \ten \id_W)\circ(\tau_{V,\mathfrak{b}} \ten \id_W)\circ(\id_V \ten \pi^*_W). \]
    Further details can be found in \cite[Section 4.4]{AndreaThesis}. 
  \item[iv).] For any Lie algebra $\mathfrak{g}$ the category of left $\mathfrak{g}$-modules ${}_\mathfrak{g}\mathcal{M}$ is a symmetric monoidal category with respect to the tensor flip. Then, pre-Cartier structures of ${}_\mathfrak{g}\mathcal{M}$ are in 1:1-correspondence with $\mathfrak{g}$-invariant tensors $C\in\mathfrak{g}\otimes\mathfrak{g}$, where the infinitesimal braiding appears as the representation of $C$. The category is Cartier if and only if $C$ is a symmetric tensor. This can be deduced along the same lines of \cite[Proposition~XX.4.2]{Kassel95}. In particular, in case $\mathfrak{g}=\mathfrak{g}_S\oplus\mathfrak{a}$ is finite-dimensional, complex and reductive, one may consider $C=c+a$, where $c$ is the Casimir element of the semisimple part $\mathfrak{g}_S$ and $a\in\mathfrak{a}\otimes\mathfrak{a}$ is an arbitrary tensor of the abelian part. Then $C$ gives rise to a pre-Cartier structure on ${}_\mathfrak{g}\mathcal{M}$, which is Cartier if and only if $a$ is symmetric. 
  \end{enumerate}
\end{example}
\noindent
In order to obtain more examples of pre-Cartier categories, we introduce pre-Cartier
quasi-bialgebras in the subsequent section.  In the following we
sometimes abbreviate a pre-Cartier category by $(\mathcal{C}, \otimes,
\sigma, t)$ or simply $(\mathcal{C},\sigma,t)$, if the associativity
constraint, unit constraints or even the monoidal structure are clear
from the context. 
\begin{remark}
\label{remark-precartier-categories}
	We observe that the definition of pre-Cartier category has
        some immediate consequences.
	\begin{enumerate}
		\item[i.)] Any infinitesimal braiding satisfies
                  $t_{X,I}=0=t_{I,X}$ for any object $X$. This follows immediately
                  from \eqref{eq:pre-cartier-one} and
                  \eqref{eq:pre-cartier-two}.
		\item[ii.)] In case of a pre-additive symmetric
                  monoidal category $(\mathcal{C}, \otimes, I, a, \ell,r,
                  \sigma)$ satisfying \eqref{eq:cartier-category} it
                  follows that a natural transformation $t: \otimes
                  \Rightarrow \otimes$ satisfies
                  \eqref{eq:pre-cartier-one} if and only if
                  \eqref{eq:pre-cartier-two} holds.
		\item[iii.)] If $(\mathcal{C}, \otimes, I, a, \ell,r,
                  \sigma)$ is strict, one recovers
                  \cite[Definition 1.1]{ABSW}.
                  
        \item[iv.)] We mentioned in the introduction that infinitesimal $\mathcal{R}$-matrices arise as first-order deformations of $\mathcal{R}$-matrices. Similarly, notice that infinitesimal braidings are precisely first-order deformations of braidings, cf. \cite{PulSev} for the Cartier case.
	\end{enumerate}
\end{remark}
\noindent
A functor of pre-Cartier categories is supposed to be a functor of the
underlying braided monoidal structure, which also respects the
infinitesimal braidings.
\begin{definition}
	Let $(\mathcal{C}, \otimes, I, a, \ell,r, \sigma, t)$ and
        $(\mathcal{C}', \otimes', I', a', \ell',r', \sigma', t')$ be two
        pre-Cartier categories. An \textbf{infinitesimally braided
          monoidal functor} is a braided monoidal functor
        $(F,F^2,F^0)$, satisfying
		\begin{equation}\label{diag:Functor-infbraid}
			\begin{tikzcd}
				F(X) \otimes' F(Y) \arrow[r, "{F^2(X,Y)}"] \arrow[d, "{t'_{F(X),F(Y)}}"'] & F(X \otimes Y) \arrow[d, "{F(t_{X,Y})}"] \\
				F(X) \otimes' F(Y) \arrow[r, "{F^2(X,Y)}"]                                & F(X \otimes Y)                          
			\end{tikzcd}
		\end{equation}
for all objects $X,Y$ in $\mathcal{C}$.
\end{definition}
\begin{remark}
	\label{composition-of-ibmf-is-ibmf}
	One verifies that composition of infinitesimally braided
        monoidal functors provides an infinitesimally braided monoidal
        functor.
\end{remark}
\noindent
Given a monoidal category $(\mathcal{C},\otimes)$ and a monoidal
functor $F$ to a pre-Cartier category
$(\mathcal{C}',\otimes',\sigma',t')$ one might attempt to determine a
braiding $\sigma$ on $\mathcal{C}$ by enforcing
\eqref{eq:braided-monoidal-functor} and, in the same spirit, to
determine an infinitesimal braiding $t$ on $\mathcal{C}$ by declaring
\eqref{diag:Functor-infbraid}. We show that this is in fact possible
under some assumptions on the functor.
\begin{proposition}
	Let $(F, F^2, F^0) \colon \mathcal{C} \to \mathcal{C}'$ be a
        strong monoidal functor between pre-additive categories
        $\mathcal{C}$ and $\mathcal{C}'$. Assume $F$ is additive and
        fully faithful.
	\begin{itemize}
		\item[i)] If $\mathcal{C}'$ has a braiding (symmetry),
                  then so does $\mathcal{C}$ and $F$ is braided strong
                  monoidal.
		\item[ii)] If $\mathcal{C}'$ has a natural
                  transformation fulfilling \eqref{eq:pre-cartier-one}
                  and \eqref{eq:pre-cartier-two}, then so does
                  $\mathcal{C}$.
		\item[iii)] If $\mathcal{C}'$ has a natural
                  transformation fulfilling
                  \eqref{eq:cartier-category}, then so does
                  $\mathcal{C}$.
	\end{itemize}
\end{proposition}
\begin{proof}
	The proof is straightforward and it follows the same lines of
        \cite[Proposition 1.3]{ABSW}.
\end{proof}
\noindent
The above proposition can thus be used as a tool to construct
pre-Cartier categories from known examples. We will encounter a
ramification of this in the form of gauge equivalences in the
following section.

\subsection{Pre-Cartier quasi-bialgebras}\label{sec:pCbialg}

One of the most interesting examples of pre-Cartier category is given
by the representation category of quasitriangular quasi-bialgebras,
endowed with a pre-Cartier structure. In the following we introduce
the latter and discuss the induced categorical counterpart.

\noindent
Consider a quasitriangular quasi-bialgebra
$(H,\Delta,\varepsilon,\Phi,\mathcal{R})$ and the trivial topological quasi-bialgebra
	$\tilde{H}:=H[[\hbar]]$ of formal power series in a formal
	parameter $\hbar$, where we $\hbar$-linearly extended the
	algebra structure and $\Delta,\varepsilon$ and use $\Phi\in
	H^{\otimes 3}\subset\tilde{H}^{\otimes 3}$ as the
	re-associator. Assume there is a quasitriangular structure
	$\tilde{\mathcal{R}}\in(H\otimes
	H)[[\hbar]]\cong\tilde{H}\tilde{\otimes}\tilde{H}$ on
	$\tilde{H}$. It follows that its zeroth order $\mathcal{R}\in
	H\otimes H$ is a universal $\mathcal{R}$-matrix for $H$. Then,
	factorizing the invertible element $\mathcal{R}$, we obtain
	\begin{equation}
		\label{eq:quantization-problem}
		\tilde{\mathcal{R}}=\mathcal{R}(1\otimes 1+\hbar\chi+\mathcal{O}(\hbar^2))
	\end{equation}
	for an element $\chi\in H\otimes H$. Reading the axioms of
	$\tilde{\mathcal{R}}$ in first order of $\hbar$ imposes conditions on $\chi$, which we collect in the following definition.
	
\begin{definition}\label{def:preCartquasibialg}
	We call $(H,\chi)$ a \textbf{pre-Cartier
          quasi-bialgebra} if $(H,\Delta,\varepsilon,\Phi,\mathcal{R})$ is a quasitriangular quasi-bialgebra and there is an element $\chi\in H\otimes
        H$, an \textbf{infinitesimal $\mathcal{R}$-matrix}, such that
	\begin{align}
		\chi\Delta(\cdot)
		&=\Delta(\cdot)\chi,\label{eq:pC1}\\
		(\mathrm{id}\otimes\Delta)(\chi)
		&=\Phi\mathcal{R}_{12}^{-1}\Phi_{213}^{-1}\chi_{13}\Phi_{213}\mathcal{R}_{12}\Phi^{-1}
		+\Phi\chi_{12}\Phi^{-1},\label{eq:pC2}\\
		(\Delta\otimes\mathrm{id})(\chi)
		&=\Phi^{-1}\mathcal{R}_{23}^{-1}\Phi_{132}\chi_{13}\Phi_{132}^{-1}\mathcal{R}_{23}\Phi
		+\Phi^{-1}\chi_{23}\Phi\label{eq:pC3}
	\end{align}
	hold. If in addition
	\begin{equation}
		\mathcal{R}\chi=\chi^\mathrm{op}\mathcal{R}\label{pC4}
	\end{equation}
	holds, we call $(H,\chi)$ a \textbf{Cartier quasi-bialgebra}.
\end{definition}
\noindent
Thus, we can understand pre-Cartier structures on quasitriangular
quasi-bialgebras as first order deformations of the trivial
topologically free bialgebra. The emerging deformation problem is one
of the guiding principles of this article. Namely, in
\cite[Question 2.10]{ABSW} it is asked whenever, given a pre-Cartier
bialgebra $(H, \mathcal{R}, \chi)$, there is a quasitriangular
structure $\tilde{\mathcal{R}}$ on the trivial topological bialgebra
$\tilde{H}$ satisfying \eqref{eq:quantization-problem}. 
We shall give an answer to this question in the context of pre-Cartier quasi-bialgebras.

\noindent
For the trivial re-associator $\Phi=1\otimes 1\otimes 1$, Definition
\ref{def:preCartquasibialg} recovers the notion of (pre-)Cartier
bialgebra given in \cite{ABSW}.

\noindent
As another a posteriori justification of Definition
\ref{def:preCartquasibialg}, anticipated in the previous section, we
show that the representation category of a quasitriangular
quasi-bialgebra is pre-Cartier precisely when the quasi-bialgebra
is. In this sense, the following extends Proposition
\ref{prop:quasi-bialgCAT}.
\begin{theorem}
	\label{Tannaka-Krein}
	Let $(H, \Phi, \mathcal{R})$ be a quasitriangular quasi-bialgebra and denote by $({}_H\mathcal{M},\otimes,\sigma)$ the corresponding braided monoidal category. Then there is a bijective correspondence between pre-Cartier structures on $(H, \Phi, \mathcal{R})$ and pre-Cartier structures on $({}_H\mathcal{M},\otimes,\sigma)$. Explicitly, given an infinitesimal $\mathcal{R}$-matrix $\chi\in H\otimes H$, the infinitesimal braiding reads
	\begin{equation}
		t_{X,Y}\colon X\otimes Y\to X\otimes Y,\qquad
		t_{X,Y}(x\otimes y)=\chi\cdot(x\otimes y),
	\end{equation}
	while $\chi$ is recovered from $t$ via $\chi=t_{H,H}(1\otimes 1)$.
\end{theorem}
\begin{proof} 
	Suppose that $(H, \mathcal{R}, \Phi)$ is pre-Cartier and
        consider $X,Y,Z$ objects in ${}_H\mathcal{M}$.  Recall by
        Proposition \ref{prop:quasi-bialgCAT} that the category
        ${}_H\mathcal{M}$ is braided monoidal with associativity
        constraint $a_{X,Y,Z}\big((x \otimes y) \otimes z\big) = \Phi
        \cdot \big(x \otimes(y \otimes z)\big)$ and braiding
        $\sigma_{X,Y}(x \otimes y) = \mathcal{R}^{\mathrm{op}} \cdot
        (y \otimes x).$ Next, set $t_{X,Y}(x \otimes y) = \chi \cdot
        (x \otimes y)$.  Then $t$ is a natural transformation (as
        shown in \cite[Theorem 2.6]{ABSW}). We have
	\begin{equation*}
		\begin{split}
			t_{X, Y \otimes Z}\big(x \otimes(y \otimes z)\big) &=\chi \cdot \big(x \otimes(y \otimes z)\big) \\
			&=\big((\mathrm{id} \otimes \Delta)(\chi)\big) \cdot \big(x \otimes (y \otimes z)\big) \\
			& =  \big(\Phi \mathcal{R}_{12}^{-1} \Phi_{213}^{-1} \chi_{13} \Phi_{213} \mathcal{R}_{12} \Phi^{-1}   + \Phi \chi_{12} \Phi^{-1} \big) \cdot \big(x \otimes(y \otimes z)\big) \\
			& = \Big( a_{X,Y,Z} \circ (\sigma^{-1}_{X,Y} \otimes \mathrm{id}_Z) \circ a^{-1}_{Y,X,Z}  \circ(\mathrm{id}_Y \otimes t_{X,Z}) \circ a_{Y,X,Z} \circ (\sigma_{X,Y} \otimes \mathrm{id}_Z) \circ a^{-1}_{X,Y,Z} \Big) \big(x \otimes (y \otimes z)\big) \\
			&  + \Big( a_{X,Y,Z} \circ (t_{X,Y} \otimes \mathrm{id}_Z) \circ a_{X,Y,Z}^{-1} \Big)\big(x \otimes (y \otimes z) \big)
		\end{split} 
	\end{equation*}
	and similarly one proves  \eqref{eq:pre-cartier-two}.
        Conversely, suppose that ${}_H\mathcal{M}$ is pre-Cartier.
        Then it is well-known that $\Phi = a_{H,H,H}\big((1 \otimes 1)
        \otimes 1\big)$ and $\mathcal{R}= \tau_{H,H} \circ
        \sigma_{H,H}(1 \otimes 1)$ endow $H$ with a quasitriangular
        quasi-bialgebra structure (see
        e.g. \cite[XV.1.2--XV.2.2]{Kassel95}). Next, set
        $\chi = t_{H,H}(1 \otimes 1)$.  Then we have that $\chi
        \Delta(\cdot) = \Delta(\cdot) \chi$ (see \cite[Theorem 2.6]{ABSW}). Using the naturality of $t$, we have
	\begin{equation*}
		\begin{split}
			&(\mathrm{id} \otimes \Delta)(\chi) = (\mathrm{id} \otimes \Delta) (t_{H,H}(1 \otimes 1)) = t_{H,H \otimes H}\big((\mathrm{id} \otimes \Delta)(1 \otimes 1)\big) = t_{H,H \otimes H}\big(1 \otimes (1 \otimes 1)\big) \\
			&= \Big(a_{H,H,H} \circ (t_{H,H} \otimes \mathrm{id}_H) \circ a_{H,H,H}^{-1} \Big) \big(1 \otimes (1 \otimes 1)\big) \\
			& \ + \Big( a_{H,H,H} \circ (\sigma^{-1}_{H,H} \otimes \mathrm{id}_H) \circ a^{-1}_{H,H,H} \circ (\mathrm{id}_H \otimes t_{H,H}) \circ a_{H,H,H} \circ (\sigma_{H,H} \otimes \mathrm{id}_H) \circ a^{-1}_{H,H,H}\Big) \big(1 \otimes (1 \otimes 1)\big) \\
			&= \Phi \chi_{12} \Phi^{-1} + \Phi \mathcal{R}_{12}^{-1} \Phi_{213}^{-1} \chi_{13} \Phi_{213} \mathcal{R}_{12} \Phi^{-1}
		\end{split}
	\end{equation*}
	and similarly \eqref{eq:pC3} follows. Bijectivity follows in complete analogy to \cite[Theorem 2.6]{ABSW}.
\end{proof}
\noindent
The previous theorem gives examples of pre-Cartier categories which
are neither symmetric, nor satisfy \eqref{eq:cartier-category} in
general.  We note some important specifications of Theorem
\ref{Tannaka-Krein}.
\begin{remark}
	Following the same argument, one can show that:
	\begin{itemize}
		\item[(i)] There is a bijective correspondence between Cartier structures on a
		quasitriangular quasi-bialgebra $(H, \Phi, \mathcal{R})$ and Cartier structures on the corresponding braided monoidal category $({}_H\mathcal{M},\otimes,\sigma)$. 
		\item[(ii)] There is a bijective correspondence between (pre-)Cartier structures on a
		triangular quasi-bialgebra $(H, \Phi, \mathcal{R})$ and (pre-)Cartier structures on the corresponding symmetric monoidal category $({}_H\mathcal{M},\otimes,\sigma)$. 
	\end{itemize}
\end{remark}

\noindent
Examples of pre-Cartier quasi-bialgebras certainly include pre-Cartier
bialgebras, as discussed in \cite{ABSW}. We show that there are
examples beyond the bialgebra realm, which require non-trivial
re-associators.

\noindent
It turns out that the notion of pre-Cartier quasi-bialgebra is stable
under gauge transformation.

\begin{proposition}\label{prop:twist}
	Let $(H,\Delta,\varepsilon,\Phi,\mathcal{R},\chi)$ be a
        (pre-)Cartier quasi-bialgebra and $\mathcal{F}\in H\otimes H$ a
        gauge transformation. Then $(H_\mathcal{F},\chi_\mathcal{F})$
	is a (pre-)Cartier quasi-bialgebra with
	\begin{equation}\label{chiF}
		\chi_\mathcal{F}:=\mathcal{F}\chi\mathcal{F}^{-1}\in H\otimes H,
	\end{equation}
	where
        $H_\mathcal{F}=(H,\Delta_\mathcal{F},\varepsilon,\Phi_\mathcal{F},\mathcal{R}_\mathcal{F})$
        is the quasitriangular quasi-bialgebra of Proposition
        \ref{propGauge}.
\end{proposition}
\begin{proof}
	According to Proposition \ref{propGauge}
        $(H_\mathcal{F},\mathcal{R}_\mathcal{F})$ is a
        quasitriangular quasi-bialgebra. We show that \eqref{chiF}
        satisfies the axioms \eqref{eq:pC1}, \eqref{eq:pC2} and
        \eqref{eq:pC3} for
        $(H_\mathcal{F},\mathcal{R}_\mathcal{F})$. In fact,
	\begin{align*}
		\chi_\mathcal{F}\Delta_\mathcal{F}(\cdot)
		&=\mathcal{F}\chi\mathcal{F}^{-1}\mathcal{F}\Delta(\cdot)\mathcal{F}^{-1}\\
		&=\mathcal{F}\Delta(\cdot)\chi\mathcal{F}^{-1}\\
		&=\Delta_\mathcal{F}(\cdot)\chi_\mathcal{F},
	\end{align*}
	using \eqref{eq:pC1} for $(H,\mathcal{R},\chi)$. Furthermore,
	\begin{align*}
		(\mathrm{id}&\otimes\Delta_\mathcal{F})(\chi_\mathcal{F})
		=\mathcal{F}_{23}(\mathrm{id}\otimes\Delta)(\mathcal{F})(\mathrm{id}\otimes\Delta)(\chi)(\mathrm{id}\otimes\Delta)(\mathcal{F}^{-1})\mathcal{F}_{23}^{-1}\\
		&=\mathcal{F}_{23}(\mathrm{id}\otimes\Delta)(\mathcal{F})
		\bigg(\Phi \mathcal{R}_{12}^{-1}\Phi_{213}^{-1}\chi_{13}\Phi_{213}\mathcal{R}_{12}\Phi^{-1}+\Phi\chi_{12}\Phi^{-1}\bigg)
		(\mathrm{id}\otimes\Delta)(\mathcal{F}^{-1})\mathcal{F}_{23}^{-1}\\
		&=\Phi^\mathcal{F}\mathcal{F}_{12}(\Delta\otimes\mathrm{id})(\mathcal{F})\mathcal{R}_{12}^{-1}\Phi_{213}^{-1}\chi_{13}\Phi_{213}\mathcal{R}_{12}(\Delta\otimes\mathrm{id})(\mathcal{F}^{-1})\mathcal{F}^{-1}_{12}(\Phi^\mathcal{F})^{-1}\\
		&\quad+\Phi^\mathcal{F}\mathcal{F}_{12}(\Delta\otimes\mathrm{id})(\mathcal{F})\chi_{12}(\Delta\otimes\mathrm{id})(\mathcal{F}^{-1})\mathcal{F}^{-1}_{12}(\Phi^\mathcal{F})^{-1}\\
		&=\Phi^\mathcal{F} \mathcal{F}_{12}\mathcal{R}_{12}^{-1}(\Delta^\mathrm{op}\otimes\mathrm{id})(\mathcal{F})\Phi_{213}^{-1}\chi_{13}\Phi_{213}(\Delta^\mathrm{op}\otimes\mathrm{id})(\mathcal{F}^{-1})\mathcal{R}_{12}\mathcal{F}^{-1}_{12}(\Phi^\mathcal{F})^{-1}\\
		&\quad+\Phi^\mathcal{F} \mathcal{F}_{12}\chi_{12}(\Delta\otimes\mathrm{id})(\mathcal{F})(\Delta\otimes\mathrm{id})(\mathcal{F}^{-1})\mathcal{F}^{-1}_{12}(\Phi^\mathcal{F})^{-1}\\
		&=\Phi^\mathcal{F} (\mathcal{R}_{12}^\mathcal{F})^{-1}\mathcal{F}_{21}(\Delta^\mathrm{op}\otimes\mathrm{id})(\mathcal{F})\Phi_{213}^{-1}\chi_{13}\Phi_{213}(\Delta^\mathrm{op}\otimes\mathrm{id})(\mathcal{F}^{-1})\mathcal{F}^{-1}_{21}\mathcal{R}^\mathcal{F}_{12}(\Phi^\mathcal{F})^{-1}+\Phi^\mathcal{F} \chi^\mathcal{F}_{12}(\Phi^\mathcal{F})^{-1}\\
		&=\Phi^\mathcal{F} (\mathcal{R}_{12}^\mathcal{F})^{-1}\mathcal{F}_{21}(\Delta^\mathrm{op}\otimes\mathrm{id})(\mathcal{F})\Phi_{213}^{-1}\chi_{13}(\overline{\mathcal{F}}_{i(1)}\otimes\overline{\mathcal{F}}^i\otimes\overline{\mathcal{F}}_{i(2)})\mathcal{F}^{-1}_{13}\Phi^\mathcal{F}_{213}\mathcal{R}^\mathcal{F}_{12}(\Phi^\mathcal{F})^{-1}+\Phi^\mathcal{F}\chi^\mathcal{F}_{12}(\Phi^\mathcal{F})^{-1}\\
		&=\Phi^\mathcal{F} (\mathcal{R}_{12}^\mathcal{F})^{-1}\mathcal{F}_{21}(\Delta^\mathrm{op}\otimes\mathrm{id})(\mathcal{F})\Phi_{213}^{-1}(\overline{\mathcal{F}}_{i(1)}\otimes\overline{\mathcal{F}}^i\otimes\overline{\mathcal{F}}_{i(2)})\chi_{13}\mathcal{F}^{-1}_{13}\Phi^\mathcal{F}_{213}\mathcal{R}^\mathcal{F}_{12}(\Phi^\mathcal{F})^{-1}+\Phi^\mathcal{F} \chi^\mathcal{F}_{12}(\Phi^\mathcal{F})^{-1}\\
		&=\Phi^\mathcal{F} (\mathcal{R}_{12}^\mathcal{F})^{-1}\mathcal{F}_{21}(\Delta^\mathrm{op}\otimes\mathrm{id})(\mathcal{F})\Phi_{213}^{-1}(\overline{\mathcal{F}}_{i(1)}\otimes\overline{\mathcal{F}}^i\otimes\overline{\mathcal{F}}_{i(2)})\mathcal{F}^{-1}_{13}\chi^\mathcal{F}_{13}\Phi^\mathcal{F}_{213}\mathcal{R}^\mathcal{F}_{12}(\Phi^\mathcal{F})^{-1} +\Phi^\mathcal{F} \chi^\mathcal{F}_{12}(\Phi^\mathcal{F})^{-1}\\
		&=\Phi^\mathcal{F}(\mathcal{R}_{12}^\mathcal{F})^{-1}(\Phi_{213}^\mathcal{F})^{-1}\chi^\mathcal{F}_{13}\Phi^\mathcal{F}_{213}\mathcal{R}^\mathcal{F}_{12}(\Phi^\mathcal{F})^{-1}
		+\Phi^\mathcal{F} \chi^\mathcal{F}_{12}(\Phi^\mathcal{F})^{-1}
	\end{align*}
	where we wrote
        $\mathcal{F}^{-1}=\overline{\mathcal{F}}^i\otimes\overline{\mathcal{F}}_i$.
        This means that \eqref{eq:pC2} is satisfied for
        $(H_\mathcal{F},\mathcal{R}_\mathcal{F},\chi_\mathcal{F})$. In
        complete analogy we prove \eqref{eq:pC3} for
        $(H_\mathcal{F},\mathcal{R}_\mathcal{F},\chi_\mathcal{F})$.
\end{proof}
\noindent
Building on the previous twisting result and Example \ref{ex:E(n)} we
construct the following class of pre-Cartier quasi-bialgebras.
\begin{example}\label{ex:E(n)chi}
  Consider the Hopf algebras $E(n)$ of Example \ref{ex:E(n)}. As shown in
  \cite{BeaDasGru}, there is an exhaustive class
  of quasitriangular structures $\mathcal{R}_{(a_{ij})}$ where
  $(a_{ij})$ is an arbitrary $n\times n$-matrix with values in
  $\Bbbk$. Moreover, in \cite[Theorem 3.30]{BRS24} it is shown that
  \begin{equation}
    \chi_{(b_{k\ell})}:=\sum\nolimits_{k,\ell=1}^nb_{k\ell}gx_k\otimes x_\ell\in E(n)\otimes E(n),
  \end{equation}
  where $(b_{k\ell})\in M_n(\Bbbk)$, is an exhaustive family of
  infinitesimal $\mathcal{R}$-matrices for
  $(E(n),\mathcal{R}_{(a_{ij})})$. Moreover, we have that
  $(E(n),\mathcal{R}_{(a_{ij})},\chi_{(b_{k\ell})})$ is Cartier if and
  only if $(b_{k\ell})$ is skew-symmetric, i.e. $b_{\ell
    k}=-b_{k\ell}$ for all $1\leq k,\ell\leq n$.

\noindent
  We have seen in Example \ref{ex:E(n)} that $\mathcal{F}=1\otimes g$
  is a gauge transformation, leading to the
  quasitriangular quasi-bialgebras
  $(E(n),\Delta_\mathcal{F},\varepsilon,\Phi_\mathcal{F},(\mathcal{R}_{(a_{ij})})_\mathcal{F})$
  with twisted structures explicitly given in \eqref{twist:delta} and
  \eqref{twist:R}. Together with Proposition \ref{prop:twist} we
  conclude that the quasitriangular quasi-bialgebra
  $(E(n),\Delta_\mathcal{F},\varepsilon,\Phi_\mathcal{F},(\mathcal{R}_{(a_{ij})})_\mathcal{F})$
  is pre-Cartier with respect to
  \begin{equation}\label{twistedchi}
    \begin{split}
      (\chi_{(b_{ij})})_\mathcal{F}
      &=\mathcal{F}\chi_{(b_{ij})}\mathcal{F}^{-1}
      =(1\otimes g)\sum\nolimits_{k,\ell=1}^nb_{k\ell}(gx_k\otimes x_\ell)(1\otimes g)
      =-\sum\nolimits_{k,\ell=1}^nb_{k\ell}gx_k\otimes x_\ell\\
      &=\chi_{(-b_{ij})}.
    \end{split}
  \end{equation}
  Since $(b_{k\ell})\in M_n(\Bbbk)$ is arbitrary, this shows that the
  untwisted $\chi_{(b_{ij})}$ are infinitesimal $\mathcal{R}$-matrices
  for the quasitriangular quasi-bialgebras
  $(E(n),\Delta_\mathcal{F},\varepsilon,\Phi_\mathcal{F},(\mathcal{R}_{(a_{ij})})_\mathcal{F})$.
\end{example}
\noindent
We continue by discussing morphisms of pre-Cartier
quasi-bialgebras. As expected, they are morphisms of the underlying
quasitriangular quasi-bialgebra structures, respecting the pre-Cartier
structure. Moreover, an equivalence of pre-Cartier quasi-bialgebras is
an isomorphism up to gauge equivalence.

\begin{definition}
  Let $(H, \Phi, \mathcal{R}, \chi)$ and $(H', \Phi', \mathcal{R}', \chi')$ be two pre-Cartier quasi-bialgebras.
  \begin{itemize}
  \item[i).] A \textbf{morphism of pre-Cartier quasi-bialgebras} is a
    morphism of quasitriangular quasi-bialgebras $\alpha: H \to H'$
    such that $(\alpha \otimes \alpha)(\chi) = \chi'$.
  \item[ii).] $H$ and $H'$ are \textbf{equivalent} if there is a gauge
    transformation $\mathcal{F} \in H' \otimes H'$ and an isomorphism
    of pre-Cartier quasi-bialgebras $\alpha : H \to H'_\mathcal{F}$.
  \end{itemize}
\end{definition}
\noindent
It is easy to see that equivalent pre-Cartier quasi-bialgebras have
equivalent categories of representations:

\begin{theorem}
	Let $H,H'$ be two equivalent pre-Cartier quasi-bialgebras
        together with a gauge transformation $\mathcal{F} \in H'
        \otimes H'$ and an isomorphism $\alpha: H \to
        H'_\mathcal{F}$. Then the triple $(\alpha^* , F^2_\mathcal{F},
        \mathrm{id})$ is an infinitesimally braided monoidal
        equivalence from ${}_{H'_\mathcal{F}}\mathcal{M}$ to
        ${}_H\mathcal{M}$.
\end{theorem}
\begin{proof}
  The proof consists in mimicking the same steps made in
  \cite[XV.3.4--3.7]{Kassel95} for the case of quasi-bialgebras.
\end{proof}

\subsection{Infinitesimal braid relations in pre-Cartier quasi-bialgebras}\label{sec:3.3}

Let $(H, \mathcal{R}, \chi)$ be a pre-Cartier
bialgebra. Then it is shown in \cite[Theorem 2.21]{ABSW} that
\begin{equation}
\label{eq:chi13}
\mathcal{R}^{-1}_{23} \chi_{13} \mathcal{R}_{23} =  \mathcal{R}^{-1}_{12} \chi_{13} \mathcal{R}_{12} =: \overline{\chi_{13}}.
\end{equation}
We show that any pre-Cartier bialgebra
contains elements satisfying the infinitesimal braid relations
\eqref{eq:infbraid1}. \\ First, recall the following lemma, whose
proof is straightforward.
\begin{lemma}
  \label{lemma-inf-braid-rel}
  Let $A,B,C$ be three elements of an algebra satisfying
  $[A+B,C]=0=[A+C,B]$. Then $[B+C,A]=0$.
\end{lemma}
\noindent
We can now prove the following
\begin{theorem}
  \label{theorem-inf-braid-rel-in-qtpcb}
  Let $(H, \mathcal{R}, \chi)$ be a pre-Cartier
  bialgebra. Then the following identites hold:
  \begin{align}
    [\chi_{12}, \chi_{23} + \overline{\chi_{13}}]&=0 \label{eq:new-identity-one}\\
    [\chi_{23}, \chi_{12} + \overline{\chi_{13}}]&=0 \label{eq:new-identity-two}\\
    [\chi_{12} + \chi_{23} ,  \overline{\chi_{13}}]&=0. \label{eq:new-identity-three} 
  \end{align}
\end{theorem}
\begin{proof}
  We first prove \eqref{eq:new-identity-one} and \eqref{eq:new-identity-two}:
  \begin{equation*}
    \begin{split}
      \chi_{12} (\chi_{23} + \mathcal{R}_{23}^{-1} \chi_{13} \mathcal{R}_{23})  &\stackrel{\eqref{eq:pC3}}{=} \chi_{12} (\Delta \ten \id)(\chi) 
      \stackrel{\eqref{eq:pC1}}{=} (\Delta \ten \id)(\chi)\chi_{12} \stackrel{\eqref{eq:pC3}}{=}(\chi_{23} + \mathcal{R}_{23}^{-1} \chi_{13} \mathcal{R}_{23}) \chi_{12}\\
      \chi_{23} (\chi_{12} + \mathcal{R}_{12}^{-1} \chi_{13} \mathcal{R}_{12})  &\stackrel{\eqref{eq:pC2}}{=} \chi_{23} (\id \ten \Delta)(\chi) 
      \stackrel{\eqref{eq:pC1}}{=} (\id \ten \Delta)(\chi)\chi_{23} \stackrel{\eqref{eq:pC2}}{=}(\chi_{12} + \mathcal{R}_{12}^{-1} \chi_{13} \mathcal{R}_{12}) \chi_{23}.
    \end{split}
  \end{equation*}
  Finally, identity \eqref{eq:new-identity-three} follows by Lemma \ref{lemma-inf-braid-rel}.
\end{proof}
\noindent
Note that Theorem \ref{theorem-inf-braid-rel-in-qtpcb} in particular
implies that $[\chi_{12},\chi_{23}] = 0$ if and only if
$[\chi_{12},\overline{\chi_{13}} ]= 0$ if and only if $[\chi_{23},
  \overline{\chi_{13}}]=0$.  Next, we generalize Theorem
\ref{theorem-inf-braid-rel-in-qtpcb} to the case of a pre-Cartier
quasi-bialgebra. First, we generalize Equation \eqref{eq:chi13} with
the following
\begin{proposition}
  Let $(H, \Phi, \mathcal{R},\chi)$ be a pre-Cartier
  quasi-bialgebra. Then the following identity holds
  \begin{equation}
    \label{eq:pre-cartier-chi13}
    \Phi^{-1} \mathcal{R}_{23}^{-1} \Phi_{132} \chi_{13} \Phi^{-1}_{132} \mathcal{R}_{23} \Phi = \mathcal{R}_{12}^{-1} \Phi_{213}^{-1} \chi_{13} \Phi_{213} \mathcal{R}_{12}.
  \end{equation}
\end{proposition}
\begin{proof}
  We have, using Equations \eqref{eq:pC2} and ~\eqref{eq:qtqb1}:
  \begin{align*}
    \mathcal{R}_{23}\big( \Phi (\chi_{12} + \mathcal{R}_{12}^{-1} \Phi^{-1}_{213} \chi_{13} \Phi_{213} \mathcal{R}_{12})\Phi^{-1}\big)  &= \mathcal{R}_{23} \big( (\id \ten \Delta)(\chi)\big)\\
    &= \big( (\id \ten \Delta^{\mathrm{op}}(\chi))\big) \mathcal{R}_{23} \\
    &= \big((\id \ten \tau) (\id \ten \Delta)(\chi)\big)\mathcal{R}_{23} \\
    &= (\id \ten \tau)\big( \Phi (\chi_{12} + \mathcal{R}_{12}^{-1} \Phi^{-1}_{213} \chi_{13} \Phi_{213} \mathcal{R}_{12})\Phi^{-1} \big)\mathcal{R}_{23}\\
    &= \big(\Phi_{132} (\chi_{13} + \mathcal{R}^{-1}_{13} \Phi^{-1}_{231}\chi_{12} \Phi_{231} \mathcal{R}_{13}) \Phi_{132}^{-1} \big) \mathcal{R}_{23},
  \end{align*}
  i.e. the following identity holds
  \begin{equation}
    \label{eq-preliminary-t13-quasi1}
    \mathcal{R}_{23} \Phi\chi_{12}\Phi^{-1} + \mathcal{R}_{23} \Phi\mathcal{R}_{12}^{-1} \Phi^{-1}_{213} \chi_{13} \Phi_{213} \mathcal{R}_{12}\Phi^{-1} = \Phi_{132}\chi_{13}\Phi_{132}^{-1}\mathcal{R}_{23} + \Phi_{132}\mathcal{R}^{-1}_{13} \Phi^{-1}_{231}\chi_{12} \Phi_{231} \mathcal{R}_{13} \Phi_{132}^{-1}\mathcal{R}_{23}.
  \end{equation}
  Next, using Equations \eqref{eq:qtqb3} and \eqref{eq:pC1} we note that
  \begin{equation}
    \label{eq-preliminary-t13-quasi2}
    \chi_{12} \Phi_{231} \mathcal{R}_{13} \Phi_{132}^{-1}\mathcal{R}_{23} = \chi_{12} \big( (\Delta \ten \id)(\mathcal{R})\big) \Phi^{-1}= \big( (\Delta \ten \id)(\mathcal{R})\big)\chi_{12}\Phi^{-1} = \Phi_{231} \mathcal{R}_{13} \Phi_{132}^{-1}\mathcal{R}_{23}\Phi\chi_{12}\Phi^{-1}.
  \end{equation}
  Inserting the right hand side of Equation
  \eqref{eq-preliminary-t13-quasi2} in the last summand of Equation
  \eqref{eq-preliminary-t13-quasi1} proves the claim.
\end{proof}
\noindent
We are now ready to construct two sets of elements inside $H^{\ten 3}$
satisfying the infinitesimal braid relations \eqref{eq:infbraid1}.
\begin{theorem}
  Let $(H, \Phi, \mathcal{R},\chi)$ be a pre-Cartier quasi-bialgebra
  and consider the following elements of $H^{\ten 3}$
  \begin{align*}
    \Theta^{12} &= \chi_{12} \qquad \qquad \qquad \qquad \qquad \qquad \qquad \quad \ \ \Upsilon^{12} = \Phi \chi_{12} \Phi^{-1}\\
    \Theta^{23} &= \Phi^{-1} \chi_{23} \Phi \qquad \qquad \qquad \qquad \qquad \qquad \  \ \ \  \Upsilon^{23} = \chi_{23}  \\
    \Theta^{13} &= \Phi^{-1} \mathcal{R}_{23}^{-1} \Phi_{132} \chi_{13} \Phi^{-1}_{132} \mathcal{R}_{23} \Phi \qquad \qquad \qquad  \Upsilon^{13} = \Phi \mathcal{R}_{12}^{-1} \Phi_{213}^{-1} \chi_{13} \Phi_{213} \mathcal{R}_{12} \Phi^{-1}.
  \end{align*}
  Then the sets $\{ \Theta^{12}, \Theta^{23}, \Theta^{13}\}$ and $\{
  \Upsilon^{12}, \Upsilon^{23}, \Upsilon^{13}\}$ satisfy the
  infinitesimal braid relations \eqref{eq:infbraid1}.
\end{theorem}
\begin{proof}
  We have, using Equations \eqref{eq:pC1}, \eqref{eq:pC3}:
  \begin{align*}
    \Theta^{12} (\Theta^{23} + \Theta^{13}) &= \chi_{12} (\Phi^{-1} \chi_{23} \Phi + \Phi^{-1} \mathcal{R}_{23}^{-1} \Phi_{132} \chi_{13}) = \chi_{12} (\Delta \ten \id)(\chi) \\
    &= (\Delta \ten \id)(\chi)\chi_{12}  = (\Phi^{-1} \chi_{23} \Phi + \Phi^{-1} \mathcal{R}_{23}^{-1} \Phi_{132} \chi_{13}) \chi_{12} = (\Theta^{23} + \Theta^{13})\Theta^{12} .
  \end{align*}
  Next, using Equation \eqref{eq:pre-cartier-chi13}, we have
  \[\Phi^{-1} \chi_{23} \Phi (\chi_{12} + \mathcal{R}_{12}^{-1} \Phi_{213}^{-1} \chi_{13} \Phi_{213} \mathcal{R}_{12} ) = (\chi_{12} + \mathcal{R}_{12}^{-1} \Phi_{213}^{-1} \chi_{13} \Phi_{213} \mathcal{R}_{12} )\Phi^{-1} \chi_{23} \Phi,\]
  which is equivalent to 
  \[ \chi_{23} \Phi (\chi_{12} + \mathcal{R}_{12}^{-1} \Phi_{213}^{-1} \chi_{13} \Phi_{213} \mathcal{R}_{12} )  \Phi^{-1}= \Phi(\chi_{12} + \mathcal{R}_{12}^{-1} \Phi_{213}^{-1} \chi_{13} \Phi_{213} \mathcal{R}_{12} )\Phi^{-1} \chi_{23} .\]
  Now, using Equations \eqref{eq:pC1}, \eqref{eq:pC2}:
  \begin{align*}
    \chi_{23} \Phi (\chi_{12} + \mathcal{R}_{12}^{-1} \Phi_{213}^{-1} \chi_{13} \Phi_{213} \mathcal{R}_{12} )  \Phi^{-1}&= \chi_{23} (\Phi\chi_{12}\Phi^{-1} + \Phi\mathcal{R}_{12}^{-1} \Phi_{213}^{-1} \chi_{13} \Phi_{213} \mathcal{R}_{12}  \Phi^{-1})\\
    &= \chi_{23} (\id \ten \Delta) (\chi) =  (\id \ten \Delta) (\chi) \chi_{23} \\
    &=  (\Phi\chi_{12}\Phi^{-1} + \Phi\mathcal{R}_{12}^{-1} \Phi_{213}^{-1} \chi_{13} \Phi_{213} \mathcal{R}_{12}  \Phi^{-1})\chi_{23}\\
    &=  \Phi(\chi_{12} + \mathcal{R}_{12}^{-1} \Phi_{213}^{-1} \chi_{13} \Phi_{213} \mathcal{R}_{12} ) \Phi^{-1}\chi_{23},
  \end{align*}
  proving that $[\Theta^{23}, \Theta^{12} + \Theta^{13}]=0$.  Finally,
  $[\Theta^{13}, \Theta^{12} + \Theta^{23}]=0$ follows from Lemma
  \ref{lemma-inf-braid-rel}.\\ In complete analogy one can prove that
  the infinitesimal braid relations \eqref{eq:infbraid1} hold for the
  set $\{ \Upsilon^{12}, \Upsilon^{23}, \Upsilon^{13}\}$.
\end{proof}
\noindent
Note that such two different sets satisfying the infinitesimal braid
relation refers to the two different ways of bracketing the threefold
tensor product $H \ten H \ten H$.

\subsection{The quasi counterpart of the infinitesimal quantum Yang-Baxter equation}\label{sec:3.4}

Recall that any quasitriangular quasi-bialgebra gives rise to a
solution of the \emph{quasi} counterpart of the quantum Yang-Baxter
equation, i.e.
\begin{equation}
\label{eq:QYBE}
\mathcal{R}_{12}\Phi_{231}\mathcal{R}_{13}\Phi^{
-1}_{132}\mathcal{R}_{23}\Phi = \Phi_{321}\mathcal{R}_{23}\Phi^{-1}_{312}\mathcal{R}_{13}\Phi_{213}\mathcal{R}_{12}.
\end{equation}
holds.  Furthermore, in \cite{ABSW} the authors show that any pre-Cartier bialgebra gives rise to a solution of the
infinitesimal quantum Yang-Baxter equation, which is the following
\[ 
\mathcal{R}_{12}\chi_{12}\mathcal{R}_{13}\mathcal{R}_{23}+ \mathcal{R}_{12}\mathcal{R}_{13}\chi_{13}\mathcal{R}_{23} + \mathcal{R}_{12}\mathcal{R}_{13}\mathcal{R}_{23}\chi_{23} = 
\mathcal{R}_{23}\chi_{23}\mathcal{R}_{13}\mathcal{R}_{12} + \mathcal{R}_{23}\mathcal{R}_{13}\chi_{13}\mathcal{R}_{12} + \mathcal{R}_{23}\mathcal{R}_{13}\mathcal{R}_{12}\chi_{12}.
\]
We generalize the latter in the following
\begin{theorem}
  Let $(H, \Phi, \mathcal{R}, \chi)$ be a pre-Cartier
  quasi-bialgebra. Then the quasi counterpart of the infinitesimal
  quantum Yang-Baxter equation
  \begin{equation}
    \label{eq:QQYBE}
    \begin{split}
      &\mathcal{R}_{12}\chi_{12}\Phi_{231}\mathcal{R}_{13}\Phi^{
        -1}_{132}\mathcal{R}_{23}\Phi + \mathcal{R}_{12}\Phi_{231}\mathcal{R}_{13}\chi_{13}\Phi^{
        -1}_{132}\mathcal{R}_{23}\Phi + \mathcal{R}_{12}\Phi_{231}\mathcal{R}_{13}\Phi^{
        -1}_{132}\mathcal{R}_{23}\chi_{23}\Phi = \\
      &\Phi_{321}\mathcal{R}_{23}\chi_{23}\Phi^{-1}_{312}\mathcal{R}_{13}\Phi_{213}\mathcal{R}_{12} + \Phi_{321}\mathcal{R}_{23}\Phi^{-1}_{312}\mathcal{R}_{13}\chi_{13}\Phi_{213}\mathcal{R}_{12} + \Phi_{321}\mathcal{R}_{23}\Phi^{-1}_{312}\mathcal{R}_{13}\Phi_{213}\mathcal{R}_{12}\chi_{12}
    \end{split}
  \end{equation}
  holds.
\end{theorem}
\begin{proof}
  We compute, using Equations \eqref{eq:qtqb3} and \eqref{eq:pC1}:
  \begin{align*}
    \chi_{12} \Phi_{231} \mathcal{R}_{13} \Phi_{132}^{-1} \mathcal{R}_{23} \Phi &= \chi_{12} \big( (\Delta \ten \id) (\mathcal{R})\big) = \big( (\Delta \ten \id) (\mathcal{R})\big) \chi_{12} = \Phi_{231} \mathcal{R}_{13} \Phi_{132}^{-1} \mathcal{R}_{23} \Phi \chi_{12},
  \end{align*}
  and therefore, using Equation \eqref{eq:QYBE} we have 
  \begin{align*}
    \mathcal{R}_{12}  \Phi_{231} \mathcal{R}_{13} \Phi_{132}^{-1} \mathcal{R}_{23} \Phi  &= \mathcal{R}_{12} \Phi_{231} \mathcal{R}_{13} \Phi_{132}^{-1} \mathcal{R}_{23} \Phi \chi_{12} = \Phi_{321} \mathcal{R}_{23} \Phi^{-1}_{312} \mathcal{R}_{13} \Phi_{213} \mathcal{R}_{12} \chi_{12}
  \end{align*}
  i.e.  the first summand of the left hand side of \eqref{eq:QQYBE} is
  equal to the third summand of the right hand side of
  \eqref{eq:QQYBE}. Similarly, we compute, using Equations
  \eqref{eq:qtqb2} and \eqref{eq:pC1}:
  \begin{align*}
    \chi_{23} \Phi_{312}^{-1} \mathcal{R}_{13} \Phi_{213} \mathcal{R}_{12} \Phi^{-1} = \chi_{23} \big( (\id \ten \Delta)(\mathcal{R})\big) =  \big( (\id \ten \Delta)(\mathcal{R})\big)\chi_{23} = \Phi_{312}^{-1} \mathcal{R}_{13} \Phi_{213} \mathcal{R}_{12} \Phi^{-1} \chi_{23}
  \end{align*}
  and therefore, using Equation \eqref{eq:QYBE} we have
  \begin{align*}
    \Phi_{321} \mathcal{R}_{23} \chi_{23} \Phi_{312}^{-1} \mathcal{R}_{13} \Phi_{213} \mathcal{R}_{12} \Phi^{-1} = \Phi_{321} \mathcal{R}_{23}\Phi_{312}^{-1} \mathcal{R}_{13} \Phi_{213} \mathcal{R}_{12} \Phi^{-1} \chi_{23} = \mathcal{R}_{12}\Phi_{231}\mathcal{R}_{13}\Phi^{
      -1}_{132}\mathcal{R}_{23} \chi_{23},
  \end{align*}
  i.e. the third summand of the left hand side of \eqref{eq:QQYBE} is equal to the first summand of the right hand side of \eqref{eq:QQYBE}.  Finally, we have, using Equation \eqref{eq:QYBE}:
  \begin{align*}
    \mathcal{R}_{12} \Phi_{231} \mathcal{R}_{13} = \Phi_{321} \mathcal{R}_{23} \Phi^{-1}_{312} \mathcal{R}_{13} \Phi_{213} \mathcal{R}_{12} \Phi^{-1} \mathcal{R}_{23}^{-1} \Phi_{132}.
  \end{align*}
  Therefore, using Equation \eqref{eq:pre-cartier-chi13} we have 
  \begin{align*}
    \mathcal{R}_{12} \Phi_{231} \mathcal{R}_{13} \chi_{13} \Phi_{132}^{-1} \mathcal{R}_{23} \Phi &= \Phi_{321} \mathcal{R}_{23} \Phi^{-1}_{312} \mathcal{R}_{13} \Phi_{213} \mathcal{R}_{12} \Phi^{-1} \mathcal{R}_{23}^{-1} \Phi_{132} \chi_{13} \Phi_{132}^{-1} \mathcal{R}_{23} \Phi \\
    &= \Phi_{321} \mathcal{R}_{23} \Phi^{-1}_{312} \mathcal{R}_{13} \Phi_{213} \mathcal{R}_{12}\mathcal{R}_{12}^{-1} \Phi_{213}^{-1} \chi_{13} \Phi_{213} \mathcal{R}_{12}\\
    &= \Phi_{321} \mathcal{R}_{23} \Phi^{-1}_{312} \mathcal{R}_{13} \chi_{13} \Phi_{213} \mathcal{R}_{12}
  \end{align*}
  showing that the second term of the left hand side of
  \eqref{eq:QQYBE} is equal to the second term of the right hand side
  of \eqref{eq:QQYBE}, hence concluding the proof.
\end{proof}

\subsection{Cartier rings and their representations}\label{sec:3.5}

It is well-known (see e.g. \cite[X.6]{Kassel95}) that for any $n \geq
2$ and any left module $V$ over a quasitriangular bialgebra, there is
a canonical representation $\rho_n^{\mathcal{R}}: \mathsf{B}_n \to
\mathrm{End}(V^{\ten n})$ of the braid group on $n$ strands
$\mathsf{B}_n$. In this section we provide a more general result for
any left module of a pre-Cartier bialgebra. The fact
that the infinitesimal braiding of a pre-Cartier bialgebra is not
necessarily invertible, together with the fact that the corresponding
module category is linear, suggests that the correct framework to consider is the one of rings instead of groups. This reasoning
motivates the following
\begin{definition}
	\label{definition-symmetric-pre-cartier-ring}
	Let $n \geq 2$. The $n$-th \textbf{Cartier ring} is the ring
        $\mathsf{C}_n$ generated by $n-1$ invertible elements
        $\beta_1, \ldots, \beta_{n-1}$ and by $n-1$ non-invertible
        elements $\gamma_{1}, \ldots, \gamma_{n-1}$ subject to
        relations
	\begin{equation*}
		\begin{split}
			\beta_i \beta_j &= \beta_j \beta_i \qquad \text{for all } i,j \text{ such that } |i-j| \geq 2 \\
			\gamma_i \gamma_j &= \gamma_j \gamma_i \qquad \text{for all } i,j \text{ such that } |i-j| \geq 2 \\
			\beta_i \gamma_j &= \gamma_j \beta_i \qquad \text{for all } i,j \text{ such that } |i-j| \geq 2 \\
			\beta_i \beta_{i+1} \beta_i &= \beta_{i+1} \beta_i  \beta_{i+1} \qquad \text{for all } i=1, \ldots, n-2 \\
			\beta_i (\gamma_i \beta_{i+1} + \beta_{i+1} \gamma_{i+1}) \beta_i + \beta_i \beta_{i+1}\beta_i\gamma_i  &= \beta_{i+1}(\gamma_{i+1} \beta_i + \beta_i \gamma_i)\beta_{i+1} + \beta_{i+1}\beta_i \beta_{i+1} \gamma_{i+1}\qquad \text{for all } i=1, \ldots, n-2.
		\end{split}
	\end{equation*}
\end{definition}
\noindent
Note that the subgroup of $\mathsf{C}_n$ generated by $\beta_1, \ldots, \beta_{n-1}$ is the
braid group on $n$ strands.  By construction of the rings
$\mathsf{C}_n$, we easily get the following
\begin{theorem}
\label{theorem-Cartier-rings}
	Let $(H, \mathcal{R}, \chi)$ be a pre-Cartier
        bialgebra, $V \in {}_H \mathcal{M}$ and $n \geq 2$. Then we
        have a representation of the $n$-th Cartier ring
	\begin{equation*}
		\begin{split}
			\rho_{n}^{\mathcal{R}, \chi} : \mathsf{C}_n &\to \mathrm{End}(V^{\otimes n})\\
			\beta_i & \mapsto \mathcal{R}_i \\
			\gamma_i & \mapsto \chi_i
		\end{split}
	\end{equation*}
	where
	\begin{equation*}
		\begin{split}
			\mathcal{R}_i (v_1 \otimes \cdots \otimes v_n) &= v_1 \otimes \cdots \otimes v_{i-1} \otimes  \mathcal{R}^\mathrm{op} \cdot (v_{i+1} \otimes v_{i}) \otimes v_{i+2} \otimes \cdots  \otimes v_{n}\\
			\chi_i(v_1 \otimes \cdots  \otimes v_n) &= v_1 \otimes \cdots \otimes v_{i-1} \otimes \chi \cdot (v_i \otimes v_{i+1}) \otimes v_{i+2} \otimes \cdots \otimes v_n.
		\end{split}
	\end{equation*}
\end{theorem}
\noindent
Note that Theorem \ref{theorem-Cartier-rings} can be generalized to
the case of any pre-Cartier quasi-bialgebra by mimicking the reasoning
of \cite[XV.4]{Kassel95}.

\section{Deforming braided monoidal categories}\label{sec:4}

In this section we construct solutions of the infinitesimal braid
relations \eqref{eq:infbraid1} and \eqref{eq:infbraid2} for any
pre-Cartier category in order to generalize the deformation theorem of Cartier, mentioned in the introduction, to the non-symmetric case. Furthermore, we are able to obtain the corresponding deformation result for pre-Cartier categories.
However, here we have to require the additional assumption
$
[t^{12},t^{23}]=0.
$
 
 \noindent
We first introduce the notion of a deformation of a braided monoidal category, see e.g.  \cite[Section 3]{Yett}, \cite[Section 2]{Yet}, or \cite[Section 5]{Car}. For this, we restrict our attention from now on to $\Bbbk$-linear categories, i.e. categories enriched over $\mathrm{Vec}_\Bbbk$.
\begin{definition}
  Let $(\mathcal{C},\ten,\Bbbk,a,\ell,r,\sigma)$ be a
  $\Bbbk$-linear braided monoidal category. A \textbf{deformation}
  of $\mathcal{C}$ is a $\Bbbk[[\hbar]]$-linear braided monoidal
  category $(\hat{\mathcal{C}},
  \hat{\otimes},\Bbbk,\hat{a},\hat{\ell},\hat{r},\hat{\sigma})$ such
  that $\mathrm{Ob}(\hat{\mathcal{C}}) = \mathrm{Ob}(\mathcal{C})$,
  $\mathrm{Hom}_{\hat{\mathcal{C}}}(X,Y) =
  \mathrm{Hom}_{\mathcal{C}}(X,Y)[[\hbar]] $, the tensor product
  $\hat{\otimes}$ coincides with $\otimes$ on objects and it is the
  $\hbar$-adic completion of $\ten$ on morphisms, and
  \begin{align*}
    \hat{a} \mod \ \hbar &= a \\
    \hat{\ell} \mod \ \hbar &= \ell \\
    \hat{r} \mod \ \hbar &= r \\
    \hat{\sigma} \mod \ \hbar &= \sigma.
  \end{align*}
  In particular, the new concatenation rule of morphisms is given by 
  \begin{align*}
    \Big(\sum_{k=0}^\infty \phi_k\hbar^k\Big)\circ \Big(\sum_{k=0}^\infty \psi_k\hbar^k\Big)
    :=
    \sum_{k=0}^\infty\Big(\sum_{i=0}^k \phi_i\circ \psi_{k-i}\Big)\hbar^k.
  \end{align*}
\end{definition}
\noindent
Similarly, one can define a deformation of a braided monoidal functor,
but this is not needed in what follows.  
\begin{remark}
  Further note that for any deformation $\hat{\mathcal{C}}$ of $\mathcal{C}$
  we have a canonical braided monoidal functor $\lim_{\hbar\to
    0}\colon \hat{\mathcal{C}} \to \mathcal{C}$, which is identity on
  objects and
  \begin{align*}
    \lim_{\hbar\to 0}\Big(\sum_{k=0}^\infty \phi_k\hbar^k\Big)=\phi_0. 
  \end{align*}
\end{remark}
\begin{example}
  \label{example-trivial-deformation}
  The trivial deformation of a $\Bbbk$-linear monoidal category
  $\mathcal{C}$ is the one obtained by $\hbar$-linearly extending the
  isomorphisms $a,\ell,r,\sigma$. We shall denote it by
  $\tilde{\mathcal{C}}$.
\end{example}

\subsection{Infinitesimal braid relations in pre-Cartier categories}\label{sec:4.1}

We begin by proving some auxiliary lemmas. They are concerning the
presentation of the infinitesimal braiding $t$ of a pre-Cartier
category $(\mathcal{C},\otimes,\sigma,t)$ as an endomorphism of the
$n$-fold tensor product for $n>2$, particularly $n=3$ and $n=4$. For
example, given three objects $X,Y,Z$ in $\mathcal{C}$, there are (at
least) two possible ways to define $t^{13}_{X,Y,Z}\colon X\otimes
Y\otimes Z\to X\otimes Y\otimes Z$, namely
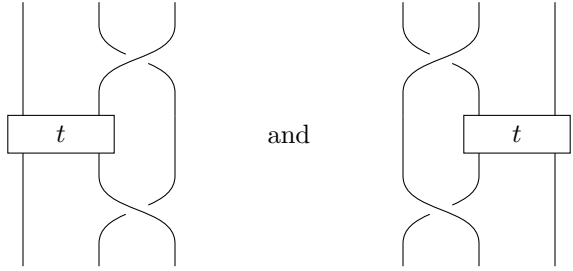
\begin{equation*}
  \begin{tikzpicture}
    \braid at (0,0) a_2^{-1} ;
    \draw [draw=black] (-0.2,-2) rectangle (1.2,-1.5);
    \node at (0.5,-1.75) {$t$};
    \draw (2,-1.5)--(2,-2);
    \braid  at (0,-2) a_2 ;
    \node at (3.5,-1.75) {\text{and}};
    \braid  at (5,0) a_1^{-1};
    \draw [draw=black] (5.8,-2) rectangle (7.2,-1.5);
    \draw (5,-1.5)--(5,-2);
    \node at (6.5,-1.75) {$t$};
    \braid  at (5,-2) a_1;
    \draw (7,0)--(7,-1.5);
    \draw (7,-2)--(7,-3.5);
  \end{tikzpicture}
\end{equation*}
where we employed the pictorial notation as in \cite[XIV]{Kassel95}, read from top to bottom, and where braidings (respectively inverse braidings) are denoted by an undercrossing (respectively overcrossings) from left to right.
In non-symmetric pre-Cartier categories there are yet two alternative
ways to describe $t^{13}$. They emerge by swapping the roles of
$\sigma$ and $\sigma^{-1}$. We show that all of these presentations
coincide in the following
\begin{lemma}
	\label{lemma-t13}
	Let $\mathcal{C}$ be a pre-Cartier category and consider the
        following morphisms
        \begin{align}
			t_{X,Y,Z}^{13, \mathrm{I}} &= (\mathrm{id}_X \otimes \sigma^{-1}_{Y,Z}) \circ (t_{X,Z} \otimes \mathrm{id}_Y) \circ (\mathrm{id}_X \otimes \sigma_{Y,Z}) \label{eq:t13-I} \\
			t_{X,Y,Z}^{13, \mathrm{II}} &= (\sigma_{X,Y}^{-1} \otimes \mathrm{id}_Z) \circ (\mathrm{id}_Y \otimes t_{X,Z}) \circ (\sigma_{X,Y} \otimes \mathrm{id}_Z) \label{eq:t13-II} \\
			t^{13, \mathrm{III}}_{X,Y,Z} &= (\sigma_{Y,X} \otimes \mathrm{id}_Z) \circ (\mathrm{id}_Y \otimes t_{X,Z}) \circ (\sigma_{Y,X}^{-1} \otimes \mathrm{id}_X)  \label{eq:t13-III} \\
			t^{13, \mathrm{IV}}_{X,Y,Z} &= (\mathrm{id}_X \otimes \sigma_{Z,Y}) \circ (t_{X,Z} \otimes \mathrm{id}_Y) \circ (\mathrm{id}_X \otimes \sigma_{Z,Y}^{-1}).   \label{eq:t13-IV} 
        \end{align}
	Then $t_{X,Y,Z}^{13, \mathrm{I}}=t_{X,Y,Z}^{13, \mathrm{II}}=
        t^{13, \mathrm{III}}_{X,Y,Z} = t^{13, \mathrm{IV}}_{X,Y,Z} =:
        t^{13}_{X,Y,Z}$.
\end{lemma}
\begin{proof}
Note that $t^{13,\mathrm{I}}=t^{13,\mathrm{III}}$ and $t^{13,\mathrm{II}}=t^{13,\mathrm{IV}}$ can be shown by solely relying on the naturality of the braiding and the hexagon equations, i.e. they hold for every morphism $X\otimes Z\to X\otimes Z$ replacing $t_{X,Z}$:
  \begin{align*}
	t^{13,\mathrm{III}}_{X,Y,Z} &= (\sigma_{Y,X} \otimes \mathrm{id}_Z)\circ (\mathrm{id}_Y \otimes t_{X,Z} )  \circ (\sigma_{Y,X}^{-1} \otimes \mathrm{id}_Z)\\
	&= (\mathrm{id}_X \otimes \sigma^{-1}_{Y,Z}) \circ (\sigma_{Y,X \otimes Z}) \circ (\mathrm{id}_Y \otimes t_{X,Z} )  \circ (\sigma_{Y,X \otimes Z}^{-1})\circ (\mathrm{id}_X \otimes \sigma_{Y,Z}) \\
	&= (\mathrm{id}_X \otimes \sigma^{-1}_{Y,Z}) \circ (\sigma_{Y,X \otimes Z}) \circ (\sigma_{Y,X \otimes Z}^{-1})  \circ (t_{X,Z} \otimes \mathrm{id}_Y) \circ (\mathrm{id}_X \otimes \sigma_{Y,Z}) \\
	&= (\mathrm{id}_X \otimes \sigma^{-1}_{Y,Z}) \circ (t_{X,Z} \otimes \mathrm{id}_Y) \circ (\mathrm{id}_X \otimes \sigma_{Y,Z}) \\
	&=t^{13, \mathrm{I}}_{X,Y,Z},
\end{align*}
and similarly for $t^{13, \mathrm{II}}_{X,Y,Z}=t^{13,\mathrm{IV}}_{X,Y,Z}$.
Using now the naturality of $t$, as well as the fact that it is an infinitesimal braiding, i.e., Equation \eqref{eq:pre-cartier-one}, we have:
  \begin{equation}
    \label{proof-aux-one}
    \begin{split}
      &(\mathrm{id}_X \otimes \sigma_{Y,Z}) \circ \Big(( t_{X,Y} \otimes \mathrm{id}_Z) +(\sigma_{X,Y}^{-1} \otimes \mathrm{id}_Z) \circ (\mathrm{id}_Y \otimes t_{X,Z}) \circ (\sigma_{X,Y} \otimes \mathrm{id}_Z) \Big) \\
      &=(\mathrm{id}_X \otimes \sigma_{Y,Z}) \circ (t_{X, Y \otimes Z}) = (t_{X, Z \otimes Y}) \circ (\mathrm{id}_X \otimes \sigma_{Y,Z})  \\
      &= \Big((t_{X,Z} \ten \mathrm{id}_{Y} ) +(\sigma_{X,Z}^{-1} \otimes \mathrm{id}_Y) \circ (\mathrm{id}_Z \ten t_{X,Y})  \circ (\sigma_{X,Z} \ten \mathrm{id}_Y) \Big) \circ (\mathrm{id}_X \otimes \sigma_{Y,Z}).
    \end{split}
  \end{equation}
  Moreover, using the naturality of $\sigma$ and Equation
  \eqref{eq:hexagon-strict-two} we have:
  \begin{equation}
    \label{proof-aux-two}
    \begin{split}
      (\mathrm{id}_Z \otimes t_{X,Y}) &\circ (\sigma_{X,Z} \otimes \mathrm{id}_Y) \circ (\mathrm{id}_X \otimes \sigma_{Y,Z}) = (\mathrm{id}_Z \otimes t_{X,Y}) \circ \sigma_{X \otimes Y,Z} \\
      &= \sigma_{X \otimes Y,Z} \circ (t_{X,Y} \otimes \mathrm{id}_Z) = (\sigma_{X,Z} \otimes \mathrm{id}_Y) \circ (\mathrm{id}_X \otimes \sigma_{Y,Z}) \circ (t_{X,Y} \otimes \mathrm{id}_Z).
    \end{split}
  \end{equation}
  Finally, we have 
  \begin{align*}
    t_{X,Y,Z}^{13, \mathrm{I}} &= (\mathrm{id}_X \otimes \sigma^{-1}_{Y,Z}) \circ (t_{X,Z} \otimes \mathrm{id}_Y) \circ (\mathrm{id}_X \otimes \sigma_{Y,Z}) \\
    & \stackrel{\eqref{proof-aux-one}}{=} (\mathrm{id}_X \otimes \sigma_{Y,Z}^{-1})  (\mathrm{id}_X \otimes \sigma_{Y,Z})  (t_{X,Y} \otimes \mathrm{id}_Z) +\\
    &+ (\mathrm{id}_{X} \otimes \sigma^{-1}_{Y,Z}) (\mathrm{id}_X \otimes \sigma_{Y,Z}) (\sigma^{-1}_{X,Y} \otimes \mathrm{id}_Z)(\mathrm{id}_Y \otimes t_{X,Z}) (\sigma_{X,Y} \otimes \mathrm{id}_Z) \\
    &- (\mathrm{id}_X \otimes \sigma^{-1}_{Y,Z}) (\sigma_{X,Z}^{-1} \otimes \mathrm{id}_Y) (\mathrm{id}_Z \otimes t_{X,Y}) (\sigma_{X,Z} \otimes \mathrm{id}_Y) (\mathrm{id}_X \otimes \sigma_{Y,Z}) \\
    &=(t_{X,Y} \otimes \mathrm{id}_Z)  + (\sigma^{-1}_{X,Y} \otimes \mathrm{id}_Z)(\mathrm{id}_Y \otimes t_{X,Z}) (\sigma_{X,Y} \otimes \mathrm{id}_Z) - (\sigma_{X \otimes Y,Z}^{-1}) (\mathrm{id}_Z \otimes t_{X,Y}) (\sigma_{X \otimes Y,Z}) \\
    & \stackrel{\eqref{proof-aux-two}} {=} (t_{X,Y} \otimes \mathrm{id}_Z) + t^{13, \mathrm{II}}_{X,Y,Z} - (t_{X,Y} \otimes \mathrm{id}_Z) \\
    &= t^{13, \mathrm{II}}_{X,Y,Z},
  \end{align*}
which concludes the proof.
\end{proof}
\begin{remark}
  In complete analogy one can prove that the non-strict counterpart of
  the morphisms \eqref{eq:t13-I}-\eqref{eq:t13-IV}
  \begin{align*}
    t_{X,Y,Z}^{13, \mathrm{I}} &=a^{-1}_{X,Y,Z} \circ (\mathrm{id}_X \otimes \sigma^{-1}_{Y,Z})\circ a_{X,Z,Y} \circ (t_{X,Z} \otimes \mathrm{id}_Y) \circ a^{-1}_{X,Z,Y} \circ (\mathrm{id}_X \otimes \sigma_{Y,Z}) \circ a_{X,Y,Z} \\
    t_{X,Y,Z}^{13, \mathrm{II}} &= (\sigma_{X,Y}^{-1} \otimes \mathrm{id}_Z) \circ a^{-1}_{Y,X,Z} \circ (\mathrm{id}_Y \otimes t_{X,Z}) \circ a_{Y,X,Z} \circ (\sigma_{X,Y} \otimes \mathrm{id}_Z)\\
    t^{13, \mathrm{III}}_{X,Y,Z} &= (\sigma_{Y,X} \otimes \mathrm{id}_Z) \circ a^{-1}_{Y,X,Z} \circ (\mathrm{id}_Y \otimes t_{X,Z})\circ a_{Y,X,Z} \circ (\sigma_{Y,X}^{-1} \otimes \mathrm{id}_X) \\
    t^{13, \mathrm{IV}}_{X,Y,Z} &= a^{-1}_{X,Y,Z} \circ (\mathrm{id}_X \otimes \sigma_{Z,Y}) \circ a_{X,Z,Y} \circ (t_{X,Z} \otimes \mathrm{id}_Y) \circ a^{-1}_{X,Z,Y} \circ (\mathrm{id}_X \otimes \sigma_{Z,Y}^{-1})\circ a_{X,Y,Z} 
  \end{align*}
  are equal.
\end{remark}
\noindent
The previous Lemma and Remark, together with Example \ref{examples-precartier-categories} ii.), imply the following result.
\begin{corollary}
Let $(\mathcal{C},\otimes,\sigma,t)$ be a pre-Cartier category. Then for all $k\in\mathbb{Z}$ the braided monoidal category $(\mathcal{C},\otimes,\sigma^{\pm1})$ is again pre-Cartier with infinitesimal braiding $t^{(k)}:=\sigma^{k}\circ t\circ\sigma^{-k}$. If $(\mathcal{C},\otimes,\sigma,t)$ is Cartier, it follows that $t^{(k)}=t$ for all $k\in\mathbb{Z}$.
\end{corollary}
\noindent
Note also that, if we define $t^{24}_{X,Y,Z,W} \coloneqq \id_X \ten
t^{13}_{Y,Z,W}$, we obtain that the corresponding ways to express
$t^{24}$ are equivalent as well.\\ Next, the situation becomes
seemingly more demanding if we consider the $4$-fold tensor product
$X\otimes Y\otimes Z\otimes W$ in a pre-Cartier category
$(\mathcal{C},\otimes,\sigma,t)$. In particular, the (at least) four
ways to construct $t^{14}_{X,Y,Z,W}$
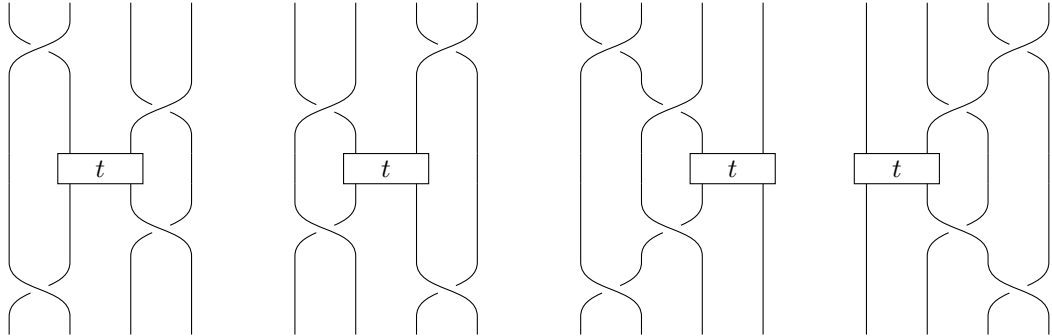
\begin{equation*}
  \begin{tikzpicture}[scale=0.8]
    \braid at (0,0) a_1^{-1} a_3^{-1};
    \draw [draw=black] (0.8,-3) rectangle (2.2,-2.5);
    \node at (1.5,-2.75) {$t$};
    \draw (0,-2.5)--(0,-3);
    \draw (3,-2.5)--(3,-3);
    \braid at (0,-3) a_3 a_1;
    \begin{scope}[shift={(0.7,0)}]
    \braid at (4,0) a_3^{-1} a_1^{-1};
    \draw [draw=black] (4.8,-3) rectangle (6.2,-2.5);
    \node at (5.5,-2.75) {$t$};
    \draw (4,-2.5)--(4,-3);
    \draw (7,-2.5)--(7,-3);
    \braid at (4,-3) a_1 a_3;
    \begin{scope}[shift={(0.7,0)}]
    \braid at (8,0) a_1^{-1} a_2^{-1};
    \draw [draw=black] (9.8,-3) rectangle (11.2,-2.5);
    \node at (10.5,-2.75) {$t$};
    \draw (8,-2.5)--(8,-3);
    \draw (9,-2.5)--(9,-3);
    \braid at (8,-3) a_2 a_1;
    \draw (11,0)--(11,-2.5);
        \draw (11,-3)--(11,-5.5);
        \begin{scope}[shift={(0.7,0)}]
    \braid at (12,0) a_3^{-1} a_2^{-1};
    \draw [draw=black] (11.8,-3) rectangle (13.2,-2.5);
    \node at (12.5,-2.75) {$t$};
    \draw (14,-2.5)--(14,-3);
    \draw (15,-2.5)--(15,-3);
    \braid at (12,-3) a_2 a_3;
\end{scope}
\end{scope}
    \end{scope}
  \end{tikzpicture}
\end{equation*}
do not seem equivalent at first glance. Nonetheless, we prove that
they are in the following
\begin{lemma}
\label{lemma-t14}
	Consider the morphisms
	\begin{equation*}
		\begin{split}
			t^{14,\mathrm{I}}_{X,Y,Z,W} &= (\sigma^{-1}_{X,Y} \otimes \mathrm{id}_{Z \otimes W})(\mathrm{id}_{Y \otimes X} \otimes \sigma^{-1}_{Z,W}) (\mathrm{id}_Y \otimes t_{X,W} \otimes \mathrm{id}_Z)(\mathrm{id}_{Y \otimes X} \otimes \sigma_{Z,W}) (\sigma_{X,Y} \otimes \mathrm{id}_{Z \otimes W})\\
			t^{14,\mathrm{II}}_{X,Y,Z,W} &= (\mathrm{id}_{X \otimes Y} \otimes \sigma^{-1}_{Z,W}) (\sigma^{-1}_{X,Y} \otimes \mathrm{id}_{W \otimes Z}) (\mathrm{id}_Y \otimes t_{X,W} \otimes \mathrm{id}_Z)(\sigma_{X,Y} \otimes \mathrm{id}_{W \otimes Z})(\mathrm{id}_{X \otimes Y} \otimes \sigma_{Z,W})  \\
			t^{14,\mathrm{III}}_{X,Y,Z,W} &= (\sigma^{-1}_{X,Y} \otimes \mathrm{id}_{Z \otimes W})(\mathrm{id}_Y \otimes \sigma^{-1}_{X,Z} \otimes \mathrm{id}_W)(\mathrm{id}_{Y \otimes Z} \otimes t_{X,W})(\mathrm{id}_Y \otimes \sigma_{X,Z} \otimes \mathrm{id}_W) (\sigma_{X,Y} \otimes \mathrm{id}_{Z \otimes W})\\
			t^{14,\mathrm{IV}}_{X,Y,Z,W} &= (\mathrm{id}_{X \otimes Y} \otimes \sigma^{-1}_{Z,W}) (\mathrm{id}_X \otimes \sigma^{-1}_{Y,W} \otimes \mathrm{id}_Z)(t_{X,W} \otimes \mathrm{id}_{Y \otimes Z})(\mathrm{id}_X \otimes \sigma_{Y,W} \otimes \mathrm{id}_Z) (\mathrm{id}_{X \otimes Y} \otimes \sigma_{Z,W}).
		\end{split}
	\end{equation*}
	Then $t^{14,\mathrm{I}}_{X,Y,Z,W} = t^{14,\mathrm{II}}_{X,Y,Z,W} = t^{14,\mathrm{III}}_{X,Y,Z,W} = t^{14,\mathrm{IV}}_{X,Y,Z,W} =: t^{14}_{X,Y,Z,W}.$
\end{lemma}
\begin{proof}
  The fact that $t_{X,Y,Z,W}^{14, \mathrm{I}}=t_{X,Y,Z,W}^{14,
    \mathrm{II}}$ is straightforward. Next, using Lemma
  \ref{lemma-t13} we have
  \begin{equation*}
    \begin{split}
      t_{X,Y,Z,W}^{14, \mathrm{I}}&= (\sigma^{-1}_{X,Y} \otimes \mathrm{id}_{Z \otimes W}) (\mathrm{id}_Y \otimes t^{13, \mathrm{I}}_{X,Z,W}) (\sigma_{X,Y} \otimes \mathrm{id}_{Z \otimes W})\\
      &= (\sigma^{-1}_{X,Y} \otimes \mathrm{id}_{Z \otimes W}) (\mathrm{id}_Y \otimes t^{13, \mathrm{II}}_{X,Z,W}) (\sigma_{X,Y} \otimes \mathrm{id}_{Z \otimes W})=  t_{X,Y,Z,W}^{14, \mathrm{III}} 
    \end{split}
  \end{equation*}	
  and
  \begin{equation*}
    \begin{split}		
      t_{X,Y,Z,W}^{14, \mathrm{II}}&= (\mathrm{id}_{X \otimes Y} \otimes \sigma^{-1}_{Z,W}) (t^{13, \mathrm{II}}_{Y,X,W} \otimes \mathrm{id}_Z) (\mathrm{id}_{X \otimes Y} \otimes \sigma_{Z,W})\\
      &= (\mathrm{id}_{X \otimes Y} \otimes \sigma^{-1}_{Z,W}) (t^{13, \mathrm{I}}_{Y,X,W} \otimes \mathrm{id}_Z) (\mathrm{id}_{X \otimes Y} \otimes \sigma_{Z,W})=  t_{X,Y,Z,W}^{14, \mathrm{IV}}. 
    \end{split}
  \end{equation*}
\end{proof}
\noindent
We are now prepared to prove that any pre-Cartier category provides
solutions of the infinitesimal braid relations. This will be the main
tool to attack the deformation theorems of the subsequent
sections. In what follows, we shall use Lemmas
\ref{lemma-t13}-\ref{lemma-t14} in order to represent morphisms
$t^{13},t^{24},t^{14}$ in the most convenient form.
\begin{proposition}
  \label{proposition-inf-braid-rel-in-pre-cartier-category}
  Let $(\mathcal{C}, \otimes, I,a,\sigma,t)$ be a pre-Cartier
  category. For any objects $X,Y,Z,W$ in $\mathcal{C}$, consider the
  following endomorphims of $X \otimes Y \otimes Z$
  \begin{align}
    t^{12}_{X,Y,Z} &\coloneqq t_{X,Y} \otimes \mathrm{id}_{Z} \label{eq:inf-braid-rel-one} \\
    t^{23}_{X,Y,Z} &\coloneqq \mathrm{id}_X \otimes t_{Y,Z} \label{eq:inf-braid-rel-two} \\
    t^{13}_{X,Y,Z} &\coloneqq (\mathrm{id}_X \otimes \sigma_{Y,Z}^{-1}) \circ (t_{X,Z} \otimes \mathrm{id}_{Y}) \circ (\mathrm{id}_X \otimes \sigma_{Y,Z})\label{eq:inf-braid-rel-three} 
  \end{align}
  and the following endomorphisms of $X \otimes Y \otimes Z \otimes W$
  \begin{align}
    t^{12}_{X,Y,Z,W} &\coloneqq t^{12}_{X,Y,Z}  \otimes \mathrm{id}_{W} = t_{X,Y} \otimes \mathrm{id}_{Z \otimes W} \label{eq:inf-braid-rel-four}  \\
    t^{23}_{X,Y,Z,W} &\coloneqq t^{23}_{X,Y,Z} \otimes \mathrm{id}_W = \mathrm{id}_X \otimes t_{Y,Z} \otimes \mathrm{id}_W \label{eq:inf-braid-rel-five} \\
    t^{13}_{X,Y,Z,W} &\coloneqq t^{13}_{X,Y,Z} \otimes \mathrm{id}_W = \Big( (\mathrm{id}_X \otimes \sigma_{Y,Z}^{-1}) \circ (t_{X,Z} \otimes \mathrm{id}_{Y}) \circ (\mathrm{id}_X \otimes \sigma_{Y,Z}) \Big) \otimes \mathrm{id}_W \label{eq:inf-braid-rel-six} \\
    t^{14}_{X,Y,Z,W} &\coloneqq(\mathrm{id}_X \otimes \sigma_{Y \otimes Z,W}^{-1}) \circ (t_{X,W} \otimes \mathrm{id}_{Y \otimes Z}) \circ (\mathrm{id}_X \otimes \sigma_{Y \otimes Z,W}) \label{eq:inf-braid-rel-seven}  \\
    t^{24}_{X,Y,Z,W} &\coloneqq \mathrm{id}_X \otimes t^{13}_{Y,Z,W} = \mathrm{id}_X \otimes \Big( (\mathrm{id}_Y \otimes \sigma_{Z,W}^{-1}) \circ (t_{Y,W} \otimes \mathrm{id}_{Z}) \circ (\mathrm{id}_Y \otimes \sigma_{Z,W})\Big) \label{eq:inf-braid-rel-eight}  \\
    t^{34}_{X,Y,Z,W} &\coloneqq \mathrm{id}_X \otimes t^{23}_{Y,Z,W} = \mathrm{id}_{X \otimes Y} \otimes t_{Z,W} \label{eq:inf-braid-rel-nine} 
  \end{align}
  Then the morphisms $t^{ij}_{X,Y,Z}$ and $t^{ij}_{X,Y,Z,W}$ satisfy
  the infinitesimal braid relations
  \eqref{eq:infbraid1}-\eqref{eq:infbraid2}.
\end{proposition}
\begin{proof}
  We first show the infinitesimal braid relations \eqref{eq:infbraid1}
  for the morphisms
  \eqref{eq:inf-braid-rel-one}-\eqref{eq:inf-braid-rel-three}.
  \begin{itemize}
  \item $\mathbf{[t^{12}_{X,Y,Z}, t^{23}_{X,Y,Z} + t^{13}_{X,Y,Z}]=0}$. 
    We have, using the naturality of $t$ and Equation \eqref{eq:pre-cartier-two}:
    \begin{equation*}
      \begin{split}
	t^{12}_{X,Y,Z} \circ (t^{23}_{X,Y,Z} + t^{13}_{X,Y,Z}) &= (t_{X,Y} \otimes \mathrm{id}_Z) \circ \Big( (\mathrm{id}_X \otimes t_{Y,Z}) + (\mathrm{id}_{X} \otimes \sigma^{-1}_{Y,Z}) \circ (t_{X,Z} \otimes \mathrm{id}_Y) \circ (\mathrm{id}_X \otimes \sigma_{Y,Z})\Big) \\
	&=(t_{X,Y} \otimes \mathrm{id}_Z) \circ t_{X \otimes Y,Z} =t_{X \otimes Y,Z}  \circ (t_{X,Y} \otimes \mathrm{id}_Z)  =  (t^{23}_{X,Y,Z} + t^{13}_{X,Y,Z}) \circ t^{12}_{X,Y,Z}.
      \end{split}
    \end{equation*}
  \item $\mathbf{[t^{23}_{X,Y,Z}, t^{12}_{X,Y,Z} +
      t^{13}_{X,Y,Z}]=0}$.  We have, using the naturality of $t$ and
    Equation \eqref{eq:pre-cartier-one}:
    \begin{equation*}
      \begin{split}
	t^{23}_{X,Y,Z} \circ (t^{12}_{X,Y,Z} + t^{13}_{X,Y,Z}) &= (\mathrm{id}_X \otimes t_{Y,Z}) \circ \Big( (t_{X,Y} \otimes \mathrm{id}_Z) + (\sigma_{X,Y}^{-1} \otimes \mathrm{id}_Z)(\mathrm{id}_Y \otimes t_{X,Z})(\sigma_{X,Y} \otimes \mathrm{id}_Z) \Big) \\
	&= (\mathrm{id}_X \otimes t_{Y,Z}) \circ t_{X,Y \otimes Z} = t_{X, Y \otimes Z} \circ  (\mathrm{id}_X \otimes t_{Y,Z}) = (t^{12}_{X,Y,Z} + t^{13}_{X,Y,Z}) \circ t^{23}_{X,Y,Z}.
      \end{split}
    \end{equation*}
  \item $\mathbf{[t^{13}_{X,Y,Z}, t^{12}_{X,Y,Z} +
      t^{23}_{X,Y,Z}]=0}$. Follows by Lemma \ref{lemma-inf-braid-rel}.
  \end{itemize}
  Next, we prove the infinitesimal braid relations \eqref{eq:infbraid1}
  for the morphisms
  \eqref{eq:inf-braid-rel-four}-\eqref{eq:inf-braid-rel-nine}. Note
  that the proof for $\{ i,j,k\}=\{1,2,3\}$ and $\{i,j,k \}=\{2,3,4\}$
  follows from the previous computations, namely:
  \begin{align*}
    [t^{12}_{X,Y,Z,W}, t^{23}_{X,Y,Z,W} + t^{13}_{X,Y,Z,W}]&=0 = [t^{23}_{X,Y,Z,W}, t^{34}_{X,Y,Z,W} + t^{24}_{X,Y,Z,W}] \\
    [t^{23}_{X,Y,Z,W}, t^{12}_{X,Y,Z,W} + t^{13}_{X,Y,Z,W}]&=0 = [t^{34}_{X,Y,Z,W}, t^{23}_{X,Y,Z,W} + t^{24}_{X,Y,Z,W}]\\
    [t^{13}_{X,Y,Z,W}, t^{23}_{X,Y,Z,W} + t^{12}_{X,Y,Z,W}]&=0 = [t^{24}_{X,Y,Z,W}, t^{34}_{X,Y,Z,W} + t^{23}_{X,Y,Z,W}] .
  \end{align*}
  We now prove the cases $\{ i,j,k\}=\{1,2,4\}$ and $\{i,j,k
  \}=\{1,3,4\}$.
  \begin{itemize} 
  \item $\mathbf{[t^{12}_{X,Y,Z,W}, t^{24}_{X,Y,Z,W} +
      t^{14}_{X,Y,Z,W}]=0}$.  We have, using the naturality of $t$ and
    Equations \eqref{eq:hexagon-strict-two},
    \eqref{eq:pre-cartier-two}:
    \begin{align*}
      t^{12}_{X,Y,Z,W} &\circ (t^{24}_{X,Y,Z,W} + t^{14}_{X,Y,Z,W})= (t_{X,Y} \otimes \mathrm{id}_{Z \otimes W}) \Bigg( (\mathrm{id}_{X \otimes Y} \otimes \sigma^{-1}_{Z,W}) (\mathrm{id}_X \otimes t_{Y,W} \otimes \mathrm{id}_Z)(\mathrm{id}_{X \otimes Y} \otimes \sigma_{Z,W}) +\\
      & \qquad \qquad \qquad \qquad + (\mathrm{id}_X \otimes \sigma^{-1}_{Y \otimes Z,W})(t_{X,W} \otimes \mathrm{id}_{Y \otimes Z})(\mathrm{id}_X \otimes \sigma_{Y \otimes Z,W})\Bigg) \\
      &= (t_{X,Y} \otimes \mathrm{id}_{Z \otimes W}) \Bigg( (\mathrm{id}_{X \otimes Y} \otimes \sigma^{-1}_{Z,W}) (\mathrm{id}_X \otimes t_{Y,W} \otimes \mathrm{id}_Z)(\mathrm{id}_{X \otimes Y} \otimes \sigma_{Z,W}) + \\
      &  \qquad \qquad + (\mathrm{id}_{X \otimes Y} \otimes \sigma_{Z,W}^{-1})(\mathrm{id}_X \otimes \sigma^{-1}_{Y,W} \otimes \mathrm{id}_Z)(t_{X,W} \otimes \mathrm{id}_{Y \otimes Z})(\mathrm{id}_X \otimes \sigma_{Y, W} \otimes \mathrm{id}_Z) (\mathrm{id}_{X \otimes Y} \otimes \sigma_{Z,W})\Bigg)\\
      &= (t_{X,Y} \otimes \mathrm{id}_{Z \otimes W})(\mathrm{id}_{X \otimes Y} \otimes \sigma_{Z,W}^{-1}) \Bigg( (\mathrm{id}_X \otimes t_{Y,W} \otimes \mathrm{id}_Z) +\\
      &\qquad \qquad \qquad \qquad + (\mathrm{id}_X \otimes \sigma^{-1}_{Y,W} \otimes \mathrm{id}_Z)(t_{X,W} \otimes \mathrm{id}_{Y \otimes Z})(\mathrm{id}_X \otimes \sigma_{Y, W} \otimes \mathrm{id}_Z) \Bigg) (\mathrm{id}_{X \otimes Y} \otimes \sigma_{Z,W}) \\
      &=(t_{X,Y} \otimes \mathrm{id}_{Z \otimes W})(\mathrm{id}_{X \otimes Y} \otimes \sigma_{Z,W}^{-1}) (t_{X \otimes Y,W} \otimes \mathrm{id}_Z)  (\mathrm{id}_{X \otimes Y} \otimes \sigma_{Z,W}) \\
      &= (\mathrm{id}_{X \otimes Y} \otimes \sigma^{-1}_{Z,W}) (t_{X,Y} \otimes \mathrm{id}_{W \otimes Z}) (t_{X \otimes Y,W} \otimes \mathrm{id}_Z)  (\mathrm{id}_{X \otimes Y} \otimes \sigma_{Z,W}) \\
      &=  (\mathrm{id}_{X \otimes Y} \otimes \sigma^{-1}_{Z,W})(t_{X \otimes Y,W} \otimes \mathrm{id}_Z)(t_{X,Y} \otimes \mathrm{id}_{W \otimes Z}) (\mathrm{id}_{X \otimes Y} \otimes \sigma_{Z,W}) \\
      &=  (\mathrm{id}_{X \otimes Y} \otimes \sigma^{-1}_{Z,W})(t_{X \otimes Y,W} \otimes \mathrm{id}_Z)(\mathrm{id}_{X \otimes Y} \otimes \sigma_{Z,W})(t_{X,Y} \otimes \mathrm{id}_{Z \otimes W}) \\
      &=  (\mathrm{id}_{X \otimes Y} \otimes \sigma^{-1}_{Z,W}) \Bigg((\mathrm{id}_X \otimes t_{Y,W} \otimes \mathrm{id}_Z) +\\
      &\qquad \qquad + (\mathrm{id}_X \otimes \sigma^{-1}_{Y,W} \otimes \mathrm{id}_Z)(t_{X,W} \otimes \mathrm{id}_{Y \otimes Z})(\mathrm{id}_X \otimes \sigma_{Y, W} \otimes \mathrm{id}_Z)  \Bigg)(\mathrm{id}_{X \otimes Y} \otimes \sigma_{Z,W})( t^{12}_{X,Y,Z,W})\\
      &= (t^{24}_{X,Y,Z,W} + t^{14}_{X,Y,Z,W})\circ t^{12}_{X,Y,Z,W}.
    \end{align*}
  \item $\mathbf{[t^{24}_{X,Y,Z,W}, t^{12}_{X,Y,Z,W} +
      t^{14}_{X,Y,Z,W}]=0}$.  First, we note that
    \[ t^{12}_{X,Y,Z,W} = (\mathrm{id}_{X \otimes Y} \otimes \sigma^{-1}_{Z,W}) (t_{X,Y} \otimes \mathrm{id}_{W \otimes Z}) (\mathrm{id}_{X \otimes Y} \otimes \sigma_{Z,W}).\]
    We have, using the naturality of $t$ and Equation \eqref{eq:pre-cartier-one}:
    \begin{align*}
      &t^{24}_{X,Y,Z,W} \circ (t^{12}_{X,Y,Z,W} + t^{14}_{X,Y,Z,W}) = t^{24}_{X,Y,Z,W} \circ (\mathrm{id}_{X \otimes Y} \otimes \sigma^{-1}_{Z,W}) \circ \Bigg( (t_{X,Y} \otimes \mathrm{id}_{W \otimes Z}) + \\
      & \qquad +(\sigma_{X,Y}^{-1} \otimes \mathrm{id}_{W \otimes Z}) (\mathrm{id}_Y \otimes t_{X,W} \otimes \mathrm{id}_Z) (\sigma_{X,Y} \otimes \mathrm{id}_{W \otimes Z}) \Bigg) \circ  (\mathrm{id}_{X \otimes Y} \otimes \sigma_{Z,W}) \\
      &=t^{24}_{X,Y,Z,W} \circ (\mathrm{id}_{X \otimes Y} \otimes \sigma^{-1}_{Z,W}) \circ (t_{X, Y \otimes W} \otimes \mathrm{id}_Z) \circ (\mathrm{id}_{X \otimes Y} \otimes \sigma_{Z,W}) \\
      &= (\mathrm{id}_{X \otimes Y} \otimes \sigma^{-1}_{Z,W}) (\mathrm{id}_X \otimes t_{Y,W} \otimes \mathrm{id}_Z)(\mathrm{id}_{X \otimes Y} \otimes \sigma_{Z,W}) (\mathrm{id}_{X \otimes Y} \otimes \sigma^{-1}_{Z,W})  (t_{X, Y \otimes W} \otimes \mathrm{id}_Z)  (\mathrm{id}_{X \otimes Y} \otimes \sigma_{Z,W}) \\
      &=  (\mathrm{id}_{X \otimes Y} \otimes \sigma^{-1}_{Z,W}) (\mathrm{id}_X \otimes t_{Y,W} \otimes \mathrm{id}_Z)  (t_{X, Y \otimes W} \otimes \mathrm{id}_Z)  (\mathrm{id}_{X \otimes Y} \otimes \sigma_{Z,W}) \\
      &= (\mathrm{id}_{X \otimes Y} \otimes \sigma^{-1}_{Z,W}) (t_{X, Y \otimes W} \otimes \mathrm{id}_Z)  (\mathrm{id}_X \otimes t_{Y,W} \otimes \mathrm{id}_Z)  (\mathrm{id}_{X \otimes Y} \otimes \sigma_{Z,W}) \\
      &= (\mathrm{id}_{X \otimes Y} \otimes \sigma^{-1}_{Z,W}) (t_{X, Y \otimes W} \otimes \mathrm{id}_Z) (\mathrm{id}_{X \otimes Y} \otimes \sigma_{Z,W}) (\mathrm{id}_{X \otimes Y} \otimes \sigma^{-1}_{Z,W}) (\mathrm{id}_X \otimes t_{Y,W} \otimes \mathrm{id}_Z)  (\mathrm{id}_{X \otimes Y} \otimes \sigma_{Z,W}) \\
      &= (t^{12}_{X,Y,Z,W} + t^{14}_{X,Y,Z,W} ) \circ t^{24}_{X,Y,Z,W}.
    \end{align*}
  \item $\mathbf{[t^{14}_{X,Y,Z,W}, t^{12}_{X,Y,Z,W} +
      t^{24}_{X,Y,Z,W}]=0}$. Follows by Lemma
    \ref{lemma-inf-braid-rel}.
  \item $\mathbf{[t^{34}_{X,Y,Z,W}, t^{13}_{X,Y,Z,W} +
      t^{14}_{X,Y,Z,W}]=0}$. We have, using the naturality of $t$ and
    Equation \eqref{eq:pre-cartier-one}:
    \begin{align*}
      &t^{34}_{X,Y,Z,W} \circ (t^{13}_{X,Y,Z,W} + t^{14}_{X,Y,Z,W}) = \\
      &= (\mathrm{id}_{X \otimes Y} \otimes t_{Z,W}) \bigg( (\sigma_{X,Y}^{-1} \otimes \mathrm{id}_{ Z \otimes W}) (\mathrm{id}_Y \otimes t_{X,Z} \otimes \mathrm{id}_W)(\sigma_{X,Y} \otimes \mathrm{id}_{Z \otimes W}) + \\
      &+ (\sigma_{X,Y}^{-1} \otimes \mathrm{id}_{Z \otimes W})(\mathrm{id}_Y \otimes \sigma^{-1}_{X,Z} \otimes \mathrm{id}_W) (\mathrm{id}_{Y \otimes Z} \otimes t_{X,W})(\mathrm{id}_Y \otimes \sigma_{X,Z} \otimes \mathrm{id}_W) (\sigma_{X,Y} \otimes \mathrm{id}_{Z \otimes W}) \bigg)\\
      &=(\mathrm{id}_{X \otimes Y} \otimes t_{Z,W}) (\sigma_{X,Y}^{-1} \otimes \mathrm{id}_{Z \otimes W}) \bigg( (\mathrm{id}_Y \otimes t_{X,Z} \otimes \mathrm{id}_W) + \\
      &+ (\mathrm{id}_Y \otimes \sigma_{X,Z}^{-1} \otimes \mathrm{id}_W) (\mathrm{id}_{Y \otimes Z} \otimes t_{X,W})(\mathrm{id}_Y \otimes \sigma_{X,Z} \otimes \mathrm{id}_W) \bigg) (\sigma_{X,Y} \otimes \mathrm{id}_{Z \otimes W})\\ 
      &= (\sigma_{X,Y}^{-1} \otimes \mathrm{id}_{Z \otimes W}) (\mathrm{id}_{Y \otimes X} \otimes t_{Z,W})(	\mathrm{id}_Y \otimes t_{X,Z \otimes W})(\sigma_{X,Y} \otimes \mathrm{id}_{Z \otimes W}) \\
      &= (\sigma_{X,Y}^{-1} \otimes \mathrm{id}_{Z \otimes W}) (	\mathrm{id}_Y \otimes t_{X,Z \otimes W})(\mathrm{id}_{Y \otimes X} \otimes t_{Z,W})(\sigma_{X,Y} \otimes \mathrm{id}_{Z \otimes W})  \\
      &=(\sigma_{X,Y}^{-1} \otimes \mathrm{id}_{Z \otimes W}) (	\mathrm{id}_Y \otimes t_{X,Z \otimes W})(\sigma_{X,Y} \otimes \mathrm{id}_{Z \otimes W})(\mathrm{id}_{X \otimes Y} \otimes t_{Z,W}) \\
      &= (t^{13}_{X,Y,Z,W} + t^{14}_{X,Y,Z,W}) \circ t^{34}_{X,Y,Z,W}
    \end{align*}
  \item $\mathbf{[t^{13}_{X,Y,Z,W}, t^{14}_{X,Y,Z,W} + t^{34}_{X,Y,Z,W}]=0}$. First, we note that 
    \begin{align*}
      t^{34}_{X,Y,Z,W} = (\sigma_{X,Y}^{-1} \otimes \mathrm{id}_{Z \otimes W})(\mathrm{id}_{Y \otimes X} \otimes t_{Z,W}) (\sigma_{X,Y} \otimes \mathrm{id}_{Z \otimes W}).
    \end{align*}
    We have, using the naturality of $t$ and Equation
    \eqref{eq:pre-cartier-two}:
    \begin{align*}
      &t^{13}_{X,Y,Z,W} \circ (t^{34}_{X,Y,Z,W} + t^{14}_{X,Y,Z,W}) = t^{13}_{X,Y,Z,W} \circ (\sigma^{-1}_{X,Y} \otimes \mathrm{id}_{Z \otimes W}) \circ \Bigg( (\mathrm{id}_{Y \otimes X} \otimes t_{Z,W}) + \\
      & \qquad +(\mathrm{id}_{Y \otimes X} \otimes \sigma^{-1}_{Z,W})(\mathrm{id}_Y \otimes t_{X,W} \otimes \mathrm{id}_Z)(\mathrm{id}_{Y \otimes X} \otimes \sigma_{Z,W}) \Bigg) \circ  (\sigma_{X,Y} \otimes \mathrm{id}_{Z \otimes W}) \\
      &=t^{13}_{X,Y,Z,W} \circ (\sigma^{-1}_{X,Y} \otimes \mathrm{id}_{Z \otimes W}) \circ (\mathrm{id}_Y \otimes t_{X \otimes Z,W}) \circ(\sigma_{X,Y} \otimes \mathrm{id}_{Z \otimes W}) \\
      &=(	\sigma^{-1}_{X,Y} \otimes \mathrm{id}_{Z \otimes W})(\mathrm{id}_Y \otimes t_{X,Z} \otimes \mathrm{id}_W)(	\sigma_{X,Y} \otimes \mathrm{id}_{Z \otimes W}) (\sigma^{-1}_{X,Y} \otimes \mathrm{id}_{Z \otimes W}) (\mathrm{id}_Y \otimes t_{X \otimes Z,W})(\sigma_{X,Y} \otimes \mathrm{id}_{Z \otimes W}) \\
      &= (	\sigma^{-1}_{X,Y} \otimes \mathrm{id}_{Z \otimes W})(\mathrm{id}_Y \otimes t_{X,Z} \otimes \mathrm{id}_W)(	\sigma_{X,Y} \otimes \mathrm{id}_{Z \otimes W}) (\sigma^{-1}_{X,Y} \otimes \mathrm{id}_{Z \otimes W}) (\mathrm{id}_Y \otimes t_{X \otimes Z,W})(\sigma_{X,Y} \otimes \mathrm{id}_{Z \otimes W}) \\
      &= (	\sigma^{-1}_{X,Y} \otimes \mathrm{id}_{Z \otimes W})(\mathrm{id}_Y \otimes t_{X,Z} \otimes \mathrm{id}_W) (\mathrm{id}_Y \otimes t_{X \otimes Z,W})(\sigma_{X,Y} \otimes \mathrm{id}_{Z \otimes W})\\
      &= (	\sigma^{-1}_{X,Y} \otimes \mathrm{id}_{Z \otimes W})(\mathrm{id}_Y \otimes t_{X \otimes Z,W})(\mathrm{id}_Y \otimes t_{X,Z} \otimes \mathrm{id}_W) (\sigma_{X,Y} \otimes \mathrm{id}_{Z \otimes W})\\
      &= (	\sigma^{-1}_{X,Y} \otimes \mathrm{id}_{Z \otimes W})(\mathrm{id}_Y \otimes t_{X \otimes Z,W})(\sigma_{X,Y} \otimes \mathrm{id}_{Z \otimes W})(\sigma_{X,Y}^{-1} \otimes \mathrm{id}_{Z \otimes W})(\mathrm{id}_Y \otimes t_{X,Z} \otimes \mathrm{id}_W) (\sigma_{X,Y} \otimes \mathrm{id}_{Z \otimes W})\\
      &=(t^{34}_{X,Y,Z,W} + t^{14}_{X,Y,Z,W})  \circ t^{13}_{X,Y,Z,W} .
    \end{align*}
  \item  $\mathbf{[t^{14}_{X,Y,Z,W}, t^{13}_{X,Y,Z,W} + t^{34}_{X,Y,Z,W}]=0}$. Follows by Lemma \ref{lemma-inf-braid-rel}.
  \end{itemize}
  Finally, we prove the infinitesimal braid relations
  \eqref{eq:infbraid2}.
  \begin{itemize}
  \item $\mathbf{[t^{12}_{X,Y,Z,W},t^{34}_{X,Y,Z,W}]=0}$. Straightforward.
  \item $\mathbf{[t^{14}_{X,Y,Z,W}, t^{23}_{X,Y,Z,W}]=0}$.  We have, using the naturality of the braiding: 
    \begin{align*}
      t^{14}_{X,Y,Z,W} \circ  t^{23}_{X,Y,Z,W} &= (\mathrm{id}_X \otimes \sigma_{Y \otimes Z, W})^{-1} \circ (t_{X,W} \otimes \mathrm{id}_{Y \otimes Z}) \circ (\mathrm{id}_X \otimes \sigma_{Y \otimes Z,W}) \circ (\mathrm{id}_X \otimes t_{Y,Z} \otimes \mathrm{id}_W) \\
      &=(\mathrm{id}_X \otimes \sigma_{Y \otimes Z, W})^{-1} \circ (t_{X,W} \otimes \mathrm{id}_{Y \otimes Z}) \circ (\mathrm{id}_{X \otimes W} \otimes t_{Y,Z}) \circ (\mathrm{id}_X \otimes \sigma_{Y \otimes Z,W}) \\
      &= (\mathrm{id}_X \otimes \sigma_{Y \otimes Z, W})^{-1} \circ (\mathrm{id}_{X \otimes W} \otimes t_{Y,Z}) \circ  (t_{X,W} \otimes \mathrm{id}_{Y \otimes Z})  \circ (\mathrm{id}_X \otimes \sigma_{Y \otimes Z,W}) \\
      &= (\mathrm{id}_X \otimes t_{Y,Z} \otimes \mathrm{id}_W) \circ (\mathrm{id}_X \otimes \sigma_{Y \otimes Z, W})^{-1} \circ  (t_{X,W} \otimes \mathrm{id}_{Y \otimes Z})  \circ (\mathrm{id}_X \otimes \sigma_{Y \otimes Z,W}) \\
      &= t^{23}_{X,Y,Z,W} \circ t^{14}_{X,Y,Z,W}
    \end{align*}
  \item $\mathbf{[t^{13}_{X,Y,Z,W}, t^{24}_{X,Y,Z,W}]=0}$. We have:
    \begin{align*}
      &t^{13}_{X,Y,Z,W} \circ t^{24}_{X,Y,Z,W} 
      = \\
      &=(\mathrm{id}_X \otimes \sigma^{-1}_{Y,Z} \otimes \mathrm{id}_W)(t_{X,Z} \otimes \mathrm{id}_{Y \otimes W})( \mathrm{id}_X \otimes \sigma_{Y,Z} \otimes \mathrm{id}_W) (\mathrm{id}_X \otimes \sigma^{-1}_{Y,Z} \otimes \mathrm{id}_W)(\mathrm{id}_{X \otimes Z} \otimes t_{Y,W})(\mathrm{id}_X \otimes \sigma_{Y,Z} \otimes \mathrm{id}_W)\\
      &= (\mathrm{id}_X \otimes \sigma^{-1}_{Y,Z} \otimes \mathrm{id}_W) (t_{X,Z} \otimes t_{Y,W})(\mathrm{id}_X \otimes \sigma_{Y,Z} \otimes \mathrm{id}_W)\\
      &=(\mathrm{id}_X \otimes \sigma^{-1}_{Y,Z} \otimes \mathrm{id}_W)(\mathrm{id}_{X \otimes Z} \otimes t_{Y,W})(\mathrm{id}_X \otimes \sigma_{Y,Z} \otimes \mathrm{id}_W) (\mathrm{id}_X \otimes \sigma^{-1}_{Y,Z} \otimes \mathrm{id}_W)(t_{X,Z} \otimes \mathrm{id}_{Y \otimes W})( \mathrm{id}_X \otimes \sigma_{Y,Z} \otimes \mathrm{id}_W)\\
      &= t^{24}_{X,Y,Z,W} \circ t^{13}_{X,Y,Z,W} .
    \end{align*}
  \end{itemize}
\end{proof}

\subsection{Deformation of Cartier categories}\label{sec:CartierCase}

In this section we restrict our attention to Cartier categories. Let
us stress that these are braided and in general not symmetric, while
the infinitesimal braiding is compatible with the braiding in the
sense of \eqref{eq:cartier-category}. We show that, in this
generality, the underlying braided monoidal structure can be deformed
via the infinitesimal braiding into another braided monoidal
structure. This vastly generalizes the result in \cite[Section 5]{Car}, which
required a symmetric braiding. 
\begin{theorem}\label{thm:main}
  Suppose that $\mathcal{C}$ is a $\Bbbk$-linear Cartier category with infinitesimal
  braiding $t$, and let $\Psi$ be a Drinfeld associator. Then
  $\mathcal{C}_{\Psi,t} = (\mathcal{C}_{\Psi,t}, \tilde{\otimes}
  , \Bbbk, a^\Psi , \tilde{\ell}, \tilde{r}, \sigma^\Psi)$ is a braided
  monoidal category, where
  \begin{equation}\label{deformedstructure}
    \begin{split}
      a_{X,Y,Z}^\Psi &= \tilde{a}_{X,Y,Z} \circ 	\Psi\big(\hbar t_{X,Y} \otimes \mathrm{id}_Z,\hbar a^{-1}_{X,Y,Z} \circ (\mathrm{id}_X \otimes t_{Y,Z})  \circ a_{X,Y,Z}\big) \\
      \sigma_{X,Y}^\Psi &= \tilde{\sigma}_{X,Y} \circ e^{\frac{\hbar}{2} t_{X,Y}}
    \end{split}
  \end{equation}
are the deformed associativity constraint and braiding, respectively.
\end{theorem}
\begin{proof}
  Note first that $a^\Psi , \sigma^\Psi, \tilde{\ell}, \tilde{r}$ are
  natural isomorphisms since the $\hbar$-linear extentsion of an
  isomorphism is again an isomorphism, and composing with an
  exponential or a Drinfeld associator gives again an invertible
  morphism. Next, we show that $\mathcal{C}_{\Psi,t}$ satisfies the
  axioms of a braided monoidal category
  \eqref{eq:triangle-axiom},\eqref{eq:pentagon-axiom},\eqref{eq:hexagon-axiom-one},\eqref{eq:hexagon-axiom-two}. We
  shall omit the tilda symbol over morphisms and tensors for brevity.
  \begin{itemize}
  \item \textbf{The triangle axiom \eqref{eq:triangle-axiom}:}
    Recalling that $t_{X,\Bbbk} = t_{\Bbbk, Y}=0$ (see
    Remark \ref{remark-precartier-categories}) and that $\Psi(0,0) =
    1$ (see Remark \ref{remark-associator-A-B-commute}), we have
    \begin{equation*}
      \begin{split}
        (\id_X \ten \ell_Y) \circ a^{\Psi}_{X,\Bbbk, Y} &=  (\id_X \ten \ell_Y) \circ a_{X,Y,Z} \circ 	\Psi(\hbar t_{X,\Bbbk} \ten \id_Y,\hbar a_{X,\Bbbk,Y} \circ (\id_X \ten t_{\Bbbk,Y})  \circ a^{-1}_{X,\Bbbk,Y}) \\
        &= (\id_X \ten \ell_Y) \circ a_{X,Y,Z} \circ \id_{(X \ten\Bbbk) \ten Y} \\
        &=  (\id_X \ten \ell_Y) \circ a_{X,Y,Z} \\
        &= r_X \ten 	\id_Y
      \end{split}
    \end{equation*}
    where the last equality follows from the fact that
    $\tilde{\mathcal{C}}$ is monoidal (see Example
    \ref{example-trivial-deformation}).
  \item \textbf{The pentagon axiom \eqref{eq:pentagon-axiom}:} We need
    to show that
    \[ (\id_X \ten a^\Psi_{Y,Z,W}) \circ a^\Psi_{X, Y \ten Z,W} \circ (a^\Psi_{X,Y,Z} \ten \id_W) = a^{\Psi}_{X,Y,Z \ten W} \circ a^\Psi_{X \ten Y,Z,W}.\]
    Denoting by $t^{ij}= t^{ij}_{X,Y,Z,W}$ and recalling that they
    satisfy the infinitesimal braid relations (see Proposition
    \ref{proposition-inf-braid-rel-in-pre-cartier-category}) we have
    \begin{align*}
      \id_X \ten a^\Psi_{Y,Z,W} &= \id_X \ten \Psi(\hbar t_{Y,Z} \ten \id_W, \hbar \id_Y \ten t_{Z,W}) = \Psi(\hbar t^{23},\hbar t^{34}) \\
      a^\Psi_{X, Y \ten Z,W} &= \Psi(\hbar t_{X, Y \ten Z} \ten \id_W, \hbar \id_X \ten t_{Y \ten Z,W}) = \Psi(\hbar t^{12} + \hbar t^{13}, \hbar t^{24} + \hbar t^{34}) \\
      a^\Psi_{X,Y,Z} \ten \id_W &= \Psi(\hbar t_{X,Y} \ten \id_Z, \hbar\id_X \ten t_{Y,Z}) \ten \id_W = \Psi(\hbar t^{12},\hbar t^{23}) \\
      a^{\Psi}_{X,Y,Z \ten W} &= \Psi(\hbar t_{X, Y} \ten \id_{Z \ten W}, \hbar\id_{X} \ten t_{Y, Z \ten W}) = \Psi(\hbar t^{12},\hbar t^{23} +\hbar t^{24}) \\
      a^\Psi_{X \ten Y,Z,W} &= \Psi(\hbar t_{X \ten Y,Z} \ten W, \hbar \id_{X \ten Y} \ten t_{Z,W}) = \Psi(\hbar t^{13} + \hbar t^{23}, \hbar t^{34}).
    \end{align*}
    Therefore, the pentagon axiom holds since $\Psi$ satisfies the
    pentagon equation \eqref{eq:pentagon-equation}.
  \item \textbf{The hexagon axiom \eqref{eq:hexagon-axiom-one}:} We
    need to show that
    \begin{equation}
      \label{eq:first-hexagon-proof}
      \begin{split}
        &\Psi(\hbar t^{12}_{Y,Z,X},\hbar t^{23}_{Y,Z,X}) \circ (\tilde{\sigma}_{X, Y \otimes Z}) \circ e^{\frac{\hbar}{2}t_{X, Y \otimes Z}} \circ \Psi(\hbar t^{12}_{X,Y,Z}, \hbar t^{23}_{X,Y,Z}) \\
        & 	\qquad \qquad=  (\mathrm{id}_Y \otimes \tilde{\sigma}_{X,Z}) \circ (\mathrm{id}_Y \otimes e^{\frac{\hbar}{2} t_{X,Z}})\circ \Psi(\hbar t^{12}_{Y,X,Z}, \hbar t^{23}_{Y,X,Z} ) \circ (\tilde{\sigma}_{X,Y} \otimes \mathrm{id}_Y)\circ ( e^{\frac{\hbar}{2} t_{X,Y}} \otimes \mathrm{id}_Y).
      \end{split}
    \end{equation}
    By the naturality of the braiding we have that $t^{12}_{Y,Z,X}
    \circ \sigma_{X,Y \ten Z} = \sigma_{X,Y \ten Z} \circ
    t^{23}_{X,Y,Z}$. Also, using Equations
    \eqref{eq:hexagon-strict-one}, \eqref{eq:cartier-category} we have
    \begin{equation*}
      \begin{split}
	\sigma_{X, Y \otimes Z}^{-1} \circ  t^{23}_{Y,Z,X} \circ \sigma_{X, Y \otimes Z} &=  (\sigma_{X,Y}^{-1} \otimes  \mathrm{id}_Z) \circ (\mathrm{id}_Y \otimes \sigma_{X,Z}^{-1}) \circ (\mathrm{id}_Y \otimes t_{Z,X}) \circ (\mathrm{id}_Y \otimes \sigma_{X,Z}) \circ (\sigma_{X,Y} \otimes  \mathrm{id}_Z) \\
	&= (\sigma_{X,Y}^{-1} \otimes \mathrm{id}_Z) \circ (\mathrm{id}_Y \otimes t_{X,Z}) \circ (\sigma_{X,Y} \otimes \mathrm{id}_Y) =  t^{13}_{X,Y,Z},
      \end{split}
    \end{equation*}
    i.e. $t^{23}_{Y,Z,X} \circ \sigma_{X, Y \otimes Z} = \sigma_{X, Y
      \otimes Z} \circ t^{13}_{X,Y,Z}$.  Therefore, we may rewrite the
    left hand side of \eqref{eq:first-hexagon-proof} as
    \[   \tilde{\sigma}_{X, Y \otimes Z}  \circ \Psi(\hbar t^{23}_{X,Y,Z},\hbar t^{13}_{X,Y,Z}) \circ  e^{\frac{\hbar}{2} (t^{12}_{X,Y,Z} + t^{13}_{X,Y,Z})} \circ  \Psi(\hbar t^{12}_{X,Y,Z},\hbar t^{23}_{X,Y,Z}). \]
    On the other side, using Equation \eqref{eq:cartier-category} we
    have $t^{12}_{Y,X,Z} \circ (\sigma_{X,Y} \ten \id_Z) =
    (\sigma_{X,Y} \ten \id_Z) \circ t^{12}_{X,Y,Z}$, and by noting
    that $t^{23}_{Y,X,Z} \circ (\sigma_{X,Y} \otimes \mathrm{id}_Z) =
    (\sigma_{X,Y} \otimes \mathrm{id}_Z) \circ t^{13}_{X,Y,Z}$ we may
    rewrite the right hand side of \eqref{eq:first-hexagon-proof} as
    \begin{equation*}
      \tilde{\sigma}_{X , Y \otimes Z}\circ e^{\frac{\hbar}{2} t^{13}_{X,Y,Z}}\circ \Psi(\hbar t^{12}_{X,Y,Z},\hbar  t^{13}_{X,Y,Z})\circ  e^{\frac{\hbar}{2} t^{12}_{X,Y,Z}}. 
    \end{equation*}
    Therefore, Equation \eqref{eq:first-hexagon-proof} is equivalent
    to the following identity
    \begin{equation}
      \label{eq:first-hexagon-proof-new}
      \begin{split}
	\Psi(\hbar t^{23}_{X,Y,Z},\hbar t^{13}_{X,Y,Z}) \circ e^{\frac{\hbar}{2} (t^{12}_{X,Y,Z} + t^{13}_{X,Y,Z})} \circ \Psi(\hbar t^{12}_{X,Y,Z}, \hbar t^{23}_{X,Y,Z}) =e^{\frac{\hbar}{2} t^{13}_{X,Y,Z}}\circ\Psi(\hbar t^{12}_{X,Y,Z},\hbar t^{13}_{X,Y,Z})\circ e^{\frac{\hbar}{2} t^{12}_{X,Y,Z}}.
      \end{split}
    \end{equation}
    Setting $A=t^{23}_{X,Y,Z}, B=t^{12}_{X,Y,Z}, C=t^{13}_{X,Y,Z} ,
    \Lambda= A+B+C$ and using Equation \eqref{eq:inverse-associator},
    we have that \eqref{eq:first-hexagon-proof-new} is equivalent to
    \[ e^{\frac{\hbar}{2}(A+B+C)} = e^{\frac{\hbar}{2} A}\Psi(\hbar C,\hbar A)e^{\frac{\hbar}{2} C} \Psi(\hbar B,\hbar C) e^{\frac{\hbar}{2} B} \Psi(\hbar A,\hbar B) \]
    which holds since $A,B,C$ satisfy the infinitesimal braid
    relations \eqref{eq:infbraid1} (see Proposition
    \ref{proposition-inf-braid-rel-in-pre-cartier-category}) and
    $\Psi$ satisfies the hexagon equation \eqref{eq:hexagon-equation}.
  \item \textbf{The hexagon axiom \eqref{eq:hexagon-axiom-two}:}	We need to show that
    \begin{equation}
      \label{eq:second-hexagon-proof}
      \begin{split}
        &\Psi(\hbar t^{12}_{Z,X,Y}, \hbar t^{23}_{Z,X,Y})^{-1}\circ \tilde{\sigma}_{X \otimes Y,Z}\circ e^{\frac{\hbar}{2} t_{X \otimes Y,Z}}\circ \Psi(\hbar t^{12}_{X,Y,Z},\hbar t^{23}_{X,Y,Z})^{-1}\\
        & \qquad \qquad= (\tilde{\sigma}_{X,Z} \otimes \mathrm{id}_Y) \circ  (e^{\frac{\hbar}{2} t_{X,Z}} \otimes \mathrm{id}_Y) \circ \Psi(\hbar t^{12}_{X,Z,Y}, \hbar t^{23}_{X,Z,Y})^{-1}\circ (\mathrm{id}_X \otimes \tilde{\sigma}_{Y,Z}) \circ (\mathrm{id}_X \otimes e^{\frac{\hbar}{2}t_{Y,Z}}).
      \end{split}
    \end{equation}
    By the naturality of the braiding we have that $t^{23}_{Z,X,Y}
    \circ \sigma_{X \otimes Y,Z} = \sigma_{X \otimes Y,Z} \circ
    t^{12}_{X,Y,Z}$.  Note also, using Equations
    \eqref{eq:hexagon-strict-two} \eqref{eq:cartier-category}, that
    \begin{equation*}
      \begin{split}
	\sigma_{X \otimes Y,Z}^{-1} \circ  t^{12}_{Z,X,Y} \circ \sigma_{X \otimes Y,Z} &= (\mathrm{id}_X \otimes \sigma_{Y,Z}^{-1}) \circ  (\sigma_{X,Z}^{-1} \otimes \mathrm{id}_Y) \circ (t_{Z,X} \otimes \mathrm{id}_Y) \circ (\sigma_{X,Z} \otimes \mathrm{id}_Y) \circ (\mathrm{id}_X \otimes \sigma_{Y,Z})\\
	&= (\mathrm{id}_X \otimes \sigma_{Y,Z}^{-1}) \circ (t_{X,Z} \otimes \mathrm{id}_Y) \circ (\mathrm{id}_X \otimes \sigma_{Y,Z}) = t^{13}_{X,Y,Z}
      \end{split}
    \end{equation*}
    i.e. that $ t^{12}_{Z,X,Y} \circ \sigma_{X \otimes Y,Z} =
    \sigma_{X \otimes Y,Z} \circ t^{13}_{X,Y,Z}$. Hence, using
    Equation \eqref{eq:inverse-associator} we can rewrite the left
    hand side of \eqref{eq:second-hexagon-proof} as
    \begin{equation*}
      \tilde{\sigma}_{X \otimes Y,Z}\circ \Psi(\hbar t^{12}_{X,Y,Z},\hbar t^{13}_{X,Y,Z})\circ e^{\frac{\hbar}{2} (t^{23}_{X,Y,Z}+ t^{13}_{X,Y,Z})} \circ \Psi(\hbar t^{23}_{X,Y,Z},\hbar t^{12}_{X,Y,Z}).
    \end{equation*}
    Next, note that $t^{12}_{X,Z,Y} \circ (\mathrm{id}_X \otimes
    \sigma_{Y,Z}) = (\mathrm{id}_X \otimes \sigma_{Y,Z}) \circ
    t^{13}_{X,Y,Z}$ and, using Equation \eqref{eq:cartier-category},
    we have $t^{23}_{X,Z,Y} \circ (\mathrm{id}_X \otimes \sigma_{Y,Z})
    = (\mathrm{id}_X \otimes \sigma_{Y,Z}) \circ t^{23}_{X,Y,Z}$.
    Therefore, using Equation \eqref{eq:inverse-associator} we can
    rewrite the right hand side of \eqref{eq:second-hexagon-proof} as
    \[ \tilde{\sigma}_{X \otimes Y,Z}\circ e^{\frac{\hbar}{2}t^{13}_{X,Y,Z}} \circ \Psi(\hbar t^{23}_{X,Y,Z},\hbar t^{13}_{X,Y,Z}) \circ e^{\frac{\hbar}{2}t^{23}_{X,Y,Z}}.\]
    Therefore,  Equation \eqref{eq:second-hexagon-proof} is equivalent to the following identity
    \begin{equation}
      \label{eq:second-hexagon-proof-new}
      \Psi(\hbar t^{12}_{X,Y,Z},\hbar t^{13}_{X,Y,Z}) \circ e^{\frac{\hbar}{2}( t^{23}_{X,Y,Z}+ t^{13}_{X,Y,Z})}\circ \Psi(\hbar t^{23}_{X,Y,Z},\hbar t^{12}_{X,Y,Z}) = e^{\frac{\hbar}{2}t^{13}_{X,Y,Z}} \circ \Psi(\hbar t^{23}_{X,Y,Z},\hbar t^{13}_{X,Y,Z}) \circ e^{\frac{\hbar}{2}t^{23}_{X,Y,Z}}.
    \end{equation}
    In order to prove \eqref{eq:second-hexagon-proof-new} it suffices
    to set $A = t^{12}_{X,Y,Z}$, $B= t^{23}_{X,Y,Z}$ and
    $C=t^{13}_{X,Y,Z}$ and use the same argument used to show the
    first hexagon axiom.
  \end{itemize}
\end{proof}
\noindent
To the authors knowledge, the question of whether the deformed
category $\mathcal{C}_{\Psi,t}$ has a non-trivial infinitesimal
braiding is open.
\begin{remark}
  Note that $\mathcal{C}_{\Psi, 0} = \tilde{\mathcal{C}}$.
\end{remark}
\noindent
We exemplify Theorem \ref{thm:main} by deforming the representation
categories of the family of $E(n)$ Cartier Hopf algebras. 
\begin{example}\label{ex:E(n)quant}
Let $n>0$ be a positive integer and consider a field $\Bbbk$ of characteristic zero. Recall from Example \ref{ex:E(n)} the quasitriangular Hopf algebra $(E(n),\mathcal{R})$, with universal $\mathcal{R}$-matrix corresponding to an $n\times n$-matrix $(a_{ij})\in M_n(\Bbbk)$. In Example \ref{ex:E(n)chi} we further discussed that there is an exhaustive family of infinitesimal $\mathcal{R}$-matrices $\chi$ for $(E(n),\mathcal{R})$, determined by $n\times n$-matrices $(b_{k\ell})\in M_n(\Bbbk)$. The latter is Cartier if and only if $(b_{k\ell})\in M_n(\Bbbk)$ is skew-symmetric. We fix an arbitrary $n\times n$-matrix $(a_{ij})\in M_n(\Bbbk)$ and a skew-symmetric $(b_{k\ell})\in M_n(\Bbbk)$.

\noindent
By Theorem \ref{Tannaka-Krein}, the corresponding representation category ${}_{E(n)}\mathcal{M}$ is Cartier, with braiding and infinitesimal braiding determined on left $E(n)$-modules $M,N$ by
\begin{equation}
	\begin{split}
		\sigma_{M,N}(m\otimes n)
		&=\frac{1}{2}\big(1\otimes 1+1\otimes g+g\otimes 1-g\otimes g\big)e^{\sum_{i,j=1}^na_{ij}x_j\otimes gx_i}\cdot(n\otimes m),\\
		t_{M,N}(m\otimes n)
		&=\chi\cdot(m\otimes n)=\sum\nolimits_{k,\ell=1}^nb_{k\ell}gx_k\cdot m\otimes x_\ell\cdot n,
	\end{split}
\end{equation}
for all $m\in M$ and $n\in N$.

\noindent
Choosing a Drinfeld associator $\Psi$ we can apply Theorem \ref{thm:main} to deform the Cartier category $({}_{E(n)}\mathcal{M},\otimes,\sigma,t)$ to obtain a braided monoidal category $({}_{E(n)}\mathcal{M}_{\Psi,t},\tilde{\otimes},a^\Psi,\sigma^\Psi)$. Explicitly, using that the exponents commute, the deformed braiding reads
\begin{equation}\label{E(n)sigmachi}
	\begin{split}
		\sigma^\Psi_{M,N}(m\otimes n)
		&=\tilde{\sigma}_{M,N}\big(e^{\frac{\hbar}{2}t_{M,N}}(m\otimes n)\big)\\
		&=\frac{1}{2}\big(1\otimes 1+1\otimes g+g\otimes 1-g\otimes g\big)e^{a_{ij}x_j\otimes gx_i+\frac{\hbar}{2}b_{k\ell}x_\ell\otimes gx_k}
		\cdot(n\otimes m),
	\end{split}
\end{equation}
while the associativity constraint remains undeformed. The latter is the case, since, by Remark \ref{remark-associator-A-B-commute}, $[t_{12},t_{23}]=0$ implies
\begin{equation}
	\Psi(\hbar t_{M,N} \otimes \mathrm{id}_O,\hbar a^{-1}_{M,N,O} \circ (\mathrm{id}_M \otimes t_{N,O})  \circ a_{M,N,O})=\mathrm{id}_{(M\otimes N)\otimes O}
\end{equation}
for left $E(n)$-modules $M,N,O$.

\noindent
In Example \ref{ex:E(n)chi} we further introduced a gauge
transformation $\mathcal{F}=1\otimes g$ on $E(n)$ and we described the
resulting gauge twisted structure $(E(n)_\mathcal{F},\mathcal{R}_\mathcal{F},\chi_\mathcal{F})$. Theorem \ref{thm:main} can be
applied to the induced Cartier category 
$({}_{E(n)_\mathcal{F}}\mathcal{M},\sigma^{\mathcal{R}_\mathcal{F}},t^{\chi_\mathcal{F}})$. The deformed braiding is
\begin{equation}\label{sigmaE(n)twist}
	\sigma^{\Psi_\mathcal{F}}_{M,N}(m\otimes n)
	=-\frac{1}{2}\big(1\otimes 1-1\otimes g-g\otimes 1-g\otimes g\big)e^{-\sum_{i,j=1}^na_{ij}x_j\otimes gx_i-\frac{\hbar}{2}\sum_{k,\ell=1}^nb_{k\ell}x_\ell\otimes gx_k}\cdot(n\otimes m),
\end{equation}
while the twisted associativity constraint $a^\mathcal{F}$ remains undeformed for the same reason as before.
\end{example}
\noindent
In the previous example we notice that \eqref{sigmaE(n)twist} is the gauge transformation of \eqref{E(n)sigmachi}, or, in other words, deformation in the sense of Theorem \ref{thm:main} and gauge transformation commute. As one easily verifies, this is a general feature of the theory, meaning that Theorem \ref{thm:main} is compatible with equivalences in the form of gauge transformations.

\subsection{Deformation of pre-Cartier categories}\label{sec:preCartierCase}

In Example \ref{ex:E(n)quant} it
turns out that the associativity constraint is undeformed under Theorem
\ref{thm:main}. We show that this is the case, since the commutator
\begin{equation}\label{chicom}
	[\chi_{12},\chi_{23}]=0
\end{equation}
vanishes. Even more, it turns out that condition \ref{chicom} is
sufficient to imply a quantization theorem, even without imposing the
Cartier condition \eqref{eq:cartier-category}. This is the content of
this section and our second main deformation result.

\begin{theorem}
  \label{theorem-quantization-of-pre-cartier-categories}
  Let $(\mathcal{C},\sigma,t)$ be a $\Bbbk$-linear pre-Cartier category satisfying 
  \begin{equation}
    \label{eq:commutation-12-23}
	  [t^{12}_{X,Y,Z}, t^{23}_{X,Y,Z}]=0
  \end{equation}
  for any objects $X,Y,Z$. Then $(\mathcal{C},\Bbbk, \tilde{\otimes} ,
  \tilde{a} , \tilde{\ell}, \tilde{r}, \hat{\sigma})$ is a braided
  monoidal category, where $\hat{\sigma}_{X,Y} = \tilde{\sigma}_{X,Y}
  \circ e^{\frac{\hbar}{2} t_{X,Y}}.$
\end{theorem}
\begin{proof}
  First of all, note that Equation \eqref{eq:commutation-12-23}
  together with the fact that the morphisms $t^{ij}_{X,Y,Z}$ satisfy
  the infinitesimal braid relations \eqref{eq:infbraid1} (see
  Proposition \ref{proposition-inf-braid-rel-in-pre-cartier-category})
  imply
  \begin{align}
    [t^{12}_{X,Y,Z}, t^{13}_{X,Y,Z}]&=0 \label{eq:commutation-12-13}\\
    [t^{23}_{X,Y,Z}, t^{13}_{X,Y,Z}]&=0.\label{eq:commutation-23-13}
  \end{align}
  The triangle and the pentagon axioms hold since $\tilde{\mathcal{C}}$ is a monoidal
  category (see Example \ref{example-trivial-deformation}). Next, we
  show the hexagons axioms \eqref{eq:hexagon-strict-one}
  \eqref{eq:hexagon-strict-two}. We have, using Equations
  \eqref{eq:pre-cartier-one} \eqref{eq:commutation-12-13}:
  \begin{equation*}
    \begin{split}
     \hat{\sigma}_{X, Y \otimes Z} &= (\tilde{\sigma}_{X, Y \otimes Z}) \circ (e^{\frac{\hbar}{2} t_{X,Y \otimes Z}}) \\
      &= (\tilde{\sigma}_{X, Y \otimes Z}) \circ ( e^{\frac{\hbar}{2} (t^{12}_{X,Y,Z} + t^{13}_{X,Y,Z})})\\
      &=(\mathrm{id}_Y \otimes \tilde{\sigma}_{X,Z})\circ(\tilde{\sigma}_{X,Y} \otimes \mathrm{id}_Z) \circ  ( e^{\frac{\hbar}{2} (t^{13}_{X,Y,Z})}) \circ ( e^{\frac{\hbar}{2} (t^{12}_{X,Y,Z})}) \\
      &= (\mathrm{id}_Y \otimes \tilde{\sigma}_{X,Z})\circ ( e^{\frac{\hbar}{2} (t^{23}_{Y,X,Z})})  \circ (\tilde{\sigma}_{X,Y} \otimes \mathrm{id}_Z) \circ ( e^{\frac{\hbar}{2} (t^{12}_{X,Y,Z})}) \\
      &=(\mathrm{id}_Y \otimes \hat{\sigma}_{X,Z})\circ(\hat{\sigma}_{X,Y} \otimes \mathrm{id}_Z)
    \end{split}
  \end{equation*}
  where the fourth equality follows from the fact that $t^{23}_{Y,X,Z}
  \circ (\sigma_{X,Y} \otimes \mathrm{id}_Z) = (\sigma_{X,Y} \otimes
  \mathrm{id}_Z) \circ t^{13}_{X,Y,Z}$. Similarly, using Equations
  \eqref{eq:pre-cartier-two} \eqref{eq:commutation-23-13}, we have:
  \begin{equation*}
    \begin{split}
      \hat{\sigma}_{X \otimes Y,Z} &= (\tilde{\sigma}_{X \otimes Y, Z}) \circ (e^{\frac{\hbar}{2} t_{X \otimes Y,Z}}) \\
      &= (\tilde{\sigma}_{X \otimes Y, Z}) \circ (e^{\frac{\hbar}{2} (t^{23}_{X,Y,Z}+ t^{13}_{X,Y,Z})} ) \\
      &= (\tilde{\sigma}_{X,Z} \otimes \mathrm{id}_Y) \circ (\mathrm{id}_X \otimes \tilde{\sigma}_{Y,Z}) \circ (e^{\frac{\hbar}{2} t^{13}_{X,Y,Z}} ) \circ (e^{\frac{\hbar}{2} t^{23}_{X,Y,Z}}) \\
      &=  (\tilde{\sigma}_{X,Z} \otimes \mathrm{id}_Y) \circ (e^{\frac{\hbar}{2} t^{12}_{X,Z,Y}} ) \circ (\mathrm{id}_X \otimes \tilde{\sigma}_{Y,Z}) \circ (e^{\frac{\hbar}{2} t^{23}_{X,Y,Z}}) \\
      &= (\hat{\sigma}_{X,Z} \otimes \mathrm{id}_Y) (\mathrm{id}_X \otimes \hat{\sigma}_{Y,Z}),
    \end{split}
  \end{equation*}
  where the fourth equality follows from $t^{12}_{X,Z,Y}
  \circ (\mathrm{id}_X \otimes \sigma_{Y,Z}) = (\mathrm{id}_X \otimes
  \sigma_{Y,Z}) \circ t^{13}_{X,Y,Z}$.
\end{proof}
\begin{example}
  We can deform the representation categories of
  $(E(n),\mathcal{R}_{(a_{ij})},\chi_{(b_{k\ell})})$ discussed in
  Example \ref{ex:E(n)quant} for
  arbitrary $n\times n$-matrices, i.e. even without assuming that
  $(b_{k\ell})$ is skew-symmetric.
\end{example}
\noindent
The following remark shows that we do not expect the deformed category arising from Theorem \ref{theorem-quantization-of-pre-cartier-categories} to be symmetric in general. Even if the initial category is symmetric, such an assumption turns out to be rather restrictive.
\begin{remark}
One verifies that the deformed category constructed by Theorem \ref{theorem-quantization-of-pre-cartier-categories} is symmetric if and only if the starting category is symmetric and
\begin{equation}\label{minusCartier}
t_{Y,X} \circ \sigma_{X,Y} = -\sigma_{X,Y} \circ t_{X,Y}
\end{equation}
holds for all pair of objects $X,Y$ in $\mathcal{C}$. In particular, if we start with a Cartier symmetric category, \eqref{minusCartier} forces $t=0$. The condition \eqref{minusCartier} has been noted in \cite[Remark 2.8]{ABSW} in the case of the representation category of Sweedler's Hopf algebra.
\end{remark}

\subsection{Quantization of infinitesimal $\mathcal{R}$-matrices}\label{sec:4.4}

This final section contains the main application of Theorems \ref{thm:main} and 
\ref{theorem-quantization-of-pre-cartier-categories}, namely
sufficient conditions for quantizing the infinitesimal $\mathcal{R}$-matrix of a
pre-Cartier quasi-bialgebra.
\begin{theorem}
Let $(H, \Phi, \mathcal{R}, \chi)$ be a Cartier quasi-bialgebra and let $\Psi$ be a Drinfeld associator. Then there is a topological quasitriangular quasi-bialgebra structure on $\tilde{H} = H[[\hbar]]$ such that 
\[\tilde{\mathcal{R}} = \mathcal{R} \big(1 \ten 1 + \hbar \chi + \mathcal{O}(\hbar^2)\big).\]
Explicitly, its re-associator and universal $\mathcal{R}$-matrix are
\begin{equation}
\label{eq:quantization-theorem-formulae}
\tilde{\Phi} = \Psi(\hbar \chi_{12},\hbar \Phi^{-1}\chi_{23}\Phi) \qquad \text{and} \qquad \tilde{\mathcal{R}} = \mathcal{R} e^{\frac{\hbar}{2}\chi}.
\end{equation}
\end{theorem}
\begin{proof}
By Proposition \ref{prop:quasi-bialgCAT} (see also Remark \ref{rem:topfree}) and Theorem \ref{Tannaka-Krein} we have that the category ${}_H\mathcal{M}$ is Cartier, where the associativity constraint, braiding, and infinitesimal braiding $a, \sigma,t$ are induced respectively by $\Phi, \mathcal{R},\chi$. Applying Theorem \ref{thm:main} we obtain a deformed braided monoidal category  ${}_H\mathcal{M}_{\Psi,t}$,  with deformed associativity constraint and braiding as in \eqref{deformedstructure}. By construction, the category ${}_H\mathcal{M}_{\Psi,t}$ has the same objects of ${}_H\mathcal{M}$, hence we can apply Proposition \ref{prop:quasi-bialgCAT} to obtain a topological quasitriangular quasi-bialgebra, which readily satisfies \eqref{eq:quantization-theorem-formulae}.
\end{proof}

\noindent
It is possible to prove the algebraic analogue of Theorem \ref{theorem-quantization-of-pre-cartier-categories} using the same strategy as above. However, it turns out that here we can even give an explicit proof, where we directly verify the axioms of a quasitriangular quasi-bialgebra, without the need to pass to the categorical picture and perform the deformation there. For this reason we spell out the direct proof below.

\noindent
Let $(H,\mathcal{R},\chi)$
be a pre-Cartier bialgebra. It turns out that one of
the equivalent conditions (see Theorem
\ref{theorem-inf-braid-rel-in-qtpcb})
\begin{align}
	[\chi_{12},\mathcal{R}^{-1}_{12}\chi_{13}\mathcal{R}_{12}]
	&=0\label{com1}\\
	[\chi_{23},\mathcal{R}^{-1}_{23}\chi_{13}\mathcal{R}_{23}]
	&=0\label{com2} \\
	[\chi_{12},\chi_{23}]&=0 \label{com3}
\end{align}
is sufficient for $\tilde{\mathcal{R}}:=\mathcal{R}\exp(\hbar\chi)$
to be a universal $\mathcal{R}$-matrix for the trivial topological
bialgebra $\tilde{H}:=H[[\hbar]]$.
\begin{theorem}\label{thm:exp}
  For a pre-Cartier bialgebra $(H,\mathcal{R},\chi)$
  satisfying one of the equivalent conditions
  \eqref{com1}-\eqref{com3} there is a quasitriangular structure
  \begin{equation}
    \tilde{\mathcal{R}}:=\mathcal{R}\exp(\hbar\chi)
  \end{equation}
  on the trivial topological bialgebra $\tilde{H}:=H[[\hbar]]$.
\end{theorem}
\begin{proof}
  Using Equations \eqref{eq:qtqb1} and \eqref{eq:pC1} it follows that
  $$
  \tilde{\mathcal{R}}\Delta(\cdot)
  =\mathcal{R}\sum_{n=0}^\infty\frac{\hbar^n}{n!}\chi^n\Delta(\cdot)
  =\mathcal{R}\Delta(\cdot)\sum_{n=0}^\infty\frac{\hbar^n}{n!}\chi^n
  =\Delta^\mathrm{op}(\cdot)\mathcal{R}\sum_{n=0}^\infty\frac{\hbar^n}{n!}\chi^n
  =\Delta^\mathrm{op}(\cdot)\tilde{\mathcal{R}}.
  $$
  Further assuming \eqref{com1} shows that
  \begin{align*}
    (\mathrm{id}\otimes\Delta)(\tilde{\mathcal{R}})
    &=(\mathrm{id}\otimes\Delta)(\mathcal{R})\sum_{n=0}^\infty\frac{\hbar^n}{n!}(\mathrm{id}\otimes\Delta)(\chi^n)\\
    &\overset{\eqref{eq:qtqb2}}{=}\mathcal{R}_{13}\mathcal{R}_{12}\sum_{n=0}^\infty\frac{\hbar^n}{n!}((\mathrm{id}\otimes\Delta)(\chi))^n\\
    &\overset{\eqref{eq:pC2}}{=}\mathcal{R}_{13}\mathcal{R}_{12}\sum_{n=0}^\infty\frac{\hbar^n}{n!}(\chi_{12}+\mathcal{R}^{-1}_{12}\chi_{13}\mathcal{R}_{12})^n\\
    &\overset{\eqref{com1}}{=}\mathcal{R}_{13}\mathcal{R}_{12}\sum_{n=0}^\infty\frac{\hbar^n}{n!}\sum_{k=0}^n\binom{n}{k}(\mathcal{R}^{-1}_{12}\chi_{13}\mathcal{R}_{12})^k\chi_{12}^{n-k}\\
    &=\mathcal{R}_{13}\mathcal{R}_{12}\sum_{n=0}^\infty\frac{\hbar^n}{n!}\sum_{k=0}^n\binom{n}{k}\mathcal{R}^{-1}_{12}\chi_{13}^k\mathcal{R}_{12}\chi_{12}^{n-k}\\
    &=\mathcal{R}_{13}\sum_{n=0}^\infty\hbar^n\sum_{k=0}^n\frac{1}{k!(n-k)!}\chi_{13}^k\mathcal{R}_{12}\chi_{12}^{n-k},
  \end{align*}
  and
  \begin{align*}
    \tilde{\mathcal{R}}_{13}\tilde{\mathcal{R}}_{12}
    &=\mathcal{R}_{13}\sum_{k=0}^\infty\frac{\hbar^k}{k!}\chi_{13}^k\mathcal{R}_{12}\sum_{r=0}^\infty\frac{\hbar^r}{r!}\chi_{12}^r\\
    &=\mathcal{R}_{13}\sum_{k,r=0}^\infty\frac{\hbar^{k+r}}{k!r!}\chi_{13}^k\mathcal{R}_{12}\chi_{12}^r\\
    &=\mathcal{R}_{13}\sum_{n=0}^\infty\hbar^n\sum_{k=0}^n\frac{1}{k!(n-k)!}\chi_{13}^k\mathcal{R}_{12}\chi_{12}^{n-k}
  \end{align*}
  coincide, where in the last equation we introduced the new summation
  index $n:=k+r$ so that $r=n-k$ and the summation over $k$ goes from
  $0$ to $n$. In complete analogy one proves
  $(\Delta\otimes\mathrm{id})(\tilde{\mathcal{R}})=\tilde{\mathcal{R}}_{13}\tilde{\mathcal{R}}_{23}$
  using Equations \eqref{eq:qtqb3}, \eqref{eq:pC3}, and \eqref{com2}.
\end{proof}
\noindent
In the previous proof we have seen that \eqref{com1} ensures
$(\mathrm{id}\otimes\Delta)(\tilde{\mathcal{R}})=\tilde{\mathcal{R}}_{13}\tilde{\mathcal{R}}_{12}$,
while \eqref{com2} implies
$(\Delta\otimes\mathrm{id})(\tilde{\mathcal{R}})=\tilde{\mathcal{R}}_{13}\tilde{\mathcal{R}}_{23}$. Note
that
$\tilde{\mathcal{R}}\Delta(\cdot)=\Delta^\mathrm{op}(\cdot)\tilde{\mathcal{R}}$
holds automatically, not even assuming \eqref{com1} or \eqref{com2}.

\noindent
Another observation is that in case $(H,\mathcal{R},\chi)$ is
Cartier, the
conditions \eqref{com1} \eqref{com2} are equivalent to
\begin{align}
	[\chi_{12},\chi_{23}]
	&=0\label{com4}\\
	[\chi_{12},\chi_{13}]
	&=0\label{com5} 
\end{align}
respectively. To see this, multiply
$0=[\chi_{12},\mathcal{R}^{-1}_{12}\chi_{13}\mathcal{R}_{12}]$ by
$\mathcal{R}_{12}$ from the left and by $\mathcal{R}^{-1}_{12}$ from
the right. Then
$$
0=\mathcal{R}_{12}(\chi_{12}\mathcal{R}^{-1}_{12}\chi_{13}\mathcal{R}_{12}-\mathcal{R}^{-1}_{12}\chi_{13}\mathcal{R}_{12}\chi_{12})\mathcal{R}^{-1}_{12}
=\chi_{21}\chi_{13}-\chi_{13}\chi_{21},
$$ which gives \eqref{com4} after tensor flip. Similarly, \eqref{com2}
implies \eqref{com5} in the Cartier case.\\
\begin{remark}
  Theorem \ref{thm:exp} is a refinement of
  \cite[Proposition 4.1]{BRS24}, where in addition nilpotence ($\chi^n=0$
  for some $n>1$) is assumed (eventhough this property is not used in
  the proof), and both identities \eqref{com1} and \eqref{com2} are
  assumed, not mentioning that they are equivalent.
\end{remark}

\noindent
For instance, we can obtain deformations of the pre-Cartier bialgebras $E(n)$ in Example \ref{ex:E(n)quant} without employing any Drinfeld associator. Namely, 
$$
\tilde{\mathcal{R}}=\tau(\hat{\sigma}_{H,H}(1\otimes 1))
=\frac{1}{2}\big(1\otimes 1+1\otimes g+g\otimes 1-g\otimes g\big)e^{a_{ij}gx_i\otimes x_j+\frac{\hbar}{2}b_{k\ell}gx_k\otimes x_\ell}
$$ 
is a quasitriangular structure for the topological bialgebra $E(n)[[\hbar]]$
and
$$
\tilde{\mathcal{R}}_\mathcal{F}
=\tau(\hat{\sigma}^{\mathcal{F}}_{H,H}(1\otimes 1))
=-\frac{1}{2}\big(1\otimes 1-1\otimes g-g\otimes 1-g\otimes g\big)e^{-\sum_{i,j=1}^na_{ij}gx_i\otimes x_j-\frac{\hbar}{2}\sum_{k,\ell=1}^nb_{k\ell}gx_k\otimes x_\ell}
$$ 
is a quasitriangular structure for the topological quasi-bialgebra $(E(n)[[\hbar]],\tilde{\Delta}_\mathcal{F},\varepsilon)$ with re-associator $1\otimes 1\otimes g$.
The previous two $\mathcal{R}$-matrices are by construction related via the gauge transformation $\mathcal{F}=1\otimes g$.

~\\
~\\
\textbf{Contacts:}\\
~\\
chesposito@unisa.it\\
andrea.rivezzi@matfyz.cuni.cz\\
jonaschristoph.schnitzer@unipv.it\\
thomas.weber@matfyz.cuni.cz

\begin{thebibliography}{9}

\bibitem{ABSW}
\textsc{A. Ardizzoni, L. Bottegoni, A. Sciandra, T. Weber:}
\emph{Infinitesimal braidings and pre-Cartier bialgebras}.
Commun. Contemp. Math. \textbf{27}, 5 (2025) 2450029.

\bibitem{PaoloWess}
\textsc{P. Aschieri, M. Dimitrijevic, F. Meyer, J. Wess:}
\emph{Noncommutative geometry and gravity}.
Class. Quantum Gravity \textbf{23} (2006) 1883-1911.

\bibitem{Baxter}
\textsc{R.J. Baxter:}
\emph{Solvable eight-vertex model on an arbitrary planar lattice}.
Philos. Trans. Royal Society A \textbf{289}, 1359 (1978).

\bibitem{BeaDasGru}
\textsc{M. Beattie, S. Dăscălescu, L. Grünenfelder:}
\emph{Constructing Pointed Hopf Algebras by Ore Extensions}.
J. Alg. \textbf{225}, 2 (2000) 743-770.

\bibitem{BRW} 
\textsc{M. Bordemann, A. Rivezzi, T. Weigel:}
\emph{A gentle introduction to Drinfel'd associators}. 
Rev. Math. Phys. \textbf{37}, 04 (2025) 2430010.  

\bibitem{BRS24}
\textsc{L. Bottegoni, F. Renda, A. Sciandra:} 
\emph{Infinitesimal R-matrices for some families of Hopf algebras}. Preprint \href{https://arxiv.org/abs/2412.02350}{arXiv:2412.02350}. 

\bibitem{Bul}
\textsc{D. Bulacu,  S. Caenepeel, F. Panaite, F. Van Oystaeyen:}
\emph{Quasi-Hopf algebras: A categorical approach}.
Cambridge University Press (2019).

\bibitem{CaenDasc}
\textsc{S. Caenepeel, S. Dăscălescu:}
\emph{On pointed Hopf algebras of dimension $2n$.} 
Bull. London Math. Soc. \textbf{31} (1999) 17-24.

\bibitem{Car}
\textsc{P. Cartier:}
\emph{Construction combinatoire des invariants de Vassiliev-Kontsevich des noeuds}.
Recherche Coop\'erative sur Programme \textbf{25}, 45 (1993) 1-10.

\bibitem{CirMart}
\textsc{L. S. Cirio, J. F. Martins:}
\emph{Infinitesimal 2-braidings and differential crossed modules}. Adv.  Math.  \textbf{277} (2015) 426-491.

\bibitem{Doikou22}
\textsc{A. Doikou, A. Ghionis, B. Vlaar:}
\emph{Quasi-bialgebras from set-theoretic type solutions of the Yang–Baxter equation}.
Lett. Math. Phys. \textbf{112}, 78 (2022).

\bibitem{DrHopfAlg}
\textsc{V.G. Drinfeld:}
\emph{Hopf algebras and the quantum Yang-Baxter equation}.
Dokl. Akad. Nauk SSSR \textbf{283}, 5 (1985) 1060-1064.

\bibitem{DR89} 
\textsc{V.G. Drinfeld:}
\emph{On quasitriangular quasi-Hopf algebras and a group closely connected to Gal($\overline{\mathbb{Q}}/\mathbb{Q}$)}.
Algebra Anal. \textbf{4}, 2 (1990) 149-181.

\bibitem{DrQG} 
\textsc{V.G. Drinfeld:} 
\emph{Quantum Groups}. 
In \textsc{A. Gleason (ed.):} \emph{Proceedings of the ICM}, AMS, Rhode Island (1987) 798-820.

\bibitem{DrQH} 
\textsc{V.G. Drinfeld:} 
\emph{Quasi Hopf algebras}. 
Leningr. Math. J. \textbf{6}, 1 (1990) 1419-1457.

\bibitem{DrKZ}
\textsc{V.G. Drinfeld:} 
\emph{Quasi-Hopf algebras and Knizhnik-Zamolodchikov equations}. In: Belavin, A.A., Klimyk, A.U., Zamolodchikov, A.B. (eds) \emph{Problems of Modern Quantum Field Theory}. Research Reports in Physics. Springer (1989) 1-13.

\bibitem{ChJoSt}
\textsc{C. Esposito, J. Schnitzer, S. Waldmann:}
\emph{A universal construction of universal deformation formulas, Drinfeld twists and their positivity}.
Pac. J. Math. \textbf{291}, 2 (2017) 319–358.

\bibitem{EGNO}
\textsc{P. I. Etingof, S. Gelaki, D. Nikshych, V. Ostrik:} 
\emph{Tensor categories}. 
Math. Surv. Monogr. \textbf{205}, Providence, RI: American Mathematical Society (AMS), 2015.

\bibitem{EtiKaz}
\textsc{P.I. Etingof, D.A. Kazhdan:}
\emph{Quantization of Lie bialgebras}. 
I. Sel. Math. (N.S.) \textbf{2} (1996) 1–41.

\bibitem{EtiKazII}
\textsc{P.I Etingof, D.A. Kazhdan:}
\emph{Quantization of Lie bialgebras}. 
II. Sel. Math. (N.S.) \textbf{4} (1998) 213–231.

\bibitem{FaTa}
\textsc{L.D. Faddeev, L.A. Takhtadzhan:}
\emph{The Quantum Method of the Inverse Problem and the Heisenberg XYZ Model}.
Russ. Math. Surv. \textbf{34}, 5 (1979) 11–68.

\bibitem{FGS}
\textsc{M. Faitg, A.M. Gainutdinov, C. Schweigert:}
\emph{An adjunction theorem for Davydov-Yetter cohomology and infinitesimal braidings}.
Preprint arXiv:2411.19111.

\bibitem{HeckVen}
\textsc{I. Heckenberger, L. Vendramin:}
\emph{Bosonization of curved Lie bialgebras}. 
Bull. Belg. Math. Soc. Simon Stevin \textbf{30}, 5 (2023) 577-600.

\bibitem{Hum}
\textsc{J. E. Humphreys:}
\emph{Introduction to Lie algebras and representation theory}.
\textbf{9}, Springer Science \& Business Media, 2012.

\bibitem{Jimbo}
\textsc{M. Jimbo:}
\emph{A q-difference analogue of U(g) and the Yang-Baxter
equation}. 
Lett. Math. Phys. \textbf{10} (1985) 63-69.

\bibitem{JS93}
\textsc{A. Joyal, R. Street:} 
\emph{Braided tensor categories}. 
Adv. Math. 	\textbf{102} (1993) 20-78.

\bibitem{KKMP}
\textsc{E. Karlsson, C. Keller, L. Müller, J. Pulmann:}
\emph{Deformation Quantization via Categorical Factorization Homology}.
Preprint arXiv:2410.12516.

\bibitem{Kassel95}
\textsc{C. Kassel:} \emph{Quantum groups}. Graduate Texts in Mathematics, \textbf{155}. Springer-Verlag, New York, 1995.

\bibitem{KeLaSc}
\textsc{C. Kemp, R. Laugwitz, A. Schenkel:}
\emph{Infinitesimal 2-braidings from 2-shifted Poisson structures}.
J. Geom. Phys. \textbf{212} (2025) 105456.

\bibitem{Kemp}
\textsc{C. Kemp:}
\emph{Syllepses from 3-shifted Poisson structures and second-order integration of infinitesimal 2-braidings}.
Preprint arXiv:2505.01949.

\bibitem{Kohno}
\textsc{T. Kohno:}
\emph{Monodromy representations of braid groups and Yang-Baxter equations}. 
Ann. Inst. Fourier (Grenoble) \textbf{37}, 4 (1987) 139-160.

\bibitem{MacLan63} 
\textsc{S. Mac Lane:} 
\emph{Natural Associativity and Commutativity}.  
Rice Univ. Stud. \textbf{49} (1963) 28-46.

\bibitem{MajidFoundation}
\textsc{S. Majid:} 
\emph{Foundations of quantum group theory}.
Cambridge University Press, Cambridge, 1995.

\bibitem{Majid98} 
\textsc{S. Majid:} 
\emph{Quantum Double for Quasi-Hopf Algebras}.  
Lett. Math. Phys. \textbf{45} (1998) 1-9.

\bibitem{PulSev}
\textsc{J. Pulmann, P. \v{S}evera:}
\emph{Quantization of Poisson Hopf algebras}.
Adv. Math. \textbf{401} (2022) 108310.

\bibitem{ReTu}
\textsc{N. Reshetikhin, V.G. Turaev:}
\emph{Invariants of 3-manifolds via link polynomials and quantum groups}.
Invent. Math. \textbf{103} (1991) 547-597.

\bibitem{AndreaThesis}
\textsc{A. Rivezzi:}
\emph{Universal constructions arising from quantization of Lie bialgebras}.
Ph.D. Thesis \href{https://www.boa.unimib.it/handle/10281/475779}{www.boa.unimib.it/handle/10281/475779}, 2024.

\bibitem{SaaRiv}
\textsc{N. Saavedra-Rivano:}
\emph{Catégories Tannakiennes}.
Bull. Soc. Math. Fr. \textbf{100} (1972) 417-430.

\bibitem{Sch01}
\textsc{P. Schauenburg:}
 \emph{Turning monoidal categories into strict ones}.
 New York J. Math \textbf{7}, 102 (2001) 257-265.
 
\bibitem{Sev}
\textsc{P. \v{S}evera:} 
\emph{Quantization of Lie bialgebras revisited}.
Sel. Math., New Ser. \textbf{22}, 3 (2016) 1563–1581.

\bibitem{Thomas}
\textsc{T. Weber:}
\emph{Braided Cartan calculi and submanifold algebras}. 
J. Geom. Phys. \textbf{150} (2020) 103612.

\bibitem{Yett}
\textsc{D.  N. Yetter:}
\emph{Framed Tangles and a Theorem of Deligne on Braided Deformations of Tannakian Categories.}
Deformation Theory and Quantum Groups with Applications to Mathematical Physics (M. Gerstenhaber and J.D. Stasheff, eds.) AMS Contemp. Math. vol. 134 325-350 (1992).

\bibitem{Yet}
\textsc{D.  N. Yetter:}
\emph{Braided deformations of monoidal categories and Vassiliev invariants}.
M. Kapranov, E. Getzler
(Eds.), Higher Category Theory, AMS Contemporary Mathematics, \textbf{320}, AMS, Providence, RI, 1998,
117-134.
\end{thebibliography}
\end{document}